\documentclass[10pt]{article}
\usepackage{graphicx,float,latexsym,amsmath,amsfonts,amssymb,subeqnarray,colordvi,epstopdf,bbm}
\usepackage[latin1]{inputenc}
\usepackage{a4wide}
\usepackage[hidelinks]{hyperref}
\usepackage{authblk}

\def\R{\mathbb{R}}

\def\imi{\textbf{\hskip1pt i\hskip1pt}}

\title{A Finite Volume - Alternating Direction Implicit Approach for the Calibration of Stochastic Local Volatility Models}
\author[1]{Maarten Wyns\footnote{Corresponding author.}}
\author[2]{Jacques Du Toit}
\affil[1]{Department of Mathematics and Computer Science,
University of Antwerp, Middelheimlaan 1, B-2020 Antwerp, Belgium.
\mbox{Email}: \texttt{maarten.wyns@uantwerpen.be};}
\affil[2]{Numerical Algorithms Group, Peter House, Oxford street, Manchester, M1 5AN, United Kingdom. \mbox{Email}:  \texttt{jacques@nag.co.uk}.}

\date{\today}

\begin{document}

\maketitle

\begin{abstract}
\noindent
Calibration of stochastic local volatility (SLV) models to their underlying local volatility model is often performed by numerically solving a two-dimensional non-linear forward Kolmogorov equation.
We propose a novel finite volume (FV) discretization in the numerical solution of general 1D and 2D forward Kolmogorov equations. 
The FV method does not require a transformation of the PDE. This constitutes a main advantage in the calibration of SLV models as the pertinent PDE coefficients are often nonsmooth. Moreover, the FV discretization has the crucial property that the total numerical mass is conserved. 
Applying the FV discretization in the calibration of SLV models yields a non-linear system of ODEs.
Numerical time stepping is performed by the Hundsdorfer--Verwer ADI scheme to increase the computational efficiency.
The non-linearity in the system of ODEs is handled by introducing an inner iteration.
Ample numerical experiments are presented that illustrate the effectiveness of the calibration procedure.
\end{abstract}
\vspace{0.2cm}\noindent

\vspace{3mm}
\normalsize

\setcounter{equation}{0}
\section{Introduction}\label{intro}

In contemporary financial mathematics, \textit{stochastic local volatility} (SLV) models constitute state-of-the-art models to describe asset price processes, notably foreign exchange (FX) rates, see e.g.\ \cite{L02,TF10}. 
Let $S_{\tau}$ represent the exchange rate at time $\tau \ge 0$ and let the spot value $S_{0}$ be given.
For modelling the exchange rate we consider the transformation $X_{\tau} = \log(S_{\tau}/S_{0})$ since the transformed variable $X_{\tau}$ reflects better the duality between the exchange rate $S_{\tau}$ and $1/S_{\tau}$, where the latter one is the exchange rate when the role of the domestic currency and the foreign currency are swapped. In this article we consider SLV models of the type
\begin{equation}
\left\{ \begin{array}{l}
dX_{\tau} = (r_{d} - r_{f} - \tfrac{1}{2}\sigma^{2}_{SLV} (X_{\tau},\tau) \psi^{2}(V_{\tau})) d\tau + \sigma_{SLV} (X_{\tau},\tau) \psi(V_{\tau}) dW^{(1)}_{\tau}, \\\\
dV_{\tau} = \kappa (\eta - V_{\tau}) d\tau + \xi V_{\tau}^{\alpha} dW^{(2)}_{\tau},
\end{array} \right.
\label{eq:SLVmodel}
\end{equation}
with $\psi$ a non-negative function on $\R^{+}$ such that $\psi(0)=0$, $\alpha, \kappa, \eta, \xi$ strictly positive parameters, $dW^{(1)}_{\tau} \cdot dW^{(2)}_{\tau} = \rho d\tau$, $-1 \leq \rho \leq 1$ and given spot values $X_{0}=0, V_{0}$.
The non-negative function $\sigma_{SLV}(x,\tau)$ is called the \textit{leverage function} and the constant $r_{d}$, respectively $r_{f}$, denotes the risk-free interest rate in the domestic currency, respectively in the foreign currency.
It is readily seen that the process $V_{\tau}$ is always non-negative. 
In \cite{AP06} it is shown that the boundary $V_{\tau}=0$ is attainable for $0<\alpha<1/2$ and for $\alpha=1/2$ if $2\kappa \eta < \xi^2$. For $\alpha>1/2$ it holds that $V_{\tau}=0$ is an unattainable boundary. Furthermore, $V_{\tau} = \infty$ is an unattainable boundary for all values of $\alpha>0$.
The choice $\psi(v)=\sqrt{v}, \alpha = 1/2$ corresponds to the Heston-based SLV model and the choice $\psi(v)=v, \alpha = 1$ corresponds to the SLV model described in \cite{TF10}.

SLV models constitute a natural combination of \textit{local volatility} (LV) models and \textit{stochastic volatility} (SV) models, \cite{L02,TF10}. The former type of models can be described by the stochastic differential equation (SDE)
\begin{equation}
dX_{LV,\tau} = (r_{d} - r_{f} - \tfrac{1}{2} \sigma_{LV}^{2}(X_{LV,\tau},\tau) ) d\tau + \sigma_{LV}(X_{LV,\tau},\tau) dW_{\tau}.
\label{eq:LVmodel}
\end{equation}
The function $\sigma_{LV}(x,\tau)$ is called the \textit{local volatility function} and can be determined by the Dupire formula \cite{D94} such that the LV model reproduces the known market prices for European call and put options. 
Since the LV model is completely determined by the market prices of European call and put options, it offers no flexibility in matching the market dynamics.
The underlying pure SV model corresponding with \eqref{eq:SLVmodel} is given by the system of SDEs
\begin{equation}
\left\{ \begin{array}{l}
dX_{SV,\tau} = (r_{d} - r_{f} - \tfrac{1}{2}\psi^{2}(V_{SV,\tau})) d\tau + \psi(V_{SV,\tau}) dW^{(1)}_{\tau}, \\\\
dV_{SV,\tau} = \kappa (\eta - V_{SV,\tau}) d\tau + \xi V_{SV,\tau}^{\alpha} dW^{(2)}_{\tau},
\end{array} \right.
\label{eq:SVmodel}
\end{equation}
where the function $\psi$ and the parameters need to satisfy the same restrictions as above. 
SV models are typically well-suited to reflect forward volatilities, but they are often unable to capture the volatility smile exactly, see e.g.\ \cite{VGO14}.
By combining features of the LV model with features of the SV models, SLV models are able to match the market dynamics and to reproduce the market prices for European call and put options. If the leverage function is identically equal to one, then the SLV model reduces to the SV model \eqref{eq:SVmodel}. If the parameter $\xi$ equals zero, then the SLV model reduces to a purely LV model.

In financial practice, see e.g.\ \cite{C11,TF10}, one will first determine the SV parameters $\kappa,\eta,\xi,\rho$ such that the underlying SV model matches the current market dynamics. Afterwards, the leverage function $\sigma_{SLV} (x,\tau)$ is calibrated such that the SLV model and the underlying LV model define the same fair value for non-path-dependent European options. Then, since LV models are calibrated exactly to the known market prices for vanilla options, the SLV model also reproduces the known market prices for European call and put options.
Let $p(x,v,\tau;X_{0},V_{0})$ denote the joint density of $(X_{\tau},V_{\tau})$ under the SLV model \eqref{eq:SLVmodel} and let $p_{LV}(x,\tau;X_{0})$ denote the density of $X_{LV,\tau}$ under \eqref{eq:LVmodel}. 
It is well-known, see e.g.\ \cite{G86,T11}, that both the SLV model \eqref{eq:SLVmodel} and the LV model \eqref{eq:LVmodel} define the same marginal distribution for $S_{\tau}$, i.e.\
\begin{equation}
p_{LV}(x,\tau;X_{0}) = \int_{0}^{\infty} p(x,v,\tau;X_{0},V_{0}) dv,
\label{eq:SLVmatchesLV}
\end{equation} 
for $\tau > 0$, if
\begin{equation}
\sigma^{2}_{LV} (x,\tau) = \mathbb{E}[\sigma_{SLV}^{2} (X_{\tau},\tau)\psi^{2}(V_{\tau}) \vert X_{\tau} = x] = \sigma_{SLV}^{2} (x,\tau)\mathbb{E}[\psi^{2}(V_{\tau}) \vert X_{\tau} = x].
\label{eq:SigmatoMatch}
\end{equation}
If one can determine the conditional expectation above, and if the leverage function is defined by \eqref{eq:SigmatoMatch}, then both models yield the same fair value for non-path-dependent European options and hence the SLV model is calibrated exactly to the market prices for vanilla options.
This is, however, highly non-trivial since the conditional expectation itself depends on the leverage function.

In the past years, a variety of numerical techniques, see e.g.\ \cite{C11,EKO12,H09,RMQ07,VGO14}, has been proposed in order to approximate the conditional expectation from \eqref{eq:SigmatoMatch} and to approximate the appropriate leverage function.
The authors in \cite{H09,VGO14} make use of Monte Carlo techniques, whereas in \cite{C11,EKO12,RMQ07} partial differential equation (PDE) methods are applied.

In this article, the PDE approach is considered.
For the effective calibration, the conditional expectation is then often rewritten as, cf.\ \cite{C11,EKO12,RMQ07}, 
\begin{equation}
\mathbb{E}[\psi^{2}(V_{\tau}) \vert X_{\tau} = x] = \frac{\int_{0}^{\infty} \psi^{2}(v) p(x,v,\tau;X_{0},V_{0})dv}{\int_{0}^{\infty} p(x,v,\tau;X_{0},V_{0})dv}.
\label{eq:CondExpec}
\end{equation}
It can be shown, see e.g.\ \cite{R89}, that the joint density function satisfies the \textit{forward Kolmogorov equation}
\begin{equation}
\begin{array}{lll}
\tfrac{\partial}{\partial \tau} p &=& \tfrac{\partial^{2}}{\partial x^{2}} \left( \tfrac{1}{2} \sigma^{2}_{SLV}\psi^{2}(v)p \right) + \tfrac{\partial^{2}}{\partial x \partial v} \left( \rho \xi \sigma_{SLV}\psi(v)v^{\alpha} p \right) + \tfrac{\partial^{2}}{\partial v^{2}} \left( \tfrac{1}{2} \xi^{2} v^{2\alpha} p \right) \\\\
&& - \ \tfrac{\partial}{\partial x} \left( (r_{d}-r_{f}-\tfrac{1}{2}\sigma^{2}_{SLV}\psi^{2}(v)) p \right) - \tfrac{\partial}{\partial v} \left( \kappa(\eta - v) p \right),
\end{array}
\label{eq:ForwardKolmogorov}
\end{equation}
for $x \in \mathbb{R}, v>0, \tau >0$ and with initial condition $p(x,v,0;X_{0},V_{0}) = \delta(x)\delta(v-V_{0})$ where $\delta$ denotes the Dirac delta function. Once the joint density $p$ is known, one can easily determine the leverage function by computing the integrals in \eqref{eq:CondExpec} and the SLV model is calibrated exactly to the LV model.
By combining \eqref{eq:SigmatoMatch}, \eqref{eq:CondExpec} and \eqref{eq:ForwardKolmogorov} it is readily seen that one has to solve a highly non-linear PDE in order to perform the calibration.

In financial mathematics, convection-diffusion equations of the type \eqref{eq:ForwardKolmogorov} are often discretized by means of finite difference methods, see e.g.\ \cite{C11,RMQ07}. If the parameter $\alpha$ is less or equal to $1/2$, however, it holds that $V_{\tau}=0$ is attainable and defining a proper boundary condition at $v=0$ is a non-trivial task. Moreover, finite difference methods are often not mass-conservative whereas conservation of mass is a key property of forward Kolmogorov equations.
The finite volume method proposed in \cite{EKO12} manages to deal with the issues above. The latter method, however, makes use of a transformation of the original PDE \eqref{eq:ForwardKolmogorov} which incorporates derivatives of the leverage function. As the leverage function is often nonsmooth and only known at a finite number of points, this could lead to undesirable (erratic) behaviour.

In this paper we will introduce a finite volume - alternating direction implicit discretization for the numerical solution of general, non-transformed forward Kolmogorov equations of the type \eqref{eq:ForwardKolmogorov}. 
The discretization makes use of the general \textit{method of lines} (MOL), cf.\ \cite{HV03}.
The PDE is first discretized in the spatial variables, yielding large systems of stiff ordinary differential equations (ODEs). These so-called \textit{semidiscrete systems} are subsequently solved by applying a suitable implicit time stepping method.
The spatial discretization is performed by \textit{finite volume methods} to keep the total numerical mass equal to one and to handle the boundary conditions in a natural way. Since the PDE \eqref{eq:ForwardKolmogorov} is multidimensional, the temporal discretization is performed by an \textit{alternating direction implicit} (ADI) scheme, more precisely the \textit{Hundsdorfer--Verwer} (HV) scheme. This can yield a large computational advantage in comparison with standard (non-split) implicit time stepping methods.
Finally, for the calibration of the SLV model to the LV model, an inner iteration is introduced in order to handle the non-linearity from inserting \eqref{eq:CondExpec} into \eqref{eq:ForwardKolmogorov}.

An outline of the rest of our paper is as follows.

In Section \ref{FV} a finite volume discretization is introduced for the spatial discretization of general 1D and 2D forward Kolmogorov equations. The performance of the finite volume discretization is illustrated by ample numerical experiments.

The spatial discretization results in a large system of ODEs. In Section \ref{ADI} an ADI temporal discretization scheme is applied to increase the computational efficiency in the numerical solution of this ODE system.

In Section \ref{Calibration} the finite volume discretization is used for the calibration of SLV models, yielding a large non-linear system of ODEs. The ADI scheme is applied for the numerical solution of this system of ODEs and an iteration procedure is described for handling the non-linearity.

In Section \ref{Experiments} numerical experiments are presented to illustrate the performance of the obtained calibration procedure and the final Section \ref{Conclusion} concludes.

\setcounter{equation}{0}
\section{Spatial discretization of forward Kolmogorov equations}\label{FV}

In the general MOL approach the PDE is first discretized in the spatial variables by for example finite difference (FD) or finite volume (FV) methods. 
In this section a spatial discretization is proposed for a general two-dimensional forward Kolmogorov equation of the type
\begin{equation}
\tfrac{\partial}{\partial \tau} p + \tfrac{\partial}{\partial x} \left( \mu_{1} p \right) + \tfrac{\partial}{\partial y} \left( \mu_{2} p \right) = \tfrac{\partial^{2}}{\partial x^{2}} \left(\tfrac{1}{2} \sigma_{1}^{2}p \right) + \tfrac{\partial^{2}}{\partial x \partial y} \left( \rho \sigma_{1} \sigma_{2} p \right) + \tfrac{\partial^{2}}{\partial y^{2}} \left(\tfrac{1}{2} \sigma_{2}^{2} p \right),
\label{eq:GeneralForward}
\end{equation}
with $x,y \in \R$, $\tau >0$, and where $\sigma_{1}, \sigma_{2}, \mu_{1}, \mu_{2}$ are real coefficient functions of $x, y,\tau$.
Moreover, the functions  $\sigma_{1}, \sigma_{2}$ are required to be non-negative and it is assumed that there exist values $X_{0}, Y_{0}$ such that the initial function is given by $p(x,y,0)=\delta(x-X_{0})\delta(y-Y_{0})$.
Due to the form of the coefficients it is possible that the spatial domain is naturally restricted. For example, if $\mu_{2}(x,y,\tau) = \kappa(\eta-y)$ and $\sigma_{2}(x,y,\tau) = \xi y^{\alpha}$ with $\kappa, \eta, \xi, \alpha$ strictly positive constants, then the domain in the $y$-direction is naturally restricted to $y \ge 0$, cf.\ PDE \eqref{eq:ForwardKolmogorov}.

Since the solution of forward Kolmogorov equations represents the density of an underlying stochastic process, conservation of mass is a fundamental property and the use of FV schemes is appropriate. 
While FD methods are well known in finance, FV methods are less common in financial applications and we briefly recall the basic idea.

\subsection{Introduction to finite volume discretizations}
\label{FVintro}

Finite volume methods were originally developed to solve conservation laws, or more generally to solve PDEs in conservative form. For example, consider the one-dimensional conservative PDE
\begin{equation}
\label{eq:ConservativePDE}
\tfrac{\partial}{\partial \tau} p + \tfrac{\partial}{\partial x} \left(a(p,x,\tau)p\right) = \tfrac{\partial}{\partial x}\left( b(p,x,\tau) \tfrac{\partial}{\partial x} p \right),
\end{equation}
for $x \in \Omega, \tau >0$, where $\Omega$ is an interval in $\R$.
Both sides of equation \eqref{eq:ConservativePDE} can be integrated in $x$ over an interval (more generally, a cell) $[x_{l}, x_{u}]$ in order to get
\begin{equation}
\tfrac{\partial}{\partial \tau} \int_{x_{l}}^{x_{u}} p dx = f(p,x_{l},\tau) - f(p,x_{u},\tau),
\label{eq:GeneralFlux}
\end{equation}
where $f(p,x,\tau)$ is a function given by 
\begin{equation*}
f(p,x,\tau) = a(p,x,\tau)p - b(p,x,\tau)\tfrac{\partial}{\partial x} p.
\end{equation*}
The function $f$ is typically called the \textit{flux} of $p$ and $f(p,x,\tau)\vert_{x=x_{l}}, f(p,x,\tau)\vert_{x=x_{u}}$ represent the fluxes at the left and right boundaries of the cell $[x_{l}, x_{u}]$.
Relationship \eqref{eq:GeneralFlux} shows that the total integral of $p$, which typically represents a mass, momentum or some similar quantity, changes only as a result of the flux difference over the cell. 
If equation \eqref{eq:ConservativePDE} is considered over a bounded domain $[x_{\min}, x_{\max}]$ and we assume that $f(p,x_{\min},\tau) = f(p,x_{\max},\tau)$ for all $\tau$, i.e.\ the flux at the left boundary is exactly matched by the flux at the right boundary, then the space integral of $p$ over $[x_{\min}, x_{\max}]$ is constant in time. This means that the total mass or momentum is conserved. If the spatial domain $\Omega$ of the PDE is unbounded, and if the interval $[x_{\min}, x_{\max}]$ is wide enough, then one will often have that $f(p,x_{\min},\tau) \approx f(p,x_{\max},\tau) \approx 0$ for all $\tau>0$.

To construct a numerical FV scheme we start with a discretization of the spatial domain. If the spatial domain is unbounded, it needs to be truncated to a wide, finite interval $[x_{\min}, x_{\max}]$. Then, consider the discretization
$$ x_{\min} = x_{1} < x_{2} < \cdots < x_{m} = x_{\max}, $$
of the domain of interest and denote 
\begin{equation*}
\Delta x_{i} = x_{i} - x_{i-1}, \qquad \mathrm{for} \ 2 \le i \le m,
\end{equation*}
and $\Delta x_{1} = \Delta x_{m+1}=0$. Define mid-points 
\begin{equation*}
x_{i-0.5} = x_{i} - \tfrac{1}{2} \Delta x_{i} = \frac{x_{i-1}+x_{i}}{2} \qquad \mathrm{for} \ 2 \le i \le m,
\end{equation*}
and let $\Omega_{i} = [x_{i-0.5}, x_{i+0.5}]$ be cells for $i=1,2,\ldots,m$ where we define $x_{0.5}:=x_{1}$ and $x_{m+0.5}:=x_{m}$.
This yields a \textit{vertex centred grid} with cell vertices $x_{i-0.5}$.
We can now consider the cell average $\overline{p}_{i}(\tau)$ which is defined by
\begin{equation*}
\overline{p}_{i}(\tau) = \frac{1}{x_{i+0.5}-x_{i-0.5}} \int_{\Omega_{i}} p(x,\tau) dx, 
\end{equation*}
and which is typically the quantity that FV schemes approximate.
If we assume that the grid is \textit{smooth} in the sense that $\Delta x_{i+1} - \Delta x_{i} = \mathcal{O}(\Delta x^{2})$ where $\Delta x$ denotes a maximal mesh width, then the cell average $\overline{p}_{i}(\tau)$ is a second-order approximation to $p(x_{i},\tau)$. Differentiating $\overline{p}_{i}(\tau)$ in $\tau$ and using \eqref{eq:GeneralFlux} gives
\begin{equation}
\overline{p}'_{i}(\tau) = \frac{f(p,x_{i-0.5},\tau) - f(p,x_{i+0.5},\tau)}{x_{i+0.5}-x_{i-0.5}}
\label{eq:GeneralIdeaFV}
\end{equation}
which is just another way of stating the conservation property since if we sum over all cells and pull out the derivative in time we find that
\begin{equation*}
\tfrac{\partial}{\partial \tau} \sum_{i=1}^{m} \overline{p}_{i}(\tau) (x_{i+0.5}-x_{i-0.5}) = 0,
\end{equation*}
provided $f(p,x_{\min},\tau) = f(p,x_{\max},\tau)$.

Equation \eqref{eq:GeneralIdeaFV} is typically taken as the starting point for the numerical discretization. 
Denote by 
\begin{equation*} 
P_{i}(\tau) \approx \overline{p}_{i}(\tau), \qquad \mathrm{for} \ 1 \le i \le m,
\end{equation*}
the numerical approximations for the cell averages and let $P$ be the vector that contains these approximations.
The numerical discretization is then defined by 
\begin{equation}
P'_{i}(\tau) = \frac{f_{i-0.5}(P,\tau) - f_{i+0.5}(P,\tau)}{x_{i+0.5}-x_{i-0.5}} = [f_{i-0.5}(P,\tau) - f_{i+0.5}(P,\tau)] \frac{2}{\Delta x_{i} + \Delta x_{i+1}},
\label{eq:GeneralIdeaFVDiscrete}
\end{equation}
where the $f_{i \pm 0.5}$ are numerical fluxes that approximate the exact fluxes $f(p,x_{i \pm 0.5},\tau)$.
By defining a discretization of the type \eqref{eq:GeneralIdeaFVDiscrete}, it readily follows that the total numerical integral (mass)
\begin{equation}
\sum_{i=1}^{m_{1}} P_{i}(\tau) (x_{i+0.5} - x_{i-0.5})
\label{eq:TotalNumericalIntegral}
\end{equation}
stays constant in time provided that $f_{0.5}(P,\tau) = f_{m+0.5}(P,\tau)$.
It is clear that the exact fluxes from \eqref{eq:GeneralIdeaFV} involve the unknown function $p$ at the cell boundaries $x_{i \pm 0.5}$. 
Therefore, we define the numerical fluxes by
\begin{equation}
f_{i \pm 0.5}(P,\tau) = a(P_{i \pm 0.5},x_{i \pm 0.5},\tau)P_{i \pm 0.5} - b(P_{i \pm 0.5},x_{i \pm 0.5},\tau)P_{x,i\pm 0.5},
\label{eq:NumericalFluxesGeneral}
\end{equation}
where the $P_{i \pm 0.5}(\tau)$ form approximations to the exact values $p(x_{i \pm 0.5},\tau)$ and the $P_{x,i\pm 0.5}(\tau)$ form approximations to $\tfrac{\partial}{\partial x} p(x,\tau) \vert_{x=x_{i \pm 0.5}}$. 
Since the cell average is a second-order approximation to $p$, the approximations $P_{i}$ can be used to define $P_{i \pm 0.5}$ and $P_{x,i\pm 0.5}$. 
The manner in which the latter values are computed at the cell boundaries from the surrounding $P_{i}$ plays a large part in defining the characteristics of the numerical scheme. 
Lastly, inserting \eqref{eq:NumericalFluxesGeneral} into \eqref{eq:GeneralIdeaFVDiscrete} yields a system of (possibly non-linear) ODEs which is solved with a suitable time integration procedure.
\newline

Recall that conservation of mass is a fundamental property of forward Kolmogorov equations and the use of FV schemes is appropriate. 
Forward Kolmogorov equations of the type \eqref{eq:GeneralForward} are, however, not in conservative form and hence straightforward application of standard FV schemes is not possible. Moreover, rewriting PDE \eqref{eq:GeneralForward} in conservative form would involve derivatives of the coefficient functions, which are not known in general practical applications. 
In the remainder of this section, a FV-based discretization of the spatial derivatives in the non-transformed PDE \eqref{eq:GeneralForward} is introduced such that conservation of total mass is guaranteed.
We start by explaining the discretization for a general one-dimensional forward Kolmogorov equation, and then generalise it to the two-dimensional case.

\subsection{One-dimensional forward Kolmogorov equations} \label{1DKolmogorov}
\label{FV1D}

Standard one-dimensional forward Kolmogorov equations are also not written in conservative form and their solutions represent density functions of underlying stochastic processes. In this subsection a FV-based discretization is introduced for the general one-dimensional equation 
\begin{equation}
\tfrac{\partial}{\partial \tau} p + \tfrac{\partial}{\partial x} \left( \mu p \right) =  \tfrac{\partial^{2}}{\partial x^{2}} \left( \tfrac{1}{2} \sigma^{2}p \right),
\label{eq:GeneralForward1D}
\end{equation} 
for $x\in \R$, $\tau >0$, where $\sigma, \mu$ are real functions of $x$ and $\tau$, with $\sigma$ non-negative and with initial function given by $p(x,0) = \delta(x-X_{0})$ for some real $X_{0}$.
Spatial discretization by FD or FV methods is often applied on a finite grid. By consequence, the spatial domain has to be truncated to $[x_{\min}, x_{\max}]$, where the boundaries are chosen sufficiently far away from $X_{0}$ such that the truncation error is negligible. 
Recall that the form of $\sigma, \mu$ can naturally restrict the spatial domain of the PDE to for example $x \ge 0$. 
In the latter case, the lower boundary is naturally defined as $x_{\min}=0$. 

As before, define a spatial mesh $x_{\min} = x_{1} < x_{2} < \ldots < x_{m} = x_{\max}$, let $ \Delta x_{i} = x_{i} - x_{i-1}$ be the mesh widths, with $\Delta x_{1} = \Delta x_{m+1} = 0$, and define
\begin{equation*}
x_{i-0.5} = x_{i} - \tfrac{1}{2} \Delta x_{i} = \frac{x_{i-1}+x_{i}}{2} \qquad \mathrm{for} \ 2 \le i \le m,
\end{equation*}
with $x_{0.5} = x_{1}$ and $x_{m+0.5} = x_{m}$.
This yields a vertex centred grid with cells $\Omega_{i} = [x_{i-0.5}, x_{i+0.5}]$.
Let the $P_{i}(\tau)$ denote approximations to the exact cell averages
\begin{equation*}
\overline{p}_{i}(\tau) = \frac{1}{x_{i+0.5}-x_{i-0.5}} \int_{\Omega_{i}} p(x,\tau) dx,
\end{equation*} 
and let $P$ be the vector containing these approximations.
Analogously to the previous section (see equations \eqref{eq:GeneralIdeaFV} and \eqref{eq:GeneralIdeaFVDiscrete}, as well as \cite{HV03}) we define discretizations of the form
\begin{equation}
P_{i}'(\tau) = [f_{i-0.5}(P,\tau) -  f_{i+0.5}(P,\tau)] \frac{2}{\Delta x_{i} + \Delta x_{i+1}}
\label{eq:discretizationforward1dgeneral}
\end{equation} 
where the numerical fluxes are given by
\begin{equation*}
f_{i \pm 0.5}(P,\tau) = f_{a,i \pm 0.5}(P,\tau) + f_{d,i \pm 0.5}(P,\tau),
\end{equation*}
with
\begin{equation}
f_{a,i \pm 0.5}(P,\tau) \approx \mu(x_{i \pm 0.5},\tau) p(x_{i \pm 0.5},\tau),
\label{eq:advectionflux1D}
\end{equation}
and
\begin{equation}
f_{d,i \pm 0.5}(P,\tau) \approx - \tfrac{\partial}{\partial x} \left( \tfrac{1}{2} \sigma^{2}(x,\tau) p(x,\tau) \right)\vert_{x = x_{i \pm 0.5}}.
\label{eq:diffusionflux1D}
\end{equation}
For the ease of presentation, from now on we omit the dependence of the parameters on $\tau$ and set $\mu_{i \pm 0.5} = \mu(x_{i \pm 0.5},\tau)$ and $\sigma_{i} = \sigma(x_{i},\tau)$.
Note that $f_{0,5}(P,\tau)$, respectively $f_{m+0.5}(P,\tau)$, corresponds with the flux at the boundary $x_{\min} = x_{1}$, respectively $x_{\max} = x_{m}$.

The advection part of the PDE \eqref{eq:GeneralForward1D} is written in conservative form. For the inner cell boundaries, i.e.\ for $x_{i-0.5}$ with $2 \le i \le m$, we consider the second-order central FV scheme, cf.\ \cite{HV03}, and define $f_{a,i - 0.5}(P,\tau)$ in \eqref{eq:advectionflux1D} as 
\begin{equation*}
f_{a,i - 0.5}(P,\tau) = \mu_{i - 0.5} \frac{P_{i-1}(\tau) + P_{i}(\tau)}{2}.
\end{equation*}
The diffusion part is not written in conservative form and hence it is not possible to apply standard FV schemes to this term directly. The idea of the second-order FV scheme for \eqref{eq:diffusionflux1D}, see e.g.\ \cite{HV03}, is generalised by defining
\begin{equation*}
f_{d,i - 0.5}(P,\tau) = -\left( \tfrac{1}{2} \sigma^{2}_{i} P_{i}(\tau) - \tfrac{1}{2} \sigma^{2}_{i-1} P_{i-1}(\tau) \right) \frac{1}{\Delta x_{i}}
\end{equation*}
for $2 \le i \le m$.
It is readily seen that 
$$  \left(\tfrac{1}{2} \sigma^{2}_{i}p(x_{i},\tau) - \tfrac{1}{2} \sigma^{2}_{i-1} p(x_{i-1},\tau)\right) \frac{1}{\Delta x_{i}} $$
is a second-order approximation of $ \tfrac{\partial}{\partial x} ( \tfrac{1}{2} \sigma^{2}p)$ at the point $x_{i-0.5}$ which explains the choice for this discretization.
Inserting these expressions back into \eqref{eq:discretizationforward1dgeneral} we get
\begin{eqnarray}
P_{i}'(\tau) &=& \frac{\sigma^{2}_{i-1}P_{i-1}(\tau)}{\Delta x_{i} (\Delta x_{i} + \Delta x_{i+1})} - \frac{\sigma_{i}^{2}P_{i}(\tau)}{\Delta x_{i} \Delta x_{i+1}} + \frac{\sigma^{2}_{i+1}P_{i+1}(\tau)}{\Delta x_{i+1} (\Delta x_{i} + \Delta x_{i+1})} \label{eq:Discretization1DWithoutBoundary} \\ \nonumber \\
&& + \  \left[ \mu_{i-0.5} \frac{P_{i-1}(\tau)+P_{i}(\tau)}{2} - \mu_{i+0.5} \frac{P_{i}(\tau)+P_{i+1}(\tau)}{2} \right] \frac{2}{\Delta x_{i} + \Delta x_{i+1}}, \nonumber
\end{eqnarray}
for $ 2 \le i \le m-1$.
Note that by applying the second-order central FD scheme for diffusion on non-uniform spatial grids, cf.\ \cite{IHF10}, on the term $\tfrac{\partial^{2}}{\partial x^{2}} ( \tfrac{1}{2} \sigma^{2}p)$, one would end up with the same discretization for the diffusion term.

To complete this semidiscretization, it also has to be defined at the boundaries of the truncated domain such that conservation of the total mass is guaranteed.
Given that $p$ represents a density function, it follows that
$$ \int_{-\infty}^{\infty} p(x,\tau) dx = 1, \qquad \forall \tau > 0,$$
and hence
$$ \int_{-\infty}^{\infty} \left[ \tfrac{\partial}{\partial \tau} p \right] dx = \int_{-\infty}^{\infty} \left[ \tfrac{1}{2} \tfrac{\partial^{2}}{\partial x^{2}} \left( \sigma^{2}p \right) - \tfrac{\partial}{\partial x} \left( \mu p \right) \right] dx = 0. $$
Assuming that $[x_{\min},x_{\max}]$ is chosen sufficiently wide, the condition above can be approximated by
$$ \left. \left[ \tfrac{\partial}{\partial x} \left( \tfrac{1}{2} \sigma^{2}p \right) - \left( \mu p \right) \right]\right\vert^{x= x_{\max}}_{x=x_{\min}} = 0, $$
reflecting the fact that the total flux over the interval $[x_{\min},x_{\max}]$ is zero. 
As stated above, for some choices of coefficient functions the spatial domain is naturally restricted. For example, if the PDE \eqref{eq:GeneralForward1D} stems from a non-negative process, $x_{\min}$ can be set equal to zero. The no-flux boundary condition above then still holds on the naturally restricted domain.

In theory it is possible that there is a positive flux at one of the boundaries, and exactly the same negative flux at the other boundary. 
However, since the solution represents a density function, it is more realistic to impose that no mass is coming in or going out at each of the boundaries. In light of this, we assume that the following boundary conditions hold:
\begin{eqnarray}
\left. \left[ \tfrac{\partial}{\partial x} \left( \tfrac{1}{2} \sigma^{2}p \right) - \left( \mu p \right) \right]\right\vert_{x=x_{\min}} &=& 0, \label{eq:BC1Dleft} \\ \nonumber \\
\left. \left[ \tfrac{\partial}{\partial x} \left( \tfrac{1}{2} \sigma^{2}p \right) - \left( \mu p \right) \right]\right\vert_{x=x_{\max}} &=& 0. \nonumber
\end{eqnarray}
The numerical equivalent of the first condition is to say that the flux at $x_{1} = x_{\min}$ is zero, i.e.\ $f_{1}(P,\tau) \equiv f_{0.5}(P,\tau) = 0$. This can be achieved by creating a ghost point $x_{0} = x_{1} - \Delta x_{2}$ and using \eqref{eq:BC1Dleft} to define the value of $\tfrac{1}{2} \sigma_{0}^{2}P_{0}(\tau)$ at the ghost point by
\begin{equation*}
\frac{\tfrac{1}{2}\sigma_{2}^{2}P_{2}(\tau)-\tfrac{1}{2}\sigma_{0}^{2}P_{0}(\tau)}{2\Delta x_{2}} - \mu_{1} P_{1}(\tau) = 0,
\end{equation*} 
where we make use of the fact that the cell averages form second-order approximations to the point values.
Turning to \eqref{eq:diffusionflux1D} we define the diffusive flux $f_{d,0.5}$ at $x_{1}$ as
\begin{equation*}
f_{d,0.5}(P,\tau) = - \frac{\tfrac{1}{2}\sigma_{2}^{2}P_{2}(\tau)-\tfrac{1}{2}\sigma_{0}^{2}P_{0}(\tau)}{2\Delta x_{2}} = - \mu_{1} P_{1}(\tau).
\end{equation*}
Since $x_{1}$ is the left boundary of the first cell, the flux on the boundary $x_{\min}$ stemming from the advection part (see \eqref{eq:advectionflux1D}) can be approximated as $f_{a,0.5} = \mu_{1}P_{1}(\tau)$. Inserting these expressions into \eqref{eq:discretizationforward1dgeneral} we obtain
\begin{equation}
P_{1}'(\tau) = - f_{1.5}(P,\tau) \frac{2}{\Delta x_{1} + \Delta x_{2}} = - f_{1.5}(P,\tau) \frac{2}{\Delta x_{2}}.
\label{eq:DiscretizatinBounderayLeft1D}
\end{equation}
The boundary condition at $x_{\max}$ can be handled analogously in order to get
\begin{equation}
P_{m}'(\tau) = f_{m-0.5}(P,\tau) \frac{2}{\Delta x_{m}}.
\label{eq:DiscretizatinBounderayRight1D}
\end{equation}
By performing the discretization of the boundary conditions in this way, it follows that $f_{0.5}(P,\tau)=f_{m+0.5}(P,\tau)$ and we ensure that mass is conserved in the numerical scheme. 

Combining \eqref{eq:Discretization1DWithoutBoundary}, \eqref{eq:DiscretizatinBounderayLeft1D} and \eqref{eq:DiscretizatinBounderayRight1D} we see that the total discretization can be written as a system of ODEs
\begin{equation}
P'(\tau) = A(\tau) P(\tau)
\label{eq:ODE1D}
\end{equation} 
for $\tau > 0$, with given matrix $A(\tau)$. 
Since $P_{i}(\tau)$ represent cell averages it is natural to define the initial vector as 
$$ P_{i}(0) = \left\{ \begin{array}{lll}
\tfrac{2}{\Delta x_{i} + \Delta x_{i+1}} & & \mathrm{if} \ X_{0} \in [x_{i-0.5}, \ x_{i+0.5}], \\\\
0 && \mathrm{otherwise.}
\end{array} \right. $$
In general, the exact solution of the system of ODEs \eqref{eq:ODE1D} can not be computed analytically and one relies on numerical methods in order to approximate it. Since the discretization above often leads to stiff semidiscrete systems, suitable implicit time stepping schemes such as the \textit{Crank--Nicolson scheme} are widely considered, see e.g.\ \cite{TR00}.

\subsection{Numerical experiment} \label{subsec:Experiment1D}

In this subsection the performance of the FV discretization is tested by considering two practical examples. 
As a first example, consider the SDE
$$ dS_{\tau} = (r_{d}-r_{f})S_{\tau} d\tau + \sigma_{BS} S_{\tau} dW_{\tau}, $$
with $\sigma_{BS} > 0$,
corresponding to the classical Black--Scholes model.
Then, the underlying density is known exactly and given by
\begin{equation}
p(s,\tau) = \frac{1}{\sigma_{BS}\sqrt{\tau}}\phi\left(\frac{\log(s/S_{0})-(r_{d}-r_{f}-\tfrac{1}{2}\sigma^{2}_{BS})\tau}{\sigma_{BS}\sqrt{\tau}}\right)\frac{1}{s}, \qquad \mathrm{for} \ s>0, \tau > 0,
\label{eq:DensityBS}
\end{equation}
where $\phi(x)$ is the density function of a standard normally distributed random variable.
The density function $p(s,\tau)$ from \eqref{eq:DensityBS} satisfies the PDE
$$ \tfrac{\partial}{\partial \tau} p = \tfrac{\partial^{2}}{\partial s^{2}} \left( \tfrac{1}{2} \sigma_{BS}^{2}s^{2}p \right) - \tfrac{\partial}{\partial s} \left( (r_{d} - r_{f})s p \right), $$
for $s, \tau > 0$, with $p(s,0)=\delta(s-S_{0})$. This PDE is of the form \eqref{eq:GeneralForward1D} with a natural restriction of the spatial domain. Note that for the numerical experiment we don't apply the $\log$-transformation from Section \ref{intro} in order to have non-constant coefficients which makes the problem more challenging.

\begin{figure}
\begin{center}
\includegraphics[scale=0.5]{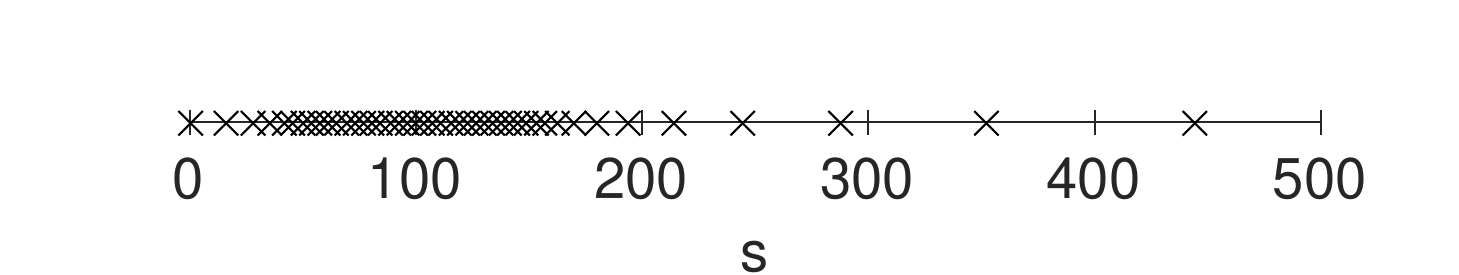}  
\caption{Illustration of the non-uniform grid around $s=S_{0}$ for the Black--Scholes example and the actual values $S_{0}=100, m=50$.}
\label{fig:Grid1D}
\end{center}
\end{figure}

Firstly, the spatial domain is truncated to $[S_{\min}, S_{\max}] = [0, 30S_{0}]$ and we define a non-uniform grid $S_{\min}=s_{1} < s_{2} < \cdots < s_{m} = S_{\max}$ by making use of a $\sinh$ transformation of a uniform underlying grid, cf.\ \cite{IHF10}. In Figure \ref{fig:Grid1D} the spatial grid is shown for the values $S_{0}=100, m = 50$ and from $s=0$ to $s=5S_{0}$ to illustrate the smaller mesh widths around the point $s=S_{0}$.
Applying the FV discretization from Subsection \ref{1DKolmogorov} then yields approximations $P_{i}(\tau)$ of the exact values $p(s_{i},\tau)$.

\begin{figure}
\begin{center}
\includegraphics[scale=0.5]{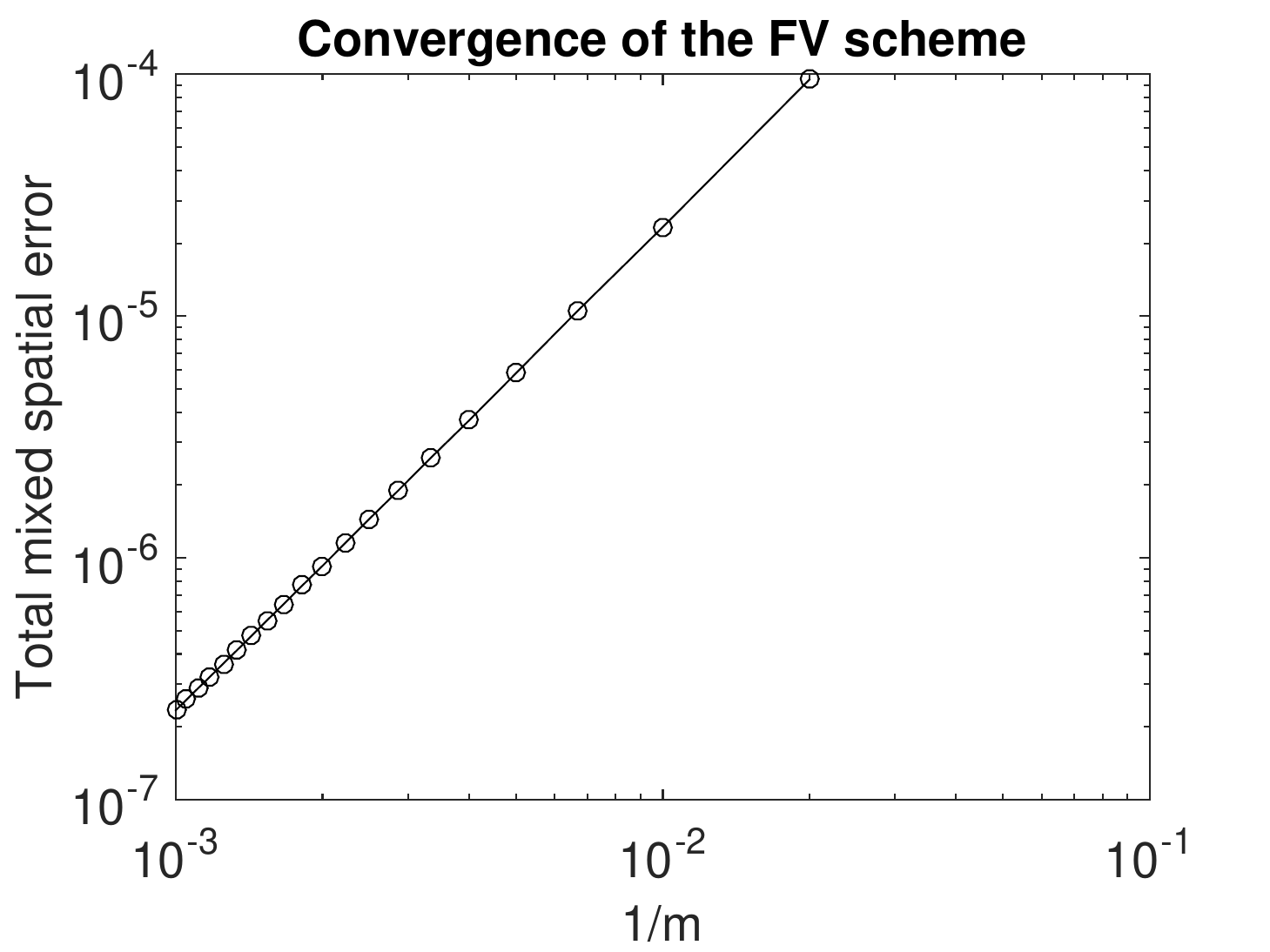}  
\includegraphics[scale=0.5]{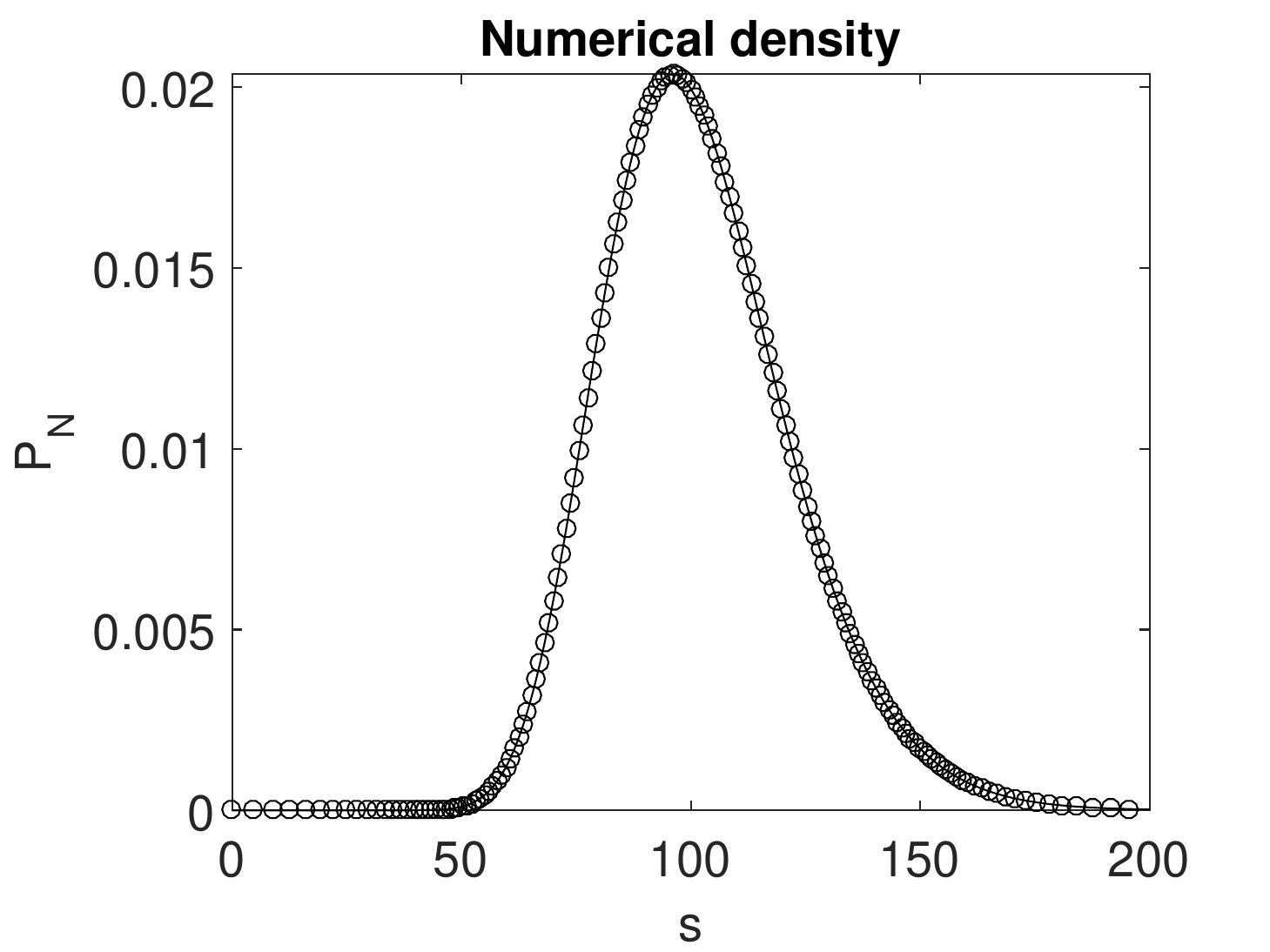}
\caption{Convergence results within the 1D Black--Scholes model. The parameter values are $r_{d} = 0.03, r_{f} = 0.01, \sigma_{BS} = 0.2, T=1$.}
\label{fig:Convergence1D_BS}
\end{center}
\end{figure}

When trying to determine the performance of a numerical method with respect to a reference solution, it is important to take note of the computational environment in which values are calculated, and to understand the impact that has on the comparison. Our calculations take place in 64 bit IEEE floating point arithmetic. Since the solution of forward Kolmogorov equations represents a probability density function, the magnitude of the solution varies dramatically over the computational domain. This is especially true of the initial condition (a Dirac delta), and also more generally with naturally bounded stochastic processes with an attainable boundary, c.f.\ the SLV model \eqref{eq:SLVmodel} with $0 < \alpha \leq 1/2$. Since IEEE floating point has a fixed-length mantissa, the density function cannot be represented to a high \textit{absolute} accuracy uniformly over the domain. In areas where the density function is large, only high \textit{relative} accuracy (correct number of digits) can be achieved.
In addition, since the numerical solution is obtained by using implicit time stepping, we should not expect high relative accuracy of the numerical density in regions where the exact solution is small.  This is because when solving the linear systems we combine many terms with very different magnitudes and then sum them up, which in IEEE arithmetic will lead to a loss of relative accuracy.
Therefore when comparing the numerical solution and the reference solution we adopt a mixed absolute-relative error metric: we use relative error when the reference solution is larger than 1, absolute error if the reference solution is less than 1, and we take the maximal error value over the whole domain.
More precisely, let
\begin{equation*}
\epsilon_{i}(m) = \left\{ \begin{array}{ll}
\left\vert \tfrac{p(s_{i},T) - P_{i}(T)}{p(s_{i},T)} \right\vert \qquad & \mathrm{if} \ p(s_{i},T) > 1, \\\\
\vert p(s_{i},T) - P_{i}(T) \vert \qquad & \mathrm{else}. 
\end{array} \right.
\end{equation*}
The total mixed spatial error is then defined by 
$$ \epsilon(m) = \max_{1 \le i \le m} \epsilon_{i}(m). $$
The value of 1 is somewhat arbitrary. The results, however, are not that sensitive to the crossover value as long as it is not too small.
For the actual experiments, the values $P_{i}(T)$ are approximated by applying the Crank--Nicolson time stepping scheme with a large number of steps such that the temporal discretization error is negligible.
In the left plot of Figure \ref{fig:Convergence1D_BS} the total mixed spatial error is shown for the relevant situation where $r_{d} = 0.03$, $r_{f} = 0.01$, $\sigma_{BS} = 0.2, T=1$ and for the number of spatial grid points $m = \{ 50, 100, \ldots, 1000\}$. 
The corresponding numerical solution for $m=200$ is shown in the right plot of Figure \ref{fig:Convergence1D_BS}.
The convergence plot clearly indicates that the FV discretization is second-order convergent with respect to the current initial-boundary value problem.

As a second example we consider the Cox--Ingersoll--Ross (CIR) process, cf.\ \cite{CIR85},
\begin{equation*}
dV_{\tau} = \kappa (\eta - V_{\tau}) d\tau + \xi \sqrt{V_{\tau}} dW_{\tau},
\end{equation*}
where $\kappa, \eta, \xi$ are strictly positive parameters. The corresponding density function is given by, see e.g.\ \cite{CIR85}, 
\begin{equation}
p(v,\tau) = c e^{-u_{0}-u_{1}}(\tfrac{u_{1}}{u_{0}})^{q/2} \rm{I_{q}}(2\sqrt{u_{0}u_{1}}),
\label{eq:CIRdensity}
\end{equation}
where
\begin{equation*}
c = \tfrac{2\kappa}{\xi^{2}(1-e^{-\kappa\tau})}, \quad u_{0} = c V_{0} e^{-\kappa \tau}, \quad u_{1} = c v, \quad q = \tfrac{2\kappa \eta}{\xi^{2}}-1,
\end{equation*}
and $\rm{I_{q}}(\cdot)$ is the modified Bessel function of the first kind of order $q$. Note that the value of $q$ is directly related with the so-called Feller condition, i.e.\ with the possibility that $V_{\tau} = 0$ is attainable. The density function \eqref{eq:CIRdensity} satisfies the forward Kolmogorov equation 
\begin{equation}
\tfrac{\partial}{\partial \tau} p = \tfrac{\partial^{2}}{\partial v^{2}} \left( \tfrac{1}{2} \xi^{2}vp \right) - \tfrac{\partial}{\partial v} \left( \kappa (\eta  - v) p \right), 
\label{eq:CIRPDE}
\end{equation}
for $v, \tau > 0$, with $p(v,0)=\delta(v-V_{0})$. It is readily seen that if the Feller condition is violated, i.e.\ if $q<0$, then the density from \eqref{eq:CIRdensity} is not defined at $v=0$ and the density function tends to infinity as $v$ tends to zero.
In addition, around $v=0$ the PDE \eqref{eq:CIRPDE} is strongly convection dominated which is very challenging for numerical discretization methods.

\begin{figure}
\begin{center}
\includegraphics[scale=0.5]{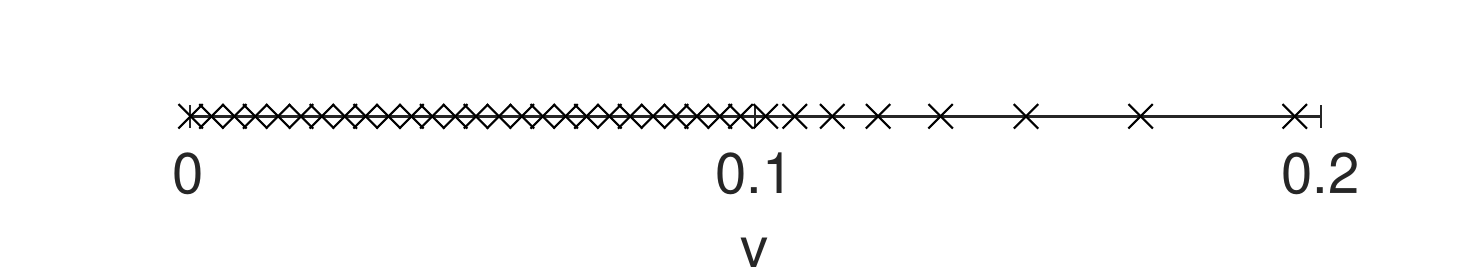}  
\caption{Illustration of the non-uniform grid around $v=0$ and $v=V_{0}$ for the CIR example and the actual values $V_{0}=0.0625, m=50$.}
\label{fig:Grid1DCIR}
\end{center}
\end{figure}

The domain is truncated to $[V_{\min}, V_{\max}] = [0, 15]$ and a non-uniform grid $0=v_{1} < v_{2} < \cdots < v_{m}=V_{\max}$ is defined by making use of a $\sinh$ transformation of a uniform underlying grid. The spatial grid is shown in Figure \ref{fig:Grid1DCIR} for the values $V_{0} = 0.0625, m=50$ and from $v=0$ to $v=0.2$ to illustrate the smaller mesh widths around $v=0$ and $v=V_{0}$.
Afterwards the FV discretization is applied which leads to approximations $P_{j}(T), \ (1 \le j \le m)$.
Recall that if $q<0$, then the density function tends to infinity as $v$ tends to zero and at $v=0$ the exact density function is not defined. By increasing the number of spatial grid points $m$, the value of the second grid point $v_{2}$ tends to zero and adequately comparing the difference between $ p(v_{2},T)$ and $P_{2}(T)$ becomes difficult.
In view of this, we opt to compute the error on similar spatial domains. 
Let $v_{low}$ be the smallest non-zero grid point, i.e.\ the point $v_{2}$, if the total number of spatial grid points $m$ is 50.
We then define the total mixed spatial error by 
$$ \max_{j_{1} \le j \le m} \epsilon_{j}(m), $$
where, for given $m$, $j_{1}$ is the lowest index such that $v_{j_{1}} \ge v_{low}$ and
\begin{equation*}
\epsilon_{j}(m) = \left\{ \begin{array}{ll}
\left\vert \tfrac{p(v_{j},T) - P_{j}(T)}{p(v_{j},T)} \right\vert \qquad & \mathrm{if} \ p(v_{j},T) > 1, \\\\
\vert p(v_{j},T) - P_{j}(T) \vert \qquad & \mathrm{else}. 
\end{array} \right.
\end{equation*}
The approximations $P_{j}(T)$ are determined by considering a large number of time steps such that the temporal discretization error is negligible. Please note that the choice $m=50$ for defining $v_{low}$ is not crucial. The conclusions of the numerical experiments are essentially unchanged as long as $v_{low}$ is defined via one of the coarsest grids considered in the experiment.

\begin{figure}
\begin{center}
\includegraphics[scale=0.5]{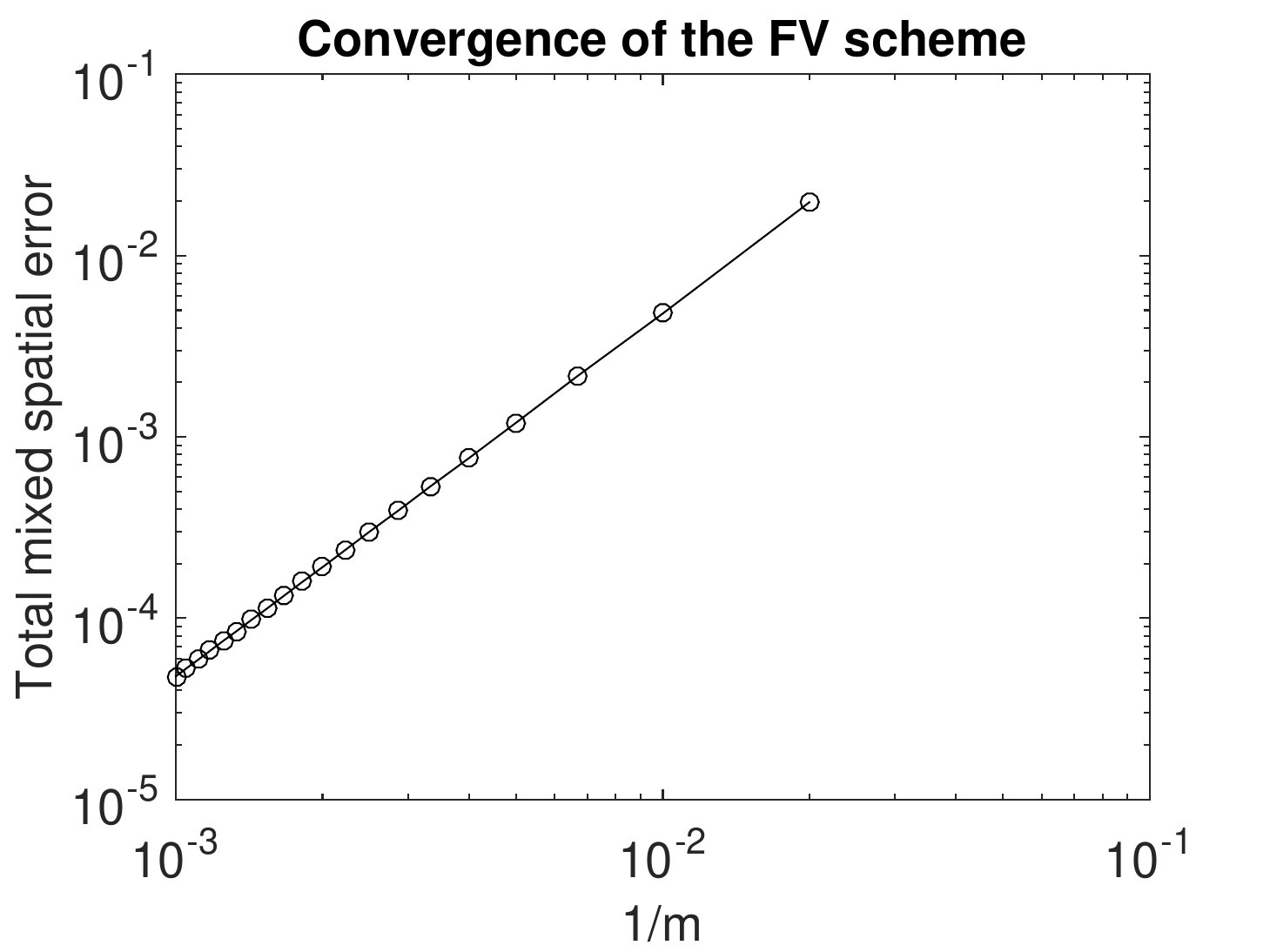}  
\includegraphics[scale=0.5]{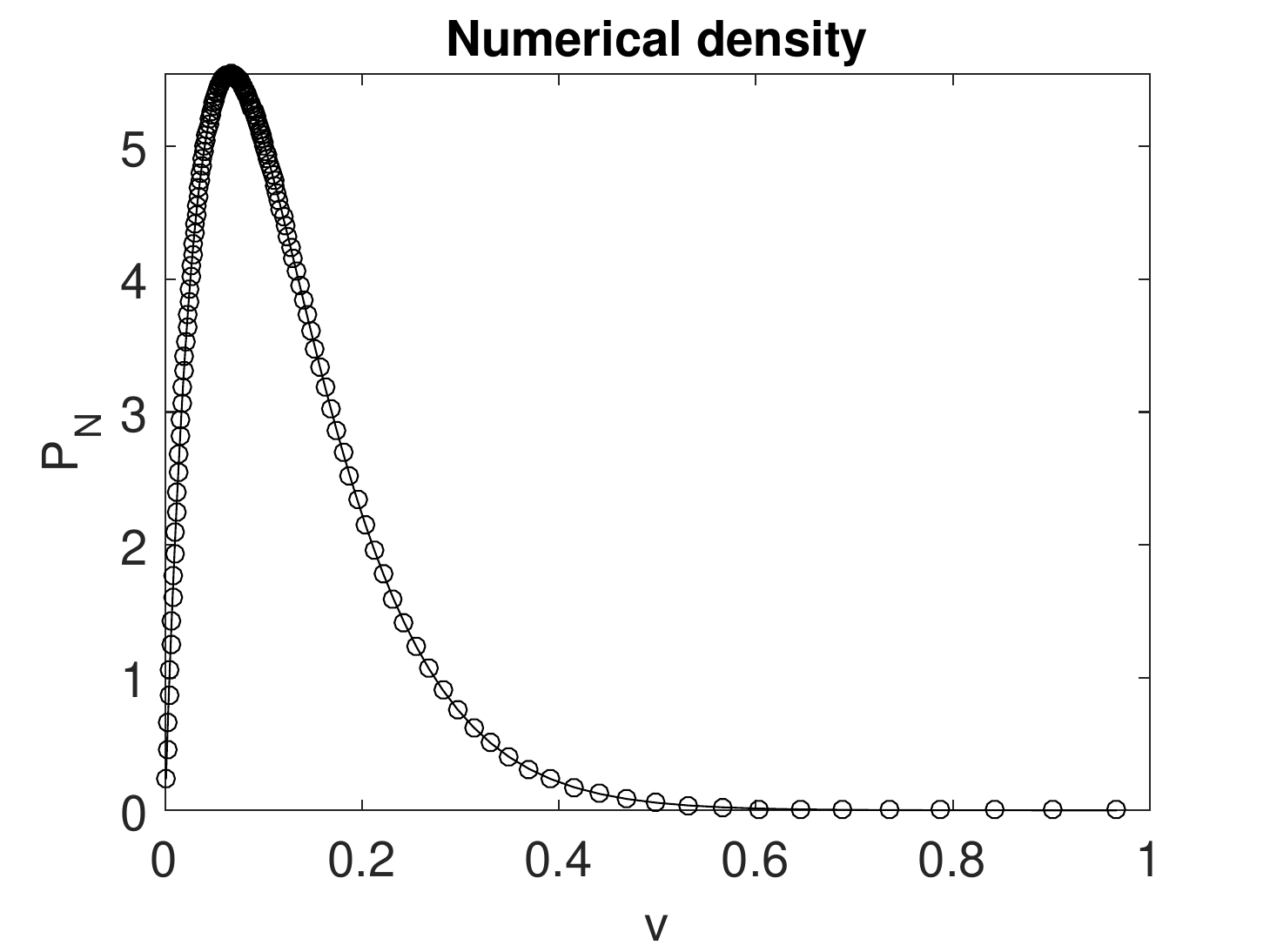}
\caption{Convergence results within the CIR model. The parameters are given by Set A.}
\label{fig:Convergence1D_CIRSetA}
\end{center}
\end{figure}

\begin{figure}
\begin{center}
\includegraphics[scale=0.5]{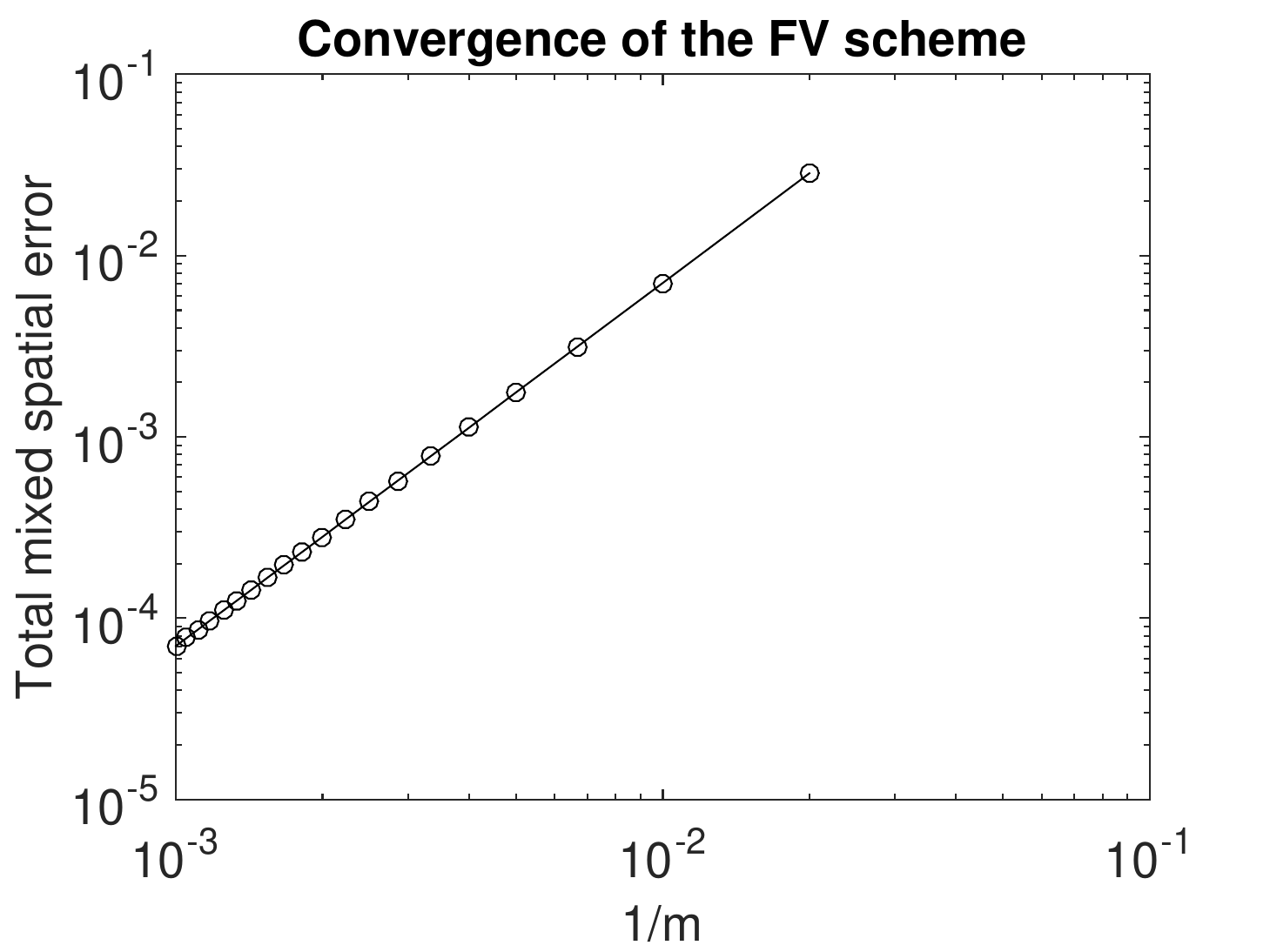}  
\includegraphics[scale=0.5]{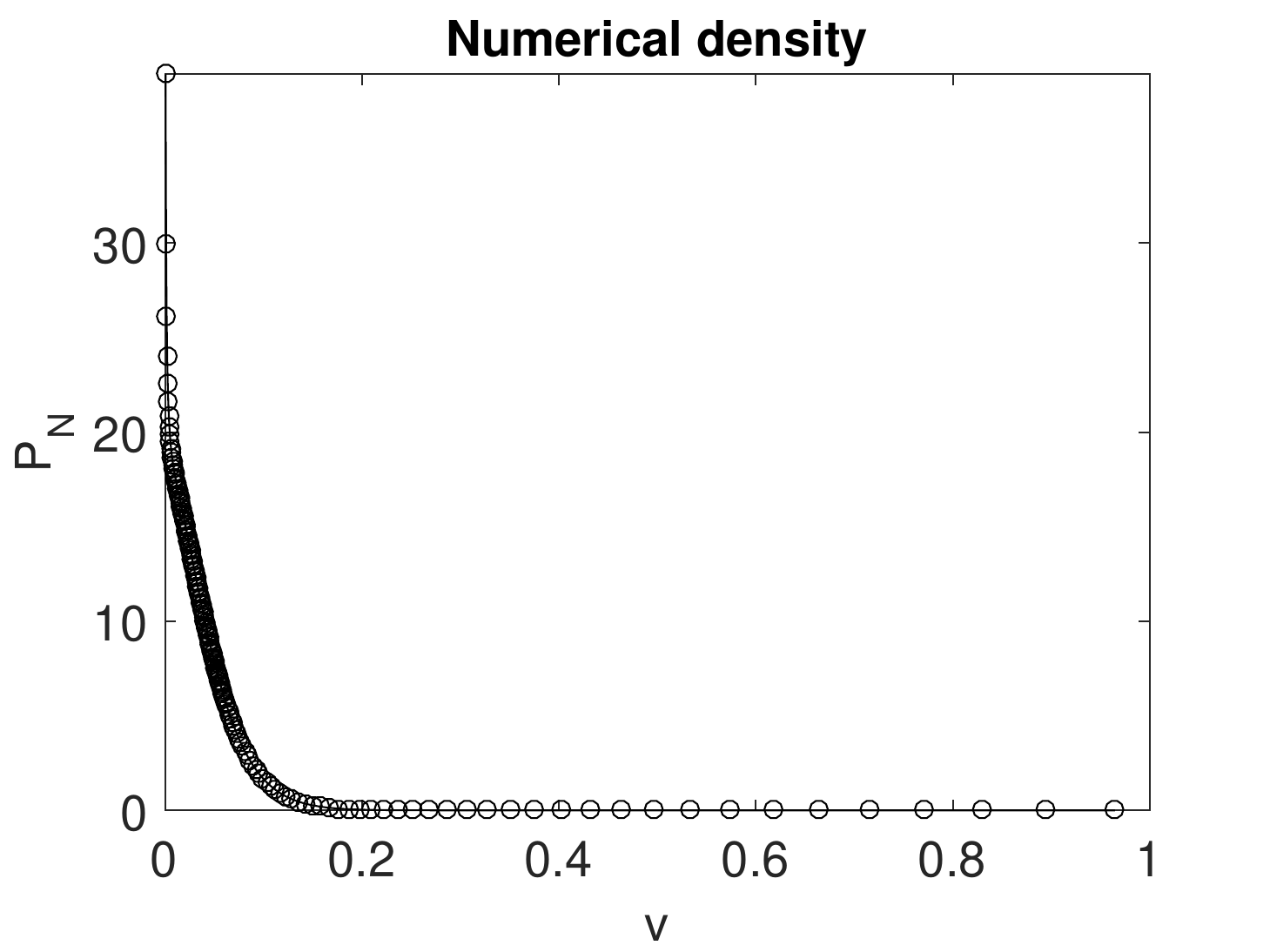}
\caption{Convergence results within the CIR model. The parameters are given by Set B.}
\label{fig:Convergence1D_CIRSetB}
\end{center}
\end{figure}

For the actual experiment we consider two sets of parameters:
\begin{center}
\begin{tabular}{|c|c|c|c|c|c|}
\hline
& $\kappa$ & $\eta$ & $\xi$ & $V_{0}$ & $T$ \\
\hline
Set A & 5 & 0.16 & 0.9 & 0.0625 & 0.25 \\
\hline
Set B & 1.15 & 0.0348 & 0.39 & 0.0348 & 0.25 \\
\hline
\end{tabular}
\end{center}
These sets are taken from \cite{FO11} and were also used in \cite{RO12}. For Set A we have $q=0.98$ and the variance process remains strictly positive. For Set B we have $q=-0.47$ and $V_{\tau}=0$ is attainable.
In the left plot of Figure \ref{fig:Convergence1D_CIRSetA}, respectively Figure \ref{fig:Convergence1D_CIRSetB}, the total mixed spatial error is shown for the parameters of Set A, respectively Set B, and for the number of spatial grid points $m = \{ 50, 100, \ldots, 1000\}$.
In the right plots, the corresponding numerical solutions are shown for $m=200$.
The convergence plots indicate that the FV discretization is convergent with respect to the current initial-boundary value problems. Additional experiments suggest that FV discretization is second-order convergent if the Feller condition is satisfied. If $q<0$ the order of convergence can drop to one.
In addition, all the experiments confirm that the total numerical mass \eqref{eq:TotalNumericalIntegral} stays constantly equal to one, even if the Feller condition is strongly violated.

\subsection{Two-Dimensional forward Kolmogorov equations} \label{2DKolmogorov}

In this subsection, the FV discretization from the one-dimensional case is used to define a spatial discretization for the general two-dimensional forward Kolmogorov equation \eqref{eq:GeneralForward}. Suppose the spatial domain is truncated to $[x_{\min}, x_{\max}] \times [y_{\min}, y_{\max}]$, and the spatial grid points in the $x$-direction, respectively $y$-direction, are given by
$$x_{\min} = x_{1} < x_{2} < \ldots < x_{m_{1}} = x_{\max},$$
respectively
$$y_{\min} = y_{1} < y_{2} < \ldots < y_{m_{2}} = y_{\max}.$$
Let $\Delta x_{i} = x_{i}-x_{i-1}$ and $\Delta y_{j} = y_{j}-y_{j-1}$ be the spatial mesh widths, where $\Delta x_{1} = \Delta x_{m_{1}+1} = \Delta y_{1} = \Delta y_{m_{2}+1} = 0$, and define volumes
$$ \Omega_{i,j} := [x_{i-0.5}, x_{i+0.5}] \times [y_{j-0.5}, y_{j+0.5}], $$
where $x_{i \pm 0.5}$, respectively $y_{j \pm 0.5}$, are defined analogously as in the one-dimensional case.

We can now consider the two-dimensional equivalent of the cell averages, volume averages, $\overline{p}_{i,j}(\tau)$ which are defined by
$$ \overline{p}_{i,j}(\tau) = \frac{1}{\vert \Omega_{i,j} \vert} \int_{\Omega_{i,j}} p(x,y,\tau) dx dy, $$
with $\vert \Omega_{i,j} \vert = (x_{i+0.5}-x_{i-0.5})(y_{j+0.5}-y_{j-0.5})$ the area of the corresponding volume.
The volume average $\overline{p}_{i,j}(\tau)$ is again a second-order approximation to $p(x_{i},y_{j},\tau)$, provided that the underlying meshes are smooth, and is the quantity that is approximated by the FV discretization.
It is readily verified that
\begin{subeqnarray}
\vert \Omega_{i,j} \vert \frac{\partial}{\partial \tau} \overline{p}_{i,j}(\tau) &=& \int_{y_{j-0.5}}^{y_{j+0.5}} \left. \left[ \frac{\partial}{\partial x} \left( \frac{1}{2} \sigma_{1}^{2} p \right) - \mu_{1}p \right] \right\vert_{x=x_{i-0.5}}^{x=x_{i+0.5}} dy \\
&& + \ \int_{x_{i-0.5}}^{x_{i+0.5}} \left. \left[ \frac{\partial}{\partial y} \left( \frac{1}{2} \sigma_{2}^{2} p \right) - \mu_{2}p \right] \right\vert_{y=y_{j-0.5}}^{y=y_{j+0.5}} dx \\
&& + \ \left. \left[ \left. \left[ \rho \sigma_{1}\sigma_{2} p \right] \right\vert_{x=x_{i-0.5}}^{x=x_{i+0.5}} \right] \right\vert_{y=y_{j-0.5}}^{y=y_{j+0.5}},
\label{eq:2DAverage}
\end{subeqnarray}
and the two-dimensional discretization is based on the approximation
\begin{subeqnarray}
\vert \Omega_{i,j} \vert \frac{\partial}{\partial \tau} \overline{p}_{i,j}(\tau) &\approx& \bigg[ \left. \left[ \frac{\partial}{\partial x} \left( \frac{1}{2} \sigma_{1}^{2} p \right) - \mu_{1}p \right] \right\vert_{x=x_{i-0.5}}^{x=x_{i+0.5}} \bigg] \bigg\vert_{y=y_{j}} \frac{\Delta y_{j} + \Delta y_{j+1}}{2}  \slabel{eq:FluxX} \\
&& + \ \bigg[ \left. \left[ \frac{\partial}{\partial y} \left( \frac{1}{2} \sigma_{2}^{2} p \right) - \mu_{2}p \right] \right\vert_{y=y_{j-0.5}}^{y=y_{j+0.5}} \bigg] \bigg\vert_{x=x_{i}} \frac{\Delta x_{i} + \Delta x_{i+1}}{2} \slabel{eq:FluxY} \\
&& + \ \left. \left[ \left. \left[ \rho \sigma_{1}\sigma_{2} p \right] \right\vert_{x=x_{i-0.5}}^{x=x_{i+0.5}} \right] \right\vert_{y=y_{j-0.5}}^{y=y_{j+0.5}}. \slabel{eq:FluxXY}
\label{eq:2DAverageApproximation}
\end{subeqnarray}
Equation \eqref{eq:2DAverage} includes several flux terms that are similar to the flux terms in Subsections \ref{FVintro} and \ref{1DKolmogorov}. It reflects the fact that the total integral of $p$ over a volume changes only as a result of the flux difference over the volume boundary. This is completely analogous with the one-dimensional interpretation of \eqref{eq:GeneralFlux}. If the total flux over the boundary of the spatial domain is zero, i.e.\ if
\begin{eqnarray*}
0 &=& \int_{y_{\min}}^{y_{\max}} \left. \left[ \frac{\partial}{\partial x} \left( \frac{1}{2} \sigma_{1}^{2} p \right) - \mu_{1}p \right] \right\vert_{x=x_{\min}}^{x=x_{\max}} dy \\
&& + \ \int_{x_{\min}}^{x_{\max}} \left. \left[ \frac{\partial}{\partial y} \left( \frac{1}{2} \sigma_{2}^{2} p \right) - \mu_{2}p \right] \right\vert_{y=y_{\min}}^{y=y_{\max}} dx \\
&& + \ \left. \left[ \left. \left[ \rho \sigma_{1}\sigma_{2} p \right] \right\vert_{x=x_{\min}}^{x=x_{\max}} \right] \right\vert_{y=y_{\min}}^{y=y_{\max}},
\end{eqnarray*}
then the total integral of $p$ over the entire domain is constant in time.

Let $\boldsymbol{P}_{i,j}(\tau)$ denote approximations to the exact value $\overline{p}_{i,j}(\tau)$, denote by $\boldsymbol{P}(\tau)$ the $m_{1}\times m_{2}$ matrix with entries $\boldsymbol{P}_{i,j}(\tau)$ and let
$$P(\tau) = \mathrm{vec}[\boldsymbol{P}(\tau)],$$ 
where $\mathrm{vec}[\cdot]$ denotes the operator that turns any given matrix into a vector by putting its successive columns below each other.
The bold notation is only introduced to indicate the subtle difference between the matrix form and the vectorised form of the approximations.
Similarly to the one-dimensional discretization in Subsections \ref{FVintro} and \ref{1DKolmogorov}, we discretize \eqref{eq:2DAverageApproximation} in the following way by introducing numerical fluxes
\begin{subeqnarray}
\boldsymbol{P}'_{i,j}(\tau) &=& \left[ f_{i-0.5,j}(P,\tau) - f_{i+0.5,j}(P,\tau) \right] \frac{2}{\Delta x_{i}+\Delta x_{i+1}} \\
        && + \ \left[ f_{i,j-0.5}(P,\tau) - f_{i,j+0.5}(P,\tau) \right] \frac{2}{\Delta y_{j}+\Delta y_{j+1}} \\
        && + \ \sum_{i_{1},j_{1} = 0}^{1} (-1)^{i_{1}+j_{1}} f_{m,i-0.5 + i_{1},j-0.5+j_{1}} \frac{2}{\Delta x_{i}+\Delta x_{i+1}}\frac{2}{\Delta y_{j}+\Delta y_{j+1}},
\label{eq:TotalDiscretization2D}
\end{subeqnarray}
for $1\le i \le m_{1}, 1 \le j \le m_{2}$.
For the ease of presentation, denote 
\begin{equation*}
\begin{array}{ll}
\mu_{1,i \pm 0.5,j} = \mu_{1}(x_{i \pm 0.5},y_{j},\tau), \quad &\sigma_{1,i,j} = \sigma_{1}(x_{i},y_{j},\tau), \\
\mu_{2,i,j \pm 0.5} = \mu_{2}(x_{i},y_{j \pm 0.5},\tau), \quad &\sigma_{2,i,j} = \sigma_{2}(x_{i},y_{j},\tau), \\
\sigma_{1,i \pm 0.5,j \pm 0.5} = \sigma_{1}(x_{i \pm 0.5},y_{j \pm 0.5},\tau), \quad & \sigma_{2,i \pm 0.5,j \pm 0.5} = \sigma_{2}(x_{i \pm 0.5},y_{j \pm 0.5},\tau).
\end{array}
\end{equation*}
Since the actual form of the fluxes in \eqref{eq:2DAverageApproximation} is completely similar to the form of the fluxes in Subsection \ref{2DKolmogorov}, we define the numerical fluxes by
\begin{equation*}
f_{i \pm 0.5,j}(P,\tau) = f_{a,i \pm 0.5,j}(P,\tau) + f_{d,i \pm 0.5,j}(P,\tau),
\end{equation*}
with
\begin{equation*}
f_{a,i - 0.5,j}(P,\tau) = \mu_{1,i - 0.5,j} \frac{\boldsymbol{P}_{i - 1,j}(\tau)+\boldsymbol{P}_{i,j}(\tau)}{2} \approx \mu_{1,i - 0.5,j} p(x_{i - 0.5},y_{j},\tau),
\end{equation*}
and
\begin{equation*}
f_{d,i - 0.5,j}(P,\tau) = \frac{\tfrac{1}{2} \sigma^{2}_{1,i-1,j} \boldsymbol{P}_{i-1,j}(\tau) - \tfrac{1}{2} \sigma^{2}_{1,i,j} \boldsymbol{P}_{i,j}(\tau) }{\Delta x_{i}} \approx - \tfrac{\partial}{\partial x} \left( \tfrac{1}{2} \sigma^{2}_{1}(x,y_{j},\tau) p(x,y_{j},\tau) \right)\vert_{x = x_{i - 0.5}},
\end{equation*}
for $2 \le i \le m_{1}, 1 \le j \le m_{2}$.
Moreover
\begin{equation*}
f_{i,j \pm 0.5}(P,\tau) = f_{a,i,j \pm 0.5}(P,\tau) + f_{d,i,j \pm 0.5}(P,\tau),
\end{equation*}
where
\begin{equation*}
f_{a,i,j - 0.5}(P,\tau) = \mu_{2,i,j-0.5} \frac{\boldsymbol{P}_{i,j-1}(\tau)+\boldsymbol{P}_{i,j}(\tau)}{2} \approx \mu_{2,i,j-0.5} p(x_{i},y_{j-0.5},\tau),
\end{equation*}
and
\begin{equation*}
f_{d,i,j - 0.5}(P,\tau) = \frac{ \tfrac{1}{2} \sigma^{2}_{2,i,j-1} \boldsymbol{P}_{i,j-1}(\tau) - \tfrac{1}{2} \sigma^{2}_{2,i,j} \boldsymbol{P}_{i,j}(\tau) }{\Delta y_{j}} \approx - \tfrac{\partial}{\partial y} \left( \tfrac{1}{2} \sigma^{2}_{2}(x_{i},y,\tau) p(x_{i},y,\tau) \right)\vert_{y = y_{j - 0.5}},
\end{equation*}
for $1 \le i \le m_{1}, 2 \le j \le m_{2}$
Finally, for the mixed spatial derivative we define
\begin{eqnarray}
f_{m,i - 0.5,j - 0.5}(P,\tau) &=& \rho \sigma_{1,i-0.5,j-0.5}\sigma_{2,i-0.5,j-0.5} \frac{\boldsymbol{P}_{i-1,j- 1}(\tau)+ \boldsymbol{P}_{i-1,j}(\tau)+\boldsymbol{P}_{i,j- 1}(\tau) + \boldsymbol{P}_{i,j}(\tau)}{4} \nonumber \\
& \approx & \rho \sigma_{1}(x_{i-0.5},y_{j - 0.5},\tau)\sigma_{2}(x_{i-0.5},y_{j - 0.5},\tau)p(x_{i-0.5},y_{j - 0.5},\tau),
\label{eq:MixedDerivativeFluxDiscretization}
\end{eqnarray}
for $1 \le i \le m_{1}+1, 1 \le j \le m_{2}+1$, where it is assumed that
\begin{equation}
\boldsymbol{P}_{0,j}(\tau) := \boldsymbol{P}_{1,j}(\tau), \ \boldsymbol{P}_{m_{1}+1,j}(\tau):= \boldsymbol{P}_{m_{1},j}(\tau), \ \boldsymbol{P}_{i,0}(\tau) := \boldsymbol{P}_{i,1}(\tau), \ \boldsymbol{P}_{i,m_{2}+1}(\tau) := \boldsymbol{P}_{i,m_{2}}(\tau),
\label{eq:BoundaryMixedDerivative}
\end{equation}
such that the general formula is naturally extended at the boundaries of the spatial domain.
The numerical flux term \eqref{eq:MixedDerivativeFluxDiscretization} can be viewed as the result of applying the discretization for the advection part first in the $x$-direction and then in the $y$-direction.

The semidiscretization is completed by defining boundary conditions and discretizations at the boundaries of the truncated domain.
Since $p$ is again a density function, it follows that
$$ \int_{-\infty}^{\infty}\int_{-\infty}^{\infty} \left[ \tfrac{\partial}{\partial \tau} p \right] dx dy = 0. $$
Inserting the right hand side of the PDE \eqref{eq:GeneralForward}, this can be rewritten as
\begin{eqnarray*}
0 &=& \int_{-\infty}^{\infty} \left( \int_{-\infty}^{\infty} \left[ \tfrac{\partial^{2}}{\partial x^{2}} \left(\tfrac{1}{2} \sigma^{2}_{1}p \right) -  \tfrac{\partial}{\partial x} \left( \mu_{1} p \right) \right] dx \right) dy \\\\
&& + \ \int_{-\infty}^{\infty} \left( \int_{-\infty}^{\infty} \left[ \tfrac{\partial^{2}}{\partial y^{2}} \left( \tfrac{1}{2} \sigma_{2}^{2} p \right) - \tfrac{\partial}{\partial y} \left( \mu_{2} p \right) \right] dy \right) dx \\\\
&& + \ \int_{-\infty}^{\infty}\int_{-\infty}^{\infty} \tfrac{\partial^{2}}{\partial x \partial y} \left( \rho \sigma_{1}\sigma_{2} p \right) dx dy.
\end{eqnarray*}
Analogously to the one-dimensional case it is assumed that the boundaries are chosen sufficiently far away from the spot value $(X_{0},Y_{0})$ or that they are defined by a natural truncation of the spatial domain. 
The condition above is then approximated by
\begin{eqnarray}
0 &=& \int_{y_{\min}}^{y_{\max}} \left. \left[ \tfrac{\partial}{\partial x} \left(\tfrac{1}{2} \sigma^{2}_{1}p \right) - \mu_{1} p \right]\right\vert_{x=x_{\min}}^{x=x_{\max}} dy \nonumber \\ \nonumber \\
&& + \ \int_{x_{\min}}^{x_{\max}} \left. \left[ \tfrac{\partial}{\partial y} \left(\tfrac{1}{2} \sigma^{2}_{2}p \right) - \mu_{2} p \right] \right\vert_{y=y_{\min}}^{y=y_{\max}} dx \nonumber \\ \nonumber \\
&& + \ \int_{y_{\min}}^{y_{\max}}\int_{x_{\min}}^{x_{\max}} \tfrac{\partial^{2}}{\partial x \partial y} \left( \rho \sigma_{1}\sigma_{2} p \right) dx dy. \label{eq:BC2D}
\end{eqnarray}
Note that by assuming that
\begin{equation}
\rho \sigma_{1}\sigma_{2}p \vert_{x=x_{\min},y=y_{\min}} = \rho \sigma_{1}\sigma_{2}p \vert_{x=x_{\min},y=y_{\max}} = \rho \sigma_{1}\sigma_{2}p \vert_{x=x_{\max},y=y_{\min}} = \rho \sigma_{1}\sigma_{2}p \vert_{x=x_{\max},y=y_{\max}} =0,
\label{eq:HoekpuntenNul}
\end{equation}
the last integral, corresponding with the mixed derivative term, is always equal to zero.
Next, we generalise the idea that there are no fluxes at the boundaries, i.e.\ that no mass is coming in or going out at the boundaries.
In light of this it is assumed that the following boundary conditions hold:
\begin{eqnarray*}
\left. \left[ \tfrac{\partial}{\partial x} \left( \tfrac{1}{2} \sigma_{1}^{2}p \right) - \left( \mu_{1} p \right) \right]\right\vert_{x=x_{\min}} &=& 0, \qquad \mathrm{for} \ y_{\min} \le y \le y_{\max}, \\\\
\left. \left[ \tfrac{\partial}{\partial x} \left( \tfrac{1}{2} \sigma_{1}^{2}p \right) - \left( \mu_{1} p \right) \right]\right\vert_{x=x_{\max}} &=& 0, \qquad \mathrm{for} \ y_{\min} \le y \le y_{\max}, \\\\
\left. \left[ \tfrac{\partial}{\partial y} \left( \tfrac{1}{2} \sigma_{2}^{2}p \right) - \left( \mu_{2} p \right) \right]\right\vert_{y=y_{\min}} &=& 0, \qquad \mathrm{for} \ x_{\min} \le x \le x_{\max}, \\\\
\left. \left[ \tfrac{\partial}{\partial y} \left( \tfrac{1}{2} \sigma_{2}^{2}p \right) - \left( \mu_{2} p \right) \right]\right\vert_{y=y_{\max}} &=& 0, \qquad \mathrm{for} \ x_{\min} \le x \le x_{\max}.
\end{eqnarray*}
By combining this boundary conditions with the assumption \eqref{eq:HoekpuntenNul} it follows that condition \eqref{eq:BC2D} is satisfied.

For the discretization of the one-dimensional fluxes in \eqref{eq:FluxX} and \eqref{eq:FluxY} at the boundaries of the truncated domain, the approach from the one-dimensional case is generalised. By using a similar discretization of the boundary conditions, one ends up with numerical fluxes which are zero at the boundaries, i.e.\ with
\begin{equation*}
f_{0.5,j}(P,\tau) = f_{m_{1}+0.5,j}(P,\tau) = f_{i,0.5}(P,\tau) = f_{i,m_{2} + 0.5}(P,\tau) = 0, 
\end{equation*}
for $1 \leq i \le m_{1}, 1 \le j \le m_{2}$.
The fluxes stemming from the mixed derivative term, see \eqref{eq:FluxXY}, are discretized at the boundaries by using \eqref{eq:MixedDerivativeFluxDiscretization} in combination with \eqref{eq:BoundaryMixedDerivative}. 
It is readily verified that if
\begin{eqnarray*}
0 &=& \rho \sigma_{1,1,1} \sigma_{2,1,1} \boldsymbol{P}_{1,1}(\tau) = \rho \sigma_{1,1,m_{2}} \sigma_{2,1,m_{2}} \boldsymbol{P}_{1,m_{2}}(\tau) \\
  &=& \rho \sigma_{1,m_{1},1} \sigma_{2,m_{1},1} \boldsymbol{P}_{m_{1},1}(\tau) = \rho \sigma_{1,m_{1},m_{2}} \sigma_{2,m_{1},m_{2}} \boldsymbol{P}_{m_{1},m_{2}}(\tau),
\end{eqnarray*} 
which is the semidiscrete version of \eqref{eq:HoekpuntenNul}, then the total numerical flux over the boundary of the spatial domain is equal to zero and the total numerical mass is kept constant in time, i.e.\
\begin{equation*} 
\sum_{i=1}^{m_{1}} \sum_{j=1}^{m_{2}} \boldsymbol{P}_{i,j}(\tau) \vert \Omega_{i,j} \vert =\mathrm{constant,} \qquad \mathrm{for} \ \tau \ge 0. 
\end{equation*}

As stated above, some processes are naturally bounded, e.g.\ the general variance process from Section \ref{intro} can never become negative. 
Suppose for example that the process corresponding with the $y$-variable in the PDE \eqref{eq:GeneralForward} is bounded from below. Then, $y_{\min}$ is naturally taken equal to this lower boundary. Moreover, it can happen that this lower boundary is attainable (cf.\ the variance process from Section \ref{intro} with $\alpha < 1/2$) and probability mass can stack up at this boundary. This can cause instabilities in the approximation of the mixed derivative term near this boundary. In order to deal with this, if for example the boundary $y_{\min}$ is attainable, the central FV scheme in the $y$-direction for the "mixed derivative fluxes" \eqref{eq:FluxXY} at $y_{\min}+\tfrac{1}{2}\Delta y_{2}$ are replaced by a first-order forward scheme. More precisely, the $f_{m,i \pm 0.5, 1.5}(P,\tau)$ from above are then replaced by
$$ f_{m,i - 0.5, 1.5}(P,\tau) =  \rho \sigma_{1,i-0.5,1.5}\sigma_{2,i-0.5,1.5} \frac{\boldsymbol{P}_{i - 1,2}(\tau)+\boldsymbol{P}_{i,2}(\tau)}{2}, $$
for $1 \leq i \leq m_{1}+1$, where
$$ \boldsymbol{P}_{0,2}(\tau) := \boldsymbol{P}_{1,2}(\tau), \quad \boldsymbol{P}_{m_{1}+1,2}(\tau) := \boldsymbol{P}_{m_{1},2}(\tau). $$  

The total spatial discretization \eqref{eq:TotalDiscretization2D} yields a large system of differential equations. By making use of the \textit{Kronecker product}, this system of differential equations can be written as a system of ODEs
\begin{equation}
P'(\tau) = A(\tau) P(\tau),
\label{eq:ODE2D}
\end{equation} 
for $\tau > 0$ and with given matrix $A(\tau)$. 
Analogously to the one-dimensional case, since the values $\boldsymbol{P}_{i,j}(\tau)$ can be seen as approximations of the cell averages $\overline{p}(x_{i},y_{j},\tau)$, it is natural to define the initial vector as $P(0) = \mathrm{vec}[\boldsymbol{P}(0)]$ where
$$ \boldsymbol{P}_{i,j}(0) = \left\{ \begin{array}{lll}
\tfrac{1}{\vert \Omega_{i,j} \vert} & & \mathrm{if} \ (X_{0},V_{0}) \in \Omega_{i,j}, \\ \\
0 & & \mathrm{else.}
\end{array} \right. $$
Note that the matrix $A$ can be split as
$$ A = A_{0} + A_{1} + A_{2},$$
where $A_{1}$, respectively $A_{2}$, represents the discretization of the spatial derivatives in the first, respectively second, spatial dimension. The matrix $A_{0}$ represents the discretization of the mixed spatial derivative term in \eqref{eq:GeneralForward}. It is readily verified that $A_{1}$ is tridiagonal, $A_{2}$ is essentially tridiagonal and $A_{0}$ has in general nine non-zero elements per row and column. 
In Section \ref{ADI} this structure is used to perform an effective time discretization of the general system of ODEs \eqref{eq:ODE2D}.

\subsection{Numerical experiment}

In this section the performance of the 2D FV discretization is tested for two practical examples. For the first experiment we consider the two-dimensional Black--Scholes model which can be described by the system of SDEs
\begin{equation*}
\left\{ \begin{array}{l}
dS_{1,\tau} = r S_{1,\tau} d\tau + \sigma_{1,BS} S_{1,\tau} dW^{(1)}_{\tau}, \\\\
dS_{2,\tau} = r S_{2,\tau} d\tau + \sigma_{2,BS} S_{2,\tau} dW^{(2)}_{\tau},
\end{array} \right.
\end{equation*}
with $dW^{(1)}_{\tau} \cdot dW^{(2)}_{\tau} = \rho d\tau$, $-1 \le \rho \le 1$ and $r, \sigma_{1,BS},\sigma_{2,BS}$ strictly positive constants. The corresponding forward Kolmogorov equation is given by
\begin{equation*}
\tfrac{\partial}{\partial \tau} p = \tfrac{\partial^{2}}{\partial s_{1}^{2}} \left( \tfrac{1}{2} \sigma_{1,BS}^{2}s_{1}^{2}p \right) + \tfrac{\partial^{2}}{\partial s_{1} \partial s_{2}} \left( \rho \sigma_{1,BS} \sigma_{2,BS} s_{1} s_{2} p \right) + \tfrac{\partial^{2}}{\partial s_{1}^{2}} \left( \tfrac{1}{2} \sigma_{2,BS}^{2} s_{2}^{2} p \right) - \tfrac{\partial}{\partial s_{1}} \left( r s_{1} p \right) - \tfrac{\partial}{\partial s_{2}} \left( r s_{2} p \right),
\end{equation*}
for $s_{1}, s_{2}, \tau > 0$ and with $p(s_{1},s_{2},0)=\delta(s_{1}-S_{1,0})\delta(s_{2}-S_{2,0})$ for given values $S_{1,0}, S_{2,0}$.
The exact solution is known analytically and can be written as
\begin{equation*}
p(s_{1},s_{2},\tau) = n_{2}(\log(s_{1}/S_{1,0}),\log(s_{2}/S_{2,0}),\tau)\frac{1}{s_{1}}\frac{1}{s_{2}}, \qquad \mathrm{for} \ s_{1}>0, s_{2}>0, \tau > 0,
\end{equation*}
where this time $n_{2}(x,y,\tau)$ is the density function of a two-dimensional normally distributed random variable with mean $\boldsymbol{\mu}$ and covariance matrix $\boldsymbol{\Sigma}$ given by
$$ \boldsymbol{\mu} = \left[ \begin{array}{c}
(r - \tfrac{1}{2} \sigma^{2}_{1,BS})\tau \\\\ (r - \tfrac{1}{2} \sigma^{2}_{2,BS})\tau
\end{array} \right], \qquad
\boldsymbol{\Sigma} = \left[ \begin{array}{ccc}
\sigma^{2}_{1,BS}\tau && \rho \sigma_{1,BS}\sigma_{2,BS}\tau  \\\\
\rho \sigma_{1,BS}\sigma_{2,BS}\tau && \sigma^{2}_{2,BS}\tau
\end{array} \right]. $$

\begin{figure}
\begin{center}
\includegraphics[scale=0.5]{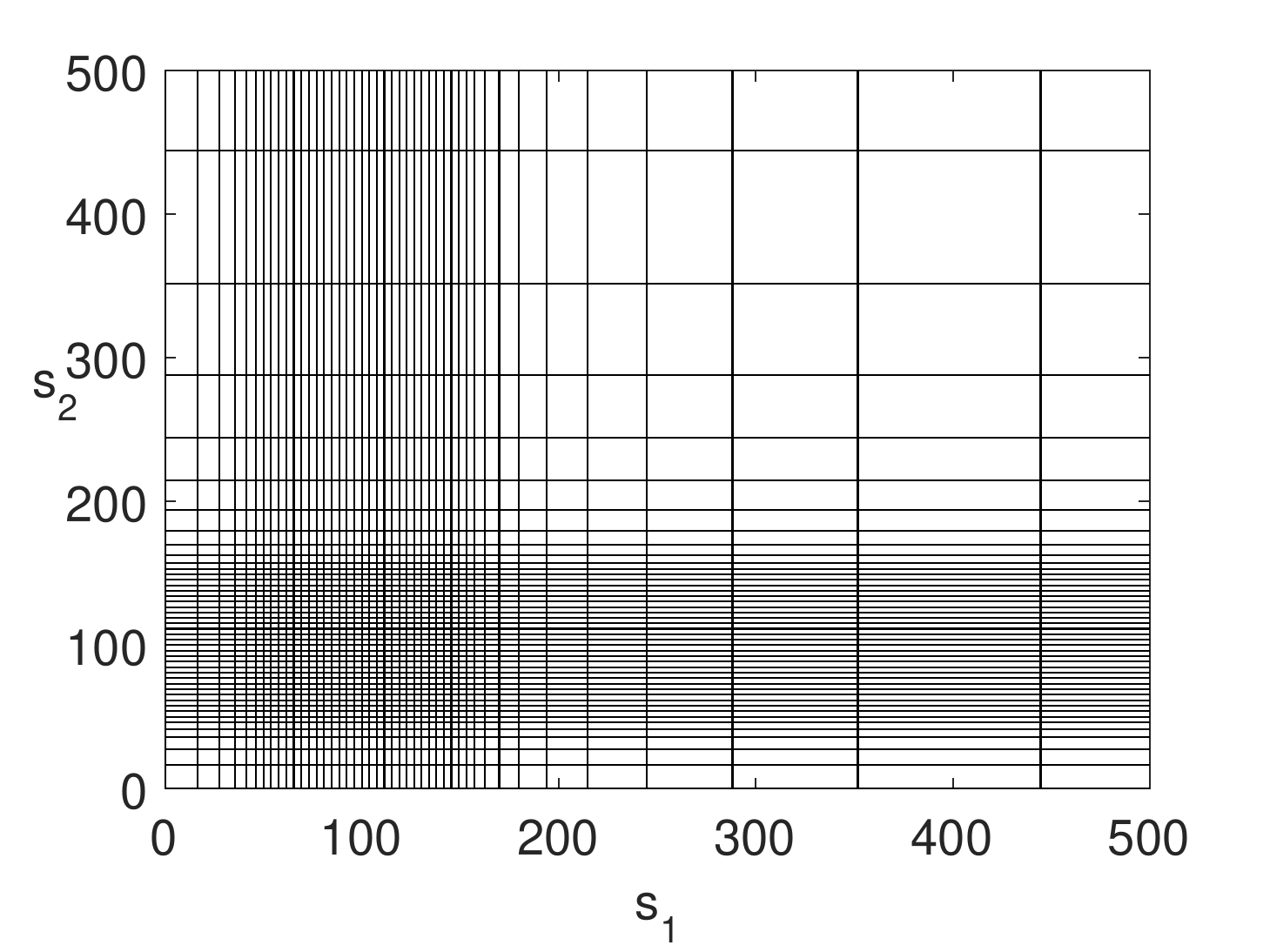}  
\caption{Illustration of the non-uniform grid around $(S_{1,0},S_{2,0})$ for the Black--Scholes example and the actual values $S_{1,0}=S_{2,0}=100$, $m_{1}=m_{2}=50$.}
\label{fig:Grid2DBS}
\end{center}
\end{figure}

\begin{figure}
\begin{center}
\includegraphics[scale=0.5]{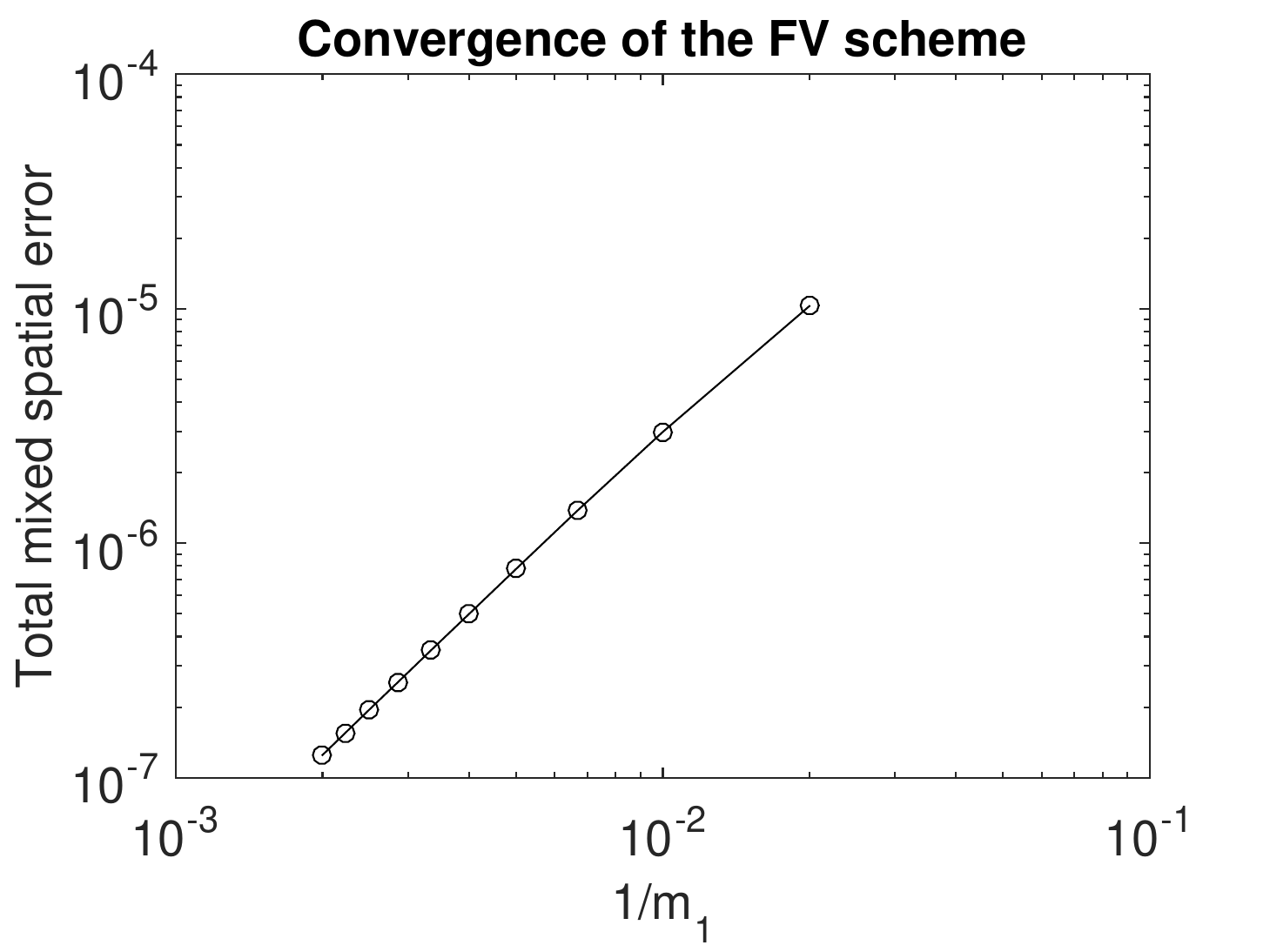}  
\includegraphics[scale=0.5]{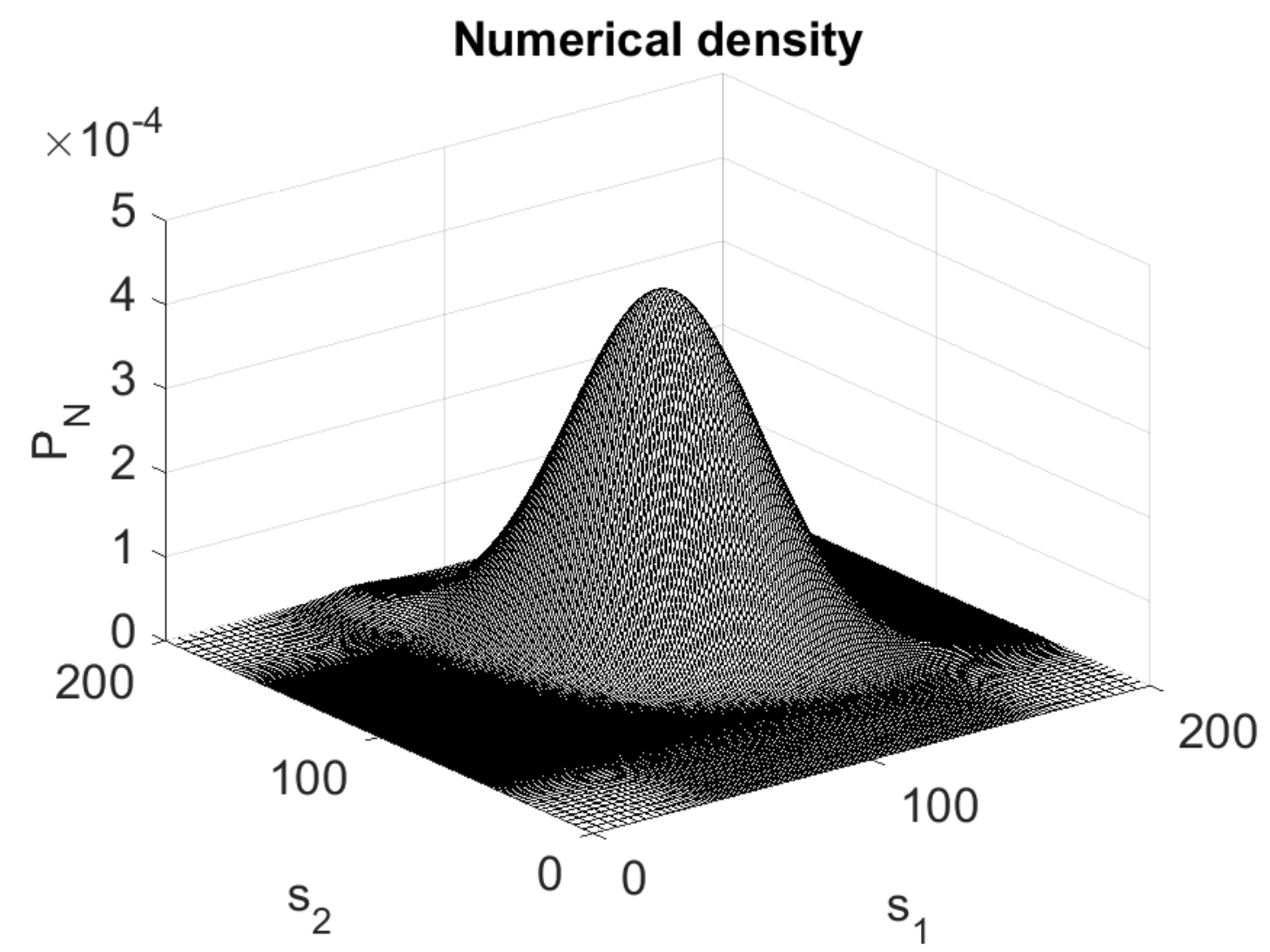} 
\caption{Convergence results within the 2D Black--Scholes model. The parameter values are $r = 0.03, \sigma_{1,BS} = 0.2, \sigma_{2,BS}=0.25, \rho=-0.7, T=1$.}
\label{fig:Convergence2DBS}
\end{center}
\end{figure}

Similarly to the domain truncation in the one-dimensional numerical experiment from Subsection \ref{subsec:Experiment1D}, the spatial domain is truncated to $[S_{1,\min}, S_{1,\max}] \times [S_{2,\min}, S_{2,\max}] = [0, 30S_{1,0}] \times [0, 30S_{2,0}]$. The Cartesian spatial grid,
$$ (s_{1,i},s_{2,j}) \qquad \mathrm{for} \ 1 \le i \le m_{1}, 1 \le j \le m_{2}, $$ 
is constructed by considering the spatial grid from Subsection \ref{subsec:Experiment1D} in both spatial dimension. In Figure \ref{fig:Grid2DBS} this spatial grid is illustrated within the region  $[0, 5S_{1,0}] \times [0, 5S_{2,0}]$ for $S_{1,0} = S_{1,0} = 100$ and $m_{1}= m_{2} = 50$.
The FV discretization from Subsection \ref{2DKolmogorov} then defines approximations $\boldsymbol{P}_{i,j}(\tau)$ to $p(s_{1,i},s_{2,j},\tau)$ and we compute the total mixed spatial error
$$ \max_{1 \le i \le m_{1}, 1 \le j \le m_{2}} \epsilon_{i,j}(m), $$
where
\begin{equation*}
\epsilon_{i,j}(m) = \left\{ \begin{array}{ll}
\left\vert \tfrac{p(s_{1,i},s_{2,j},T) - \boldsymbol{P}_{i,j}(T)}{p(s_{1,i},s_{2,j},T)} \right\vert \qquad & \mathrm{if} \ p(s_{1,i},s_{2,j},T) > 1, \\\\
\vert p(s_{1,i},s_{2,j},T) - \boldsymbol{P}_{i,j}(T) \vert \qquad & \mathrm{else}. 
\end{array} \right.
\end{equation*}
The values $\boldsymbol{P}_{i,j}(T)$ are calculated by applying the Hundsdorfer--Verwer time stepping method (see Section \ref{ADI}) with a large number of time steps such that the temporal discretization error is negligible.
In Figure \ref{fig:Convergence2DBS} the total mixed spatial error is shown for the realistic parameter values $r= 0.03$, $\sigma_{1,BS}=0.2$, $\sigma_{2,BS}=0.25$, $\rho=-0.7$, $T=1$ and for the number of spatial grid points $m_{1}=m_{2} = \{50,100, \ldots, 500 \}$. 
The convergence plot indicates that the FV discretization is second-order convergent with respect to the current initial-boundary value problem.

For the second example we consider the popular Heston model \cite{H93}, i.e.\ the SV model \eqref{eq:SVmodel} with $\psi(v)=\sqrt{v}$ and $\alpha=1/2$. The underlying density function satisfies the forward Kolmogorov equation
\begin{equation}
\tfrac{\partial}{\partial \tau} p = \tfrac{\partial^{2}}{\partial x^{2}} \left( \tfrac{1}{2} vp \right) + \tfrac{\partial^{2}}{\partial x \partial v} \left( \rho \xi v p \right) + \tfrac{\partial^{2}}{\partial v^{2}} \left( \tfrac{1}{2} \xi^{2} v p \right)
 - \tfrac{\partial}{\partial x} \left( (r_{d}-r_{f}-\tfrac{1}{2}v) p \right) - \tfrac{\partial}{\partial v} \left( \kappa(\eta - v) p \right),
\label{eq:HestonPDEForward}
\end{equation}
for $x \in \mathbb{R}, v>0, \tau >0$ and with initial condition $p(x,v,0) = \delta(x)\delta(v-V_{0})$ where $X_{0}=0$ and $V_{0}$ is given.
To the best of our knowledge, no analytical solution is available in the literature for the density function $p$ that satisfies \eqref{eq:HestonPDEForward}. In order to test the performance of our FV discretization, we compute a reference solution with an alternative discretization method described in \cite{FO11}. The latter method is based on rewriting
\begin{equation}
p(x,v,\tau) = p_{1}(x,\tau \vert V_{SV,\tau}=v)p_{2}(v,\tau), 
\label{eq:COSQuadrature}
\end{equation}
where $p_{2}(v,\tau)$ denotes the one-dimensional density of the volatility process, and $p_{1}(x,\tau \vert V_{SV,\tau}=v)$ denotes the conditional density of $X_{\tau}$ given the variance value $V_{SV,\tau} = v$. 
In \cite{FO11} the characteristic function 
$$\psi(\omega\vert V_{SV,\tau} =v) = \mathbb{E}[\exp(\imi \omega X_{SV,\tau}) \vert V_{SV,\tau} = v]$$
corresponding with $p_{1}(x,\tau \vert V_{SV,\tau}=v)$ is given in semi-analytical form and it is stated that $p_{2}(v,\tau)$ is given by \eqref{eq:CIRdensity}.
By using the COS-method \cite{FO08} we approximate the conditional density function $p_{1}(x,\tau \vert V_{SV,\tau}=v)$ and our reference solution $p_{ref}$ is then defined via \eqref{eq:COSQuadrature}.

\begin{figure}
\begin{center}
\includegraphics[scale=0.5]{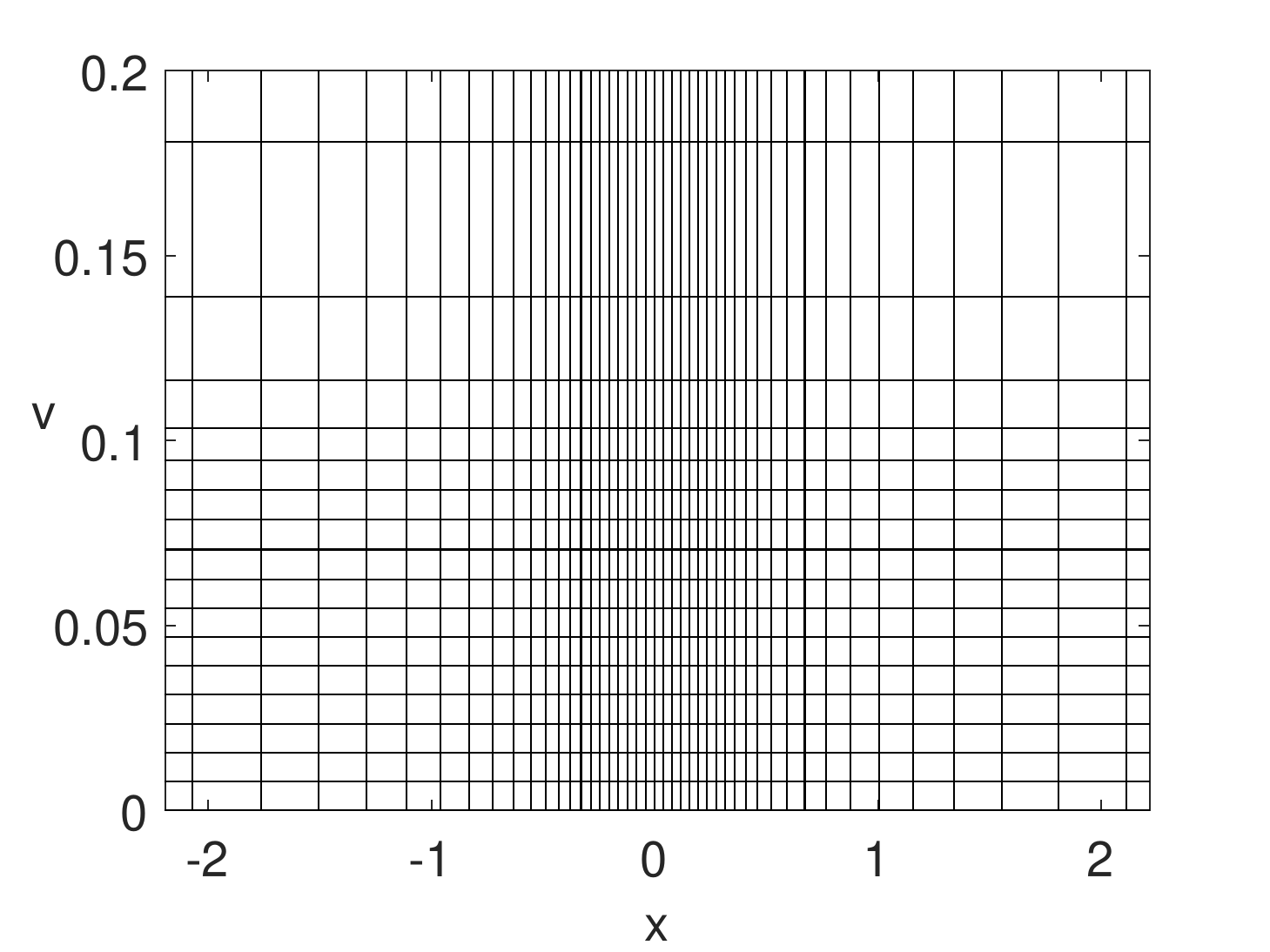}  
\caption{Illustration of the non-uniform grid around $(X_{0},V_{0})$ for the Heston example and the actual values $X_{0}=0, V_{0}=0.0625$, $m_{1}=2m_{2}=50$.}
\label{fig:Grid2DHeston}
\end{center}
\end{figure}

In the numerical experiment, we opt to truncate the domain in the $x$-direction to $[X_{0} - \log(30), X_{0} + \log(30)]$. The spatial domain in the $v$-direction is truncated to $[0, 15]$, analogously to the CIR example.
We consider spatial grids 
$$ -\log(30) = x_{1} < x_{2} < \cdots < x_{m_{1}} = \log(30), $$
$$ 0 = v_{1} < v_{2} < \cdots < v_{m_{2}} = 15, $$ 
with $m_{1}=2m_{2}$, which are similar to the ones described in \cite{HH12} and such that there exist indices $i_{0},j_{0}$ such that $(x_{i_{0}},v_{j_{0}}) = (X_{0},V_{0})$.
In Figure \ref{fig:Grid2DHeston} the spatial grid is shown for $X_{0}=0, V_{0} = 0.0625$ and the sample value $m_{1}=2m_{2}=50$.
Applying the FV discretization then yields approximations $\boldsymbol{P}_{i,j}(\tau)$ to $p(x_{i},v_{j},\tau)$.
From equation \eqref{eq:COSQuadrature} and the CIR example it follows that if $q= \tfrac{2\kappa\eta}{\xi^{2}}-1<0$, then the density function can tend to infinity as $v$ tends to zero and at $v=0$ the exact density is not defined.
Analogously to the remark in Subsection \ref{subsec:Experiment1D}, it is readily seen that adequately comparing the difference between $p(x_{i},v_{2},T)$ and $P_{i,2}(T)$ then becomes difficult for increasing values of $m_{2}$.
The error is therefore computed on similar spatial domains.
Let $v_{low}$ again be the smallest non-zero grid point in the $v$-direction if the total number of grid points in that direction is $m_{2}=50$. 
For given $m_{2}$, let $j_{1}$ be the lowest index such that $v_{j_{1}} \ge v_{low}$ and define the total mixed spatial error by
$$ \max_{1 \le i \le m_{1}, j_{1} \le j \le m_{2}} \epsilon_{i,j}(m), $$
where
\begin{equation*}
\epsilon_{i,j}(m) = \left\{ \begin{array}{ll}
\left\vert \tfrac{p_{ref}(x_{i},v_{j},T) - \boldsymbol{P}_{i,j}(T)}{p_{ref}(x_{i},v_{j},T)} \right\vert \qquad & \mathrm{if} \ p_{ref}(x_{i},v_{j},T) > 1, \\\\
\vert p_{ref}(x_{i},v_{j},T) - \boldsymbol{P}_{i,j}(T) \vert \qquad & \mathrm{else.}
\end{array} \right.
\end{equation*}
The values $P_{i,j}(T)$ are once more approximated by applying the Hundsdorfer--Verwer time stepping method (see Section \ref{ADI}) with a small time step such that the temporal discretization error is negligible. 

\begin{figure}
\begin{center}
\includegraphics[scale=0.5]{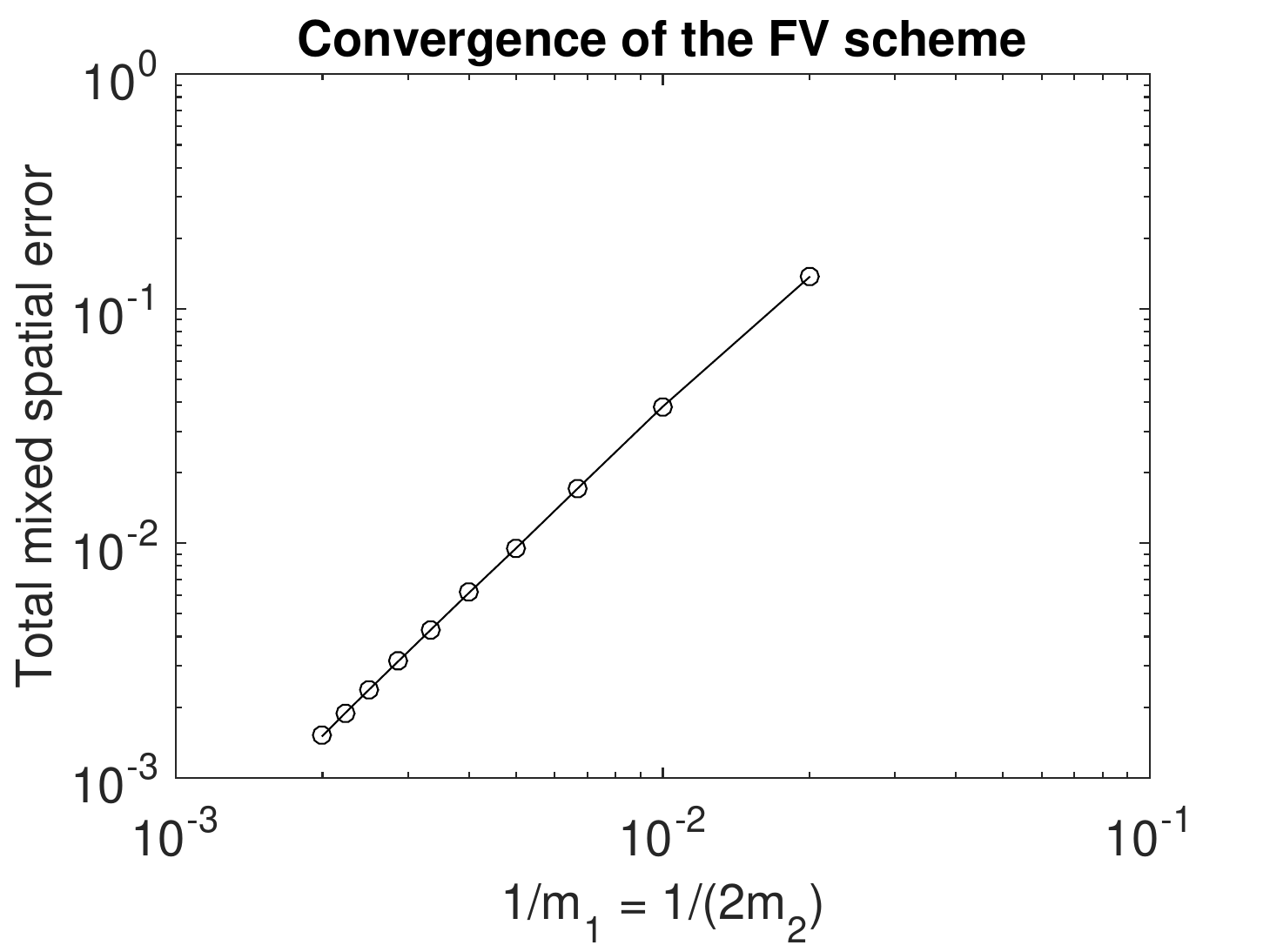}  
\includegraphics[scale=0.5]{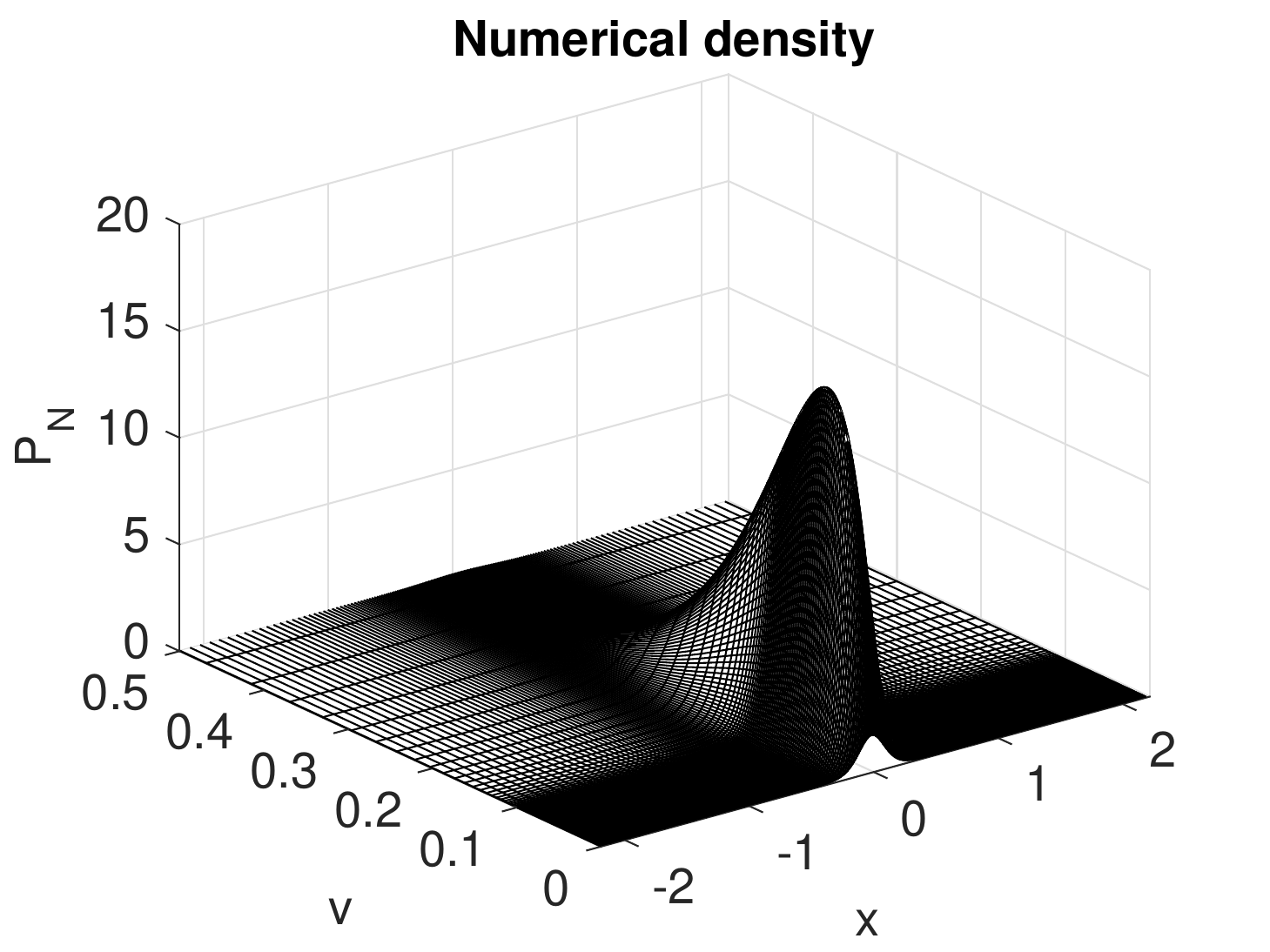}
\caption{Convergence results within the Heston model. The parameters are given by Set C.}
\label{fig:Convergence2D_HestonSetA}
\end{center}
\end{figure}

\begin{figure}
\begin{center}
\includegraphics[scale=0.5]{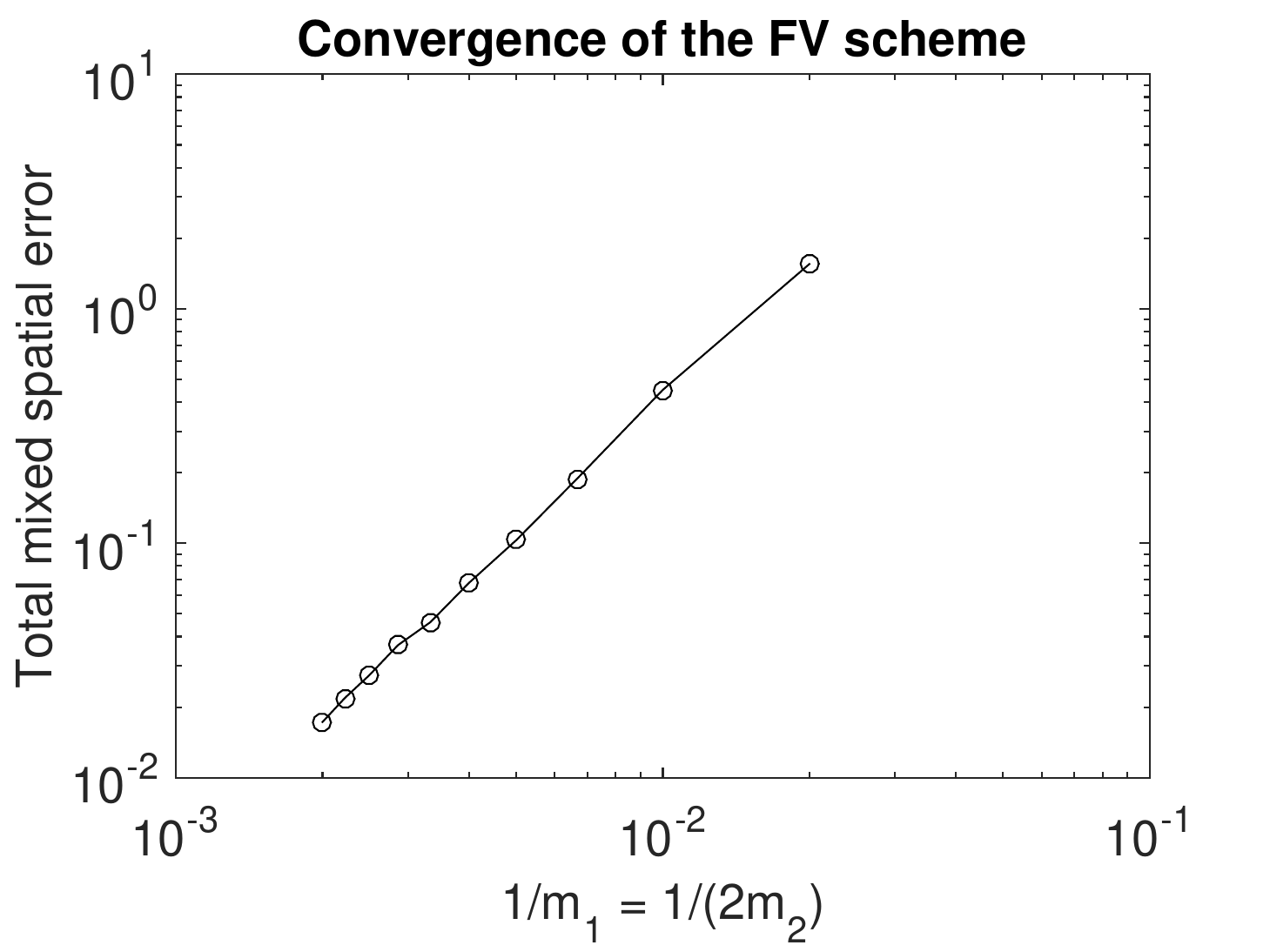}  
\includegraphics[scale=0.5]{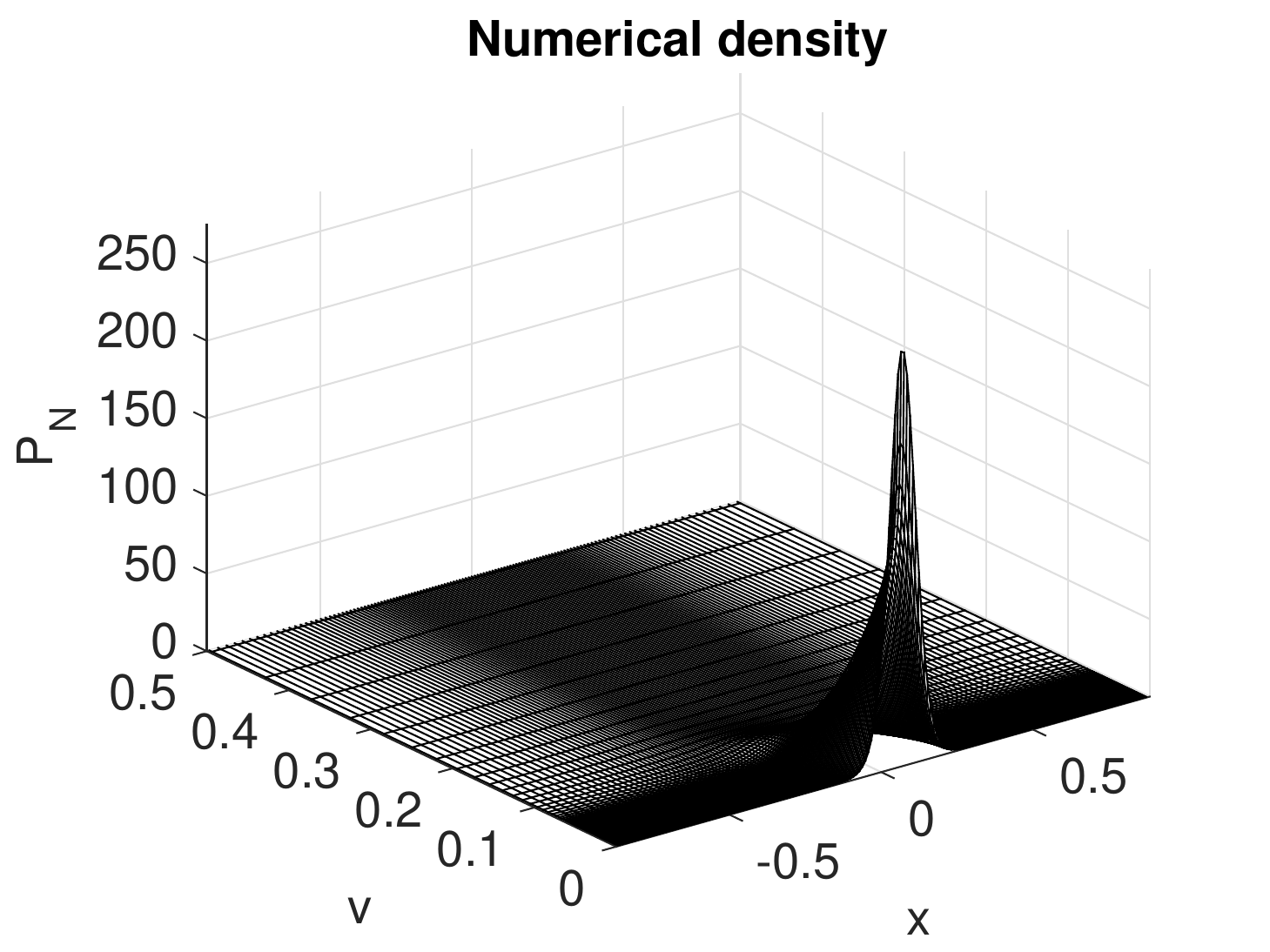}
\caption{Convergence results within the Heston model. The parameters are given by Set D.}
\label{fig:Convergence2D_HestonSetB}
\end{center}
\end{figure}

For the actual numerical experiments we consider an extension of the two sets of parameters used in Subsection \ref{subsec:Experiment1D}. The extensions are also taken from \cite{FO11}, used in \cite{RO12}, and given by
\begin{center}
\begin{tabular}{|c|c|c|c|c|c|c|c|c|c|}
\hline
& $\kappa$ & $\eta$ & $\xi$ & $\rho$ & $r_{d}$ & $r_{f}$ & $X_{0}$ & $V_{0}$ & $T$  \\
\hline
Set C & 5 & 0.16 & 0.9 & 0.1 & 0.1 & 0 & 0 & 0.0625 & 0.25 \\
\hline
Set D & 1.15 & 0.0348 & 0.39 & -0.64 & 0.04 & 0 & 0 & 0.0348 & 0.25 \\
\hline
\end{tabular}
\end{center}
Recall that for Set C we have that $q=0.98$ and for Set D it holds that $q=-0.47$.
In the left plot of Figure \ref{fig:Convergence2D_HestonSetA}, respectively Figure \ref{fig:Convergence2D_HestonSetB}, the total mixed spatial error is shown for the parameters of Set C, respectively Set D, and for the number of spatial grid points $m_{1} = 2m_{2} = \{ 50, 100, \ldots, 500\}$.
In the right plots, the corresponding numerical solutions are shown for $m_{1}=2m_{2}=200$.
The convergence plots indicate that the FV discretization is convergent with respect to the current initial-boundary value problems. Additional experiments again suggest that the FV discretization is second-order convergent if the Feller condition is satisfied. If $q<0$ the order of convergence can drop to one.
Please note that the conclusions of the numerical experiments are essentially unchanged for different values of $v_{low}$ as long as it is defined via one of the coarsest grids considered in the experiment.
The two-dimensional tests also confirm that the total numerical mass is, indeed, kept constantly equal to one, even if the Feller condition is strongly violated.

\setcounter{equation}{0}
\section{Temporal discretization by ADI methods}\label{ADI}

In general, spatial discretization by FV methods of initial-boundary value problems for time-dependent convection diffusion equations leads to large systems of stiff ODEs of the type  
$$ W'(\tau) = F(\tau, W(\tau)) \quad (0\leq \tau \leq T), \quad W(0) = W_{0}, $$
with given vector-valued function $F:[0,T]\times \R^{m} \rightarrow \R^{m}$ and given vector $W_{0} \in \R^{m}$ where $m$ is the number of volumes.
When the function $F$ is stemming from semidiscretization of a multidimensional PDE, classical implicit time stepping schemes such as the \textit{Crank--Nicolson scheme} can often be very time consuming. For the effective time discretization of such semidiscrete systems, operator splitting schemes of the Alternating Direction Implicit (ADI) type are widely considered.
In the past decades, several ADI schemes have been studied for the situation where mixed spatial derivatives are present in the convection-diffusion equation. Mixed derivative terms are ubiquitous, notably, in the field of financial mathematics. There they arise due to correlations between the underlying stochastic processes, cf.\ the general SLV model \eqref{eq:SLVmodel}.

In this paper we consider the state-of-the-art \textit{Hundsdorfer--Verwer} (HV) \textit{scheme}, \cite{H02,V99}, which is of the ADI-type.
Suppose the semidiscrete operator $F$ is stemming from the semidiscretization of a two-dimensional convection-diffusion equation and suppose $F$ can be decomposed as
$$ F(\tau,w) = F_{0}(\tau,w) + F_{1}(\tau,w) + F_{2}(\tau,w) \qquad (0 \le \tau \le T, w \in \R^{m}), $$
where $F_{0}$ represents the mixed spatial derivative term and $F_{1}$, respectively $F_{2}$, represents all spatial derivative terms in the first, respectively second, spatial direction.
Let $\theta > 0$ be a given parameter, $N \ge 1$ the number of time steps and $\tau_{n} = n \Delta \tau$ with $\Delta \tau = T/N$.
Then the HV scheme defines approximation $W_{n}$ to $W(\tau_{n})$ successively for $n=1,2,3,\ldots,N$ through
\begin{equation}
\label{eq:HV}
\left\{\begin{array}{l}
Y_0 = W_{n-1}+\Delta \tau\, F(\tau_{n-1},W_{n-1}), \\\\
Y_l = Y_{l-1}+\theta\Delta \tau \left(F_l(\tau_n,Y_l)-F_l(\tau_{n-1},W_{n-1})\right),
\quad l=1,2, \\\\
\widetilde{Y}_0 = Y_{0} + \tfrac{1}{2} \Delta \tau \left( F(\tau_{n},Y_{2}) - F(\tau_{n-1},W_{n-1})\right), \\\\
\widetilde{Y}_l = \widetilde{Y}_{l-1}+\theta \Delta \tau \,(F_l(\tau_n,\widetilde{Y}_l)-F_{l}(\tau_{n},Y_{2})), \quad l=1,2, \\\\
W_n = \widetilde{Y}_2.
\end{array}\right.
\end{equation}
The HV scheme \eqref{eq:HV} starts with an explicit Euler predictor stage, followed by two implicit but unidirectional corrector stages. Then a second explicit stage is performed, followed again by two implicit unidirectional corrector stages. The operator $F_{0}$, which contains the mixed spatial derivative term, is always treated in an explicit manner. The implicit stages only handle spatial derivatives in only one spatial dimension. This can lead to a major computational advantage in comparison to classical non-splitted implicit time stepping methods.

Recently, various positive stability results for the HV scheme have been derived relevant to multidimensional convection-diffusion equations with mixed derivative term, see e.g.\ \cite{IHM13,IHW07,IHW09}. Subsequently, it has been proved \cite{IHW15} that, under some natural stability and smoothness assumptions, the HV scheme is second-order convergent with respect to the time step whenever it is applied to semidiscrete two-dimensional convection-diffusion equations with mixed derivative term. The temporal convergence result from \cite{IHW15} has the key property that it holds uniformly in the spatial mesh width. 

In this article, the vector $P(0)$ is stemming from an initial function that is nonsmooth. It is well-known that convergence of time discretization methods can then be seriously impaired, cf.\ \cite{PVF03,W16}. A convergence analysis for the HV scheme relevant to nonsmooth data is still open in the literature. 
In order to prevent the numerical solution from undesirable behaviour, Rannacher time stepping will be applied, that is, the first few time steps of the HV scheme are replaced by twice as many half time steps of the implicit Euler scheme \cite{R84}. 
We opt to replace here the first two HV time steps by four half time steps of the implicit Euler scheme.
This choice is based on the results in \cite{W16} for the Modified Craig--Sneyd scheme which forms another popular ADI scheme for multidimensional convection-diffusion equations with mixed derivative terms.
Note that applying Rannacher time stepping to systems of ODEs stemming from multidimensional convection-diffusion equations can be time consuming if exact solvers are used. In order to limit this drawback, the implicit Euler time steps are performed with a suitable multigrid solver, \cite{HY02}.

\setcounter{equation}{0}
\section{Calibration of the SLV model to the LV model}\label{Calibration}

As stated in Section \ref{intro}, the goal of the paper is to calibrate state-of-the-art SLV models to its underlying LV model in order to reproduce the known prices for European call and put options. This is done by defining the leverage function $\sigma_{SLV}$ such that the relationship \eqref{eq:SigmatoMatch} is satisfied. By combining equations \eqref{eq:SigmatoMatch}, \eqref{eq:CondExpec} and \eqref{eq:ForwardKolmogorov}, it is readily seen that a highly non-linear PDE needs to be solved. 
In this section the FV-ADI method is used in combination with an inner iteration to approximate the corresponding leverage function and density function.

In order to use the FV-ADI discretization, one first has to define spatial and temporal grids. Since the initial function of forward Kolmogorov equations is highly non-smooth around the spot value $(X_{0},V_{0})$, and the region of interest is also situated around this value, it is natural to consider non-uniform Cartesian grids that are concentrated around the value $(X_{0},V_{0})$. If the parameter $\alpha$ from the SLV model is chosen smaller than or equal to $1/2$, then the natural boundary $V_{\tau}=0$ can be reached and probability mass stacks up at the boundary $v=0$. It is then natural to additionally require smaller mesh widths in the $v$-direction at this boundary.
The non-uniform grids define volumes of which the volume average is approximated by the FV scheme.
Denote by $m_{1}$, respectively $m_{2}$, the number of spatial grid points in the $x$-direction, respectively $v$-direction. 
We consider spatial grids 
$$ x_{\min} = x_{1} < x_{2} < \cdots < x_{m_{1}} = x_{\max}, $$
$$ 0 = v_{1} < v_{2} < \cdots < v_{m_{2}} = V_{\max}, $$ 
which are similar to the ones described in \cite{HH12} and such that there exist indices $i_{0},j_{0}$ such that $(x_{i_{0}},v_{j_{0}}) = (X_{0},V_{0})$. Denote their corresponding mesh widths by $\Delta x_{i}, \Delta v_{j}$ and define volumes
$$ \Omega_{i,j} = [x_{i-0.5}, x_{i+0.5}] \times [v_{j-0.5}, v_{j+0.5}],$$
where the values $x_{i-0.5}, v_{j-0.5}$ are defined similar as in Section \ref{FV}.
The spatial grids are the result of a continuous transformation of a uniform underlying grid and it can be shown that the pertinent meshes are smooth.
The values $x_{\min}, x_{\max}, V_{\max}$ are chosen sufficiently far away from the spot value such that the boundary conditions from Section \ref{FV} can be applied. An example of the spatial grid for the small sample values $m_{1} = 2m_{2} = 50$ was already shown in Figure \ref{fig:Grid2DHeston} for the case where $(X_{0},V_{0}) = (0, 0.0348)$.
For the discretization in time we always consider uniform grids $\tau_{n} = n \Delta \tau$ where the step size is given by $\Delta \tau = T / N$ and $N$ denotes the total number of time steps.

Once the spatial grid and corresponding volumes are defined, the FV discretization from Section \ref{FV} is applied. This yields a large system of ordinary differential equations 
\begin{equation}
P'(\tau) = A(\tau) P(\tau) = (A_{0}(\tau) + A_{1}(\tau) + A_{2}(\tau))P(\tau) \quad (0 \le \tau \le T),
\label{eq:EffectiveSemidiscrSystem}
\end{equation}
with given matrices $A_{0}(\tau), A_{1}(\tau), A_{2}(\tau)$ and initial function defined via
\begin{equation*}
\boldsymbol{P}_{i,j}(0) = 
\left\{ \begin{array}{ll}
\tfrac{2}{\Delta x_{i-1} + \Delta x_{i}}\tfrac{2}{\Delta v_{j-1} + \Delta v_{j}} & \qquad \mathrm{if} \ (i,j) = (i_{0},j_{0}), \\\\
0 & \qquad \mathrm{else}.
\end{array} \right.
\end{equation*}
The matrices $A_{0}, A_{1}$ contain, however, the unknown function $\sigma_{SLV}$ at the spatial grid points. 
The $P_{i}(\tau)$ can be used in combination with a numerical integration technique to approximate the conditional expectations \eqref{eq:CondExpec} and hence the pertinent leverage function.
We opt to perform the numerical integration with the trapezoidal rule and define approximations
\begin{equation}
\mathbb{E}_{i}(\tau) = \frac{\sum_{j=1}^{m_{2}} \psi^{2}(v_{j}) \boldsymbol{P}_{i,j}(\tau) \tfrac{\Delta v_{j} + \Delta v_{j+1}}{2}}{\sum_{j=1}^{m_{2}} \boldsymbol{P}_{i,j}(\tau) \tfrac{\Delta v_{j} + \Delta v_{j+1}}{2}} \approx \mathbb{E}[ \psi^{2}(V_{\tau}) \Vert X_{\tau} = x_{i} ],
\label{eq:SemidiscreteCondExpec}
\end{equation}
where we recall that $\Delta v_{1} = \Delta v_{m_{2}+1} = 0$.
Inserting the approximations \eqref{eq:SemidiscreteCondExpec} into \eqref{eq:EffectiveSemidiscrSystem} leads to a non-linear system of ODEs.

As a final step, the system of ODEs \eqref{eq:EffectiveSemidiscrSystem} is discretized in time with the time stepping scheme from Section \ref{ADI} and an inner iteration to handle the non-linearity, cf.\ \cite{TF10}.
By applying the HV scheme, the conditional expectations \eqref{eq:SemidiscreteCondExpec} are naturally replaced by their fully discrete versions
\begin{equation}
\mathbb{E}_{n,i} = \frac{\sum_{j=1}^{m_{2}} \psi^{2}(v_{j}) \boldsymbol{P}_{n,i,j} \tfrac{\Delta v_{j} + \Delta v_{j+1}}{2}}{\sum_{j=1}^{m_{2}} \boldsymbol{P}_{n,i,j} \tfrac{\Delta v_{j} + \Delta v_{j+1}}{2}},
\label{eq:DiscreteCondExpec1}
\end{equation}
and we define the leverage function $\sigma_{SLV}$ at the spatial and temporal grid by
\begin{equation}
\sigma_{SLV}(x_{i},\tau_{n}) = \frac{\sigma_{LV}(x_{i},\tau_{n})}{\sqrt{\mathbb{E}_{n,i}}}.
\label{eq:FullyDiscreteLeverage}
\end{equation}
It is readily seen that at the initial time level, i.e.\ at $n=0$, the expression $\eqref{eq:DiscreteCondExpec1}$ is only defined if $i = i_{0}$. To render the calibration procedure feasible we put
\begin{equation*}
\mathbb{E}_{0,i} = \psi^{2}(V_{0}) \quad \mathrm{for} \ 1 \le i \le m_{1}.
\end{equation*}
For strictly positive time levels, i.e.\ for $n > 0$, the exact density $p(x_{i},v_{j},\tau_{n})$ is always non-negative. By performing the spatial and temporal discretization, however, it is possible that some of the values $\boldsymbol{P}_{n,i,j}$ become (slightly) negative. In order to prevent the numerical solution from undesirable behaviour, the expression \eqref{eq:DiscreteCondExpec1} is replaced in the calibration procedure by
\begin{equation}
\mathbb{E}_{n,i} = \frac{\sum_{j=1}^{m_{2}} \psi^{2}(v_{j}) \vert \boldsymbol{P}_{n,i,j} \vert \tfrac{\Delta v_{j} + \Delta v_{j+1}}{2}}{\sum_{j=1}^{m_{2}} \vert \boldsymbol{P}_{n,i,j} \vert \tfrac{\Delta v_{j} + \Delta v_{j+1}}{2}} \quad \mathrm{for} \ 1 \le i \le m_{1}, \ n>0.
\label{eq:DiscreteCondExpec}
\end{equation}
Theoretically it is possible that the denominator (and hence also the nominator) of \eqref{eq:DiscreteCondExpec} equals zero and the fully discrete conditional expectation is undefined. In this case we assume that the conditional expectation is locally constant in time and set $\mathbb{E}_{n,i} = \mathbb{E}_{n-1,i}$.
Let $Q \ge 1$ be a given integer. For the actual calibration of the SLV model to the LV model, we employ the following numerical procedure.
\newline
\newline
\texttt{for $n$ is $1$ to $N$ do
\begin{itemize}
\item[] let $P_{n} = P_{n-1}$ be an initial approximation to $P(\tau_{n})$; \newline \newline 
for $q$ is $1$ to $Q$ do
\begin{itemize}
\item[(a)] approximate the conditional expectations $\mathbb{E}[ \psi^{2}(V_{\tau_{n}}) \vert X_{\tau_{n}} = x_{i} ]$ by \eqref{eq:DiscreteCondExpec};
\item[(b)] Define $\sigma_{SLV}(\cdot,\tau_{n})$ on the grid in the $x$-direction by formula \eqref{eq:FullyDiscreteLeverage};
\item[(c)] update $P_{n}$ by performing a numerical time step for \eqref{eq:EffectiveSemidiscrSystem} from $\tau_{n-1}$ to $\tau_{n}$;
\end{itemize}
end
\end{itemize}
end  
\newline \newline
}
Whenever a time step from $\tau_{n-1}$ to $\tau_{n}$ with the HV scheme is replaced by two half-time steps of the implicit Euler scheme, the inner iteration above is first performed for the substep from $\tau_{n-1}$ to $\tau_{n-1/2} = \tau_{n-1}+\Delta \tau /2$, yielding an approximation of $P(\tau_{n}-\Delta \tau/2)$ and $\sigma_{SLV}(\cdot,\tau_{n}-\Delta \tau/2)$. Next, the inner iteration is performed for the substep from  $\tau_{n}-\Delta \tau/2$ to $\tau_{n}$, yielding an approximation of $P(\tau_{n})$ and $\sigma_{SLV}(\cdot,\tau_{n})$.
Upon completion of the time stepping and iteration procedure above, the original approximation for $\sigma_{SLV}(\cdot,0)$ is replaced on the grid in the $x$-direction by $\sigma_{SLV}(\cdot,\tau_{1})$.
This appears more realistic as the original approximation was actually only valid for the index $i=i_{0}$.

\setcounter{equation}{0}
\section{Numerical experiments}\label{Experiments}

In this section, the effectiveness of the calibration procedure is illustrated by applying it to a practical example.
Here, we opt to consider the popular and challenging Heston-based SLV model, i.e.\ SLV model \eqref{eq:SLVmodel} with $\psi(v)=\sqrt{v}$ and 
$\alpha=1/2$, to describe the evolution of the EUR/USD exchange rate.

As stated in the introduction, in financial practice it is common to first determine the SV parameters of the underlying SV model and to define the LV function such that the LV model reproduces the known market prices for European call and put options.
Afterwards, the calibration procedure aims at matching the SLV model with its underlying LV model, i.e.\ at obtaining equality \eqref{eq:SLVmatchesLV}.
In this article, we assume that the SV parameters and the LV function are known and we then apply the calibration procedure from Section \ref{Calibration}. 
The performance is illustrated by comparing the densities from the LV model and from the SLV model and by comparing the corresponding option values. 

\begin{figure}
\begin{center}
\includegraphics[scale=0.5]{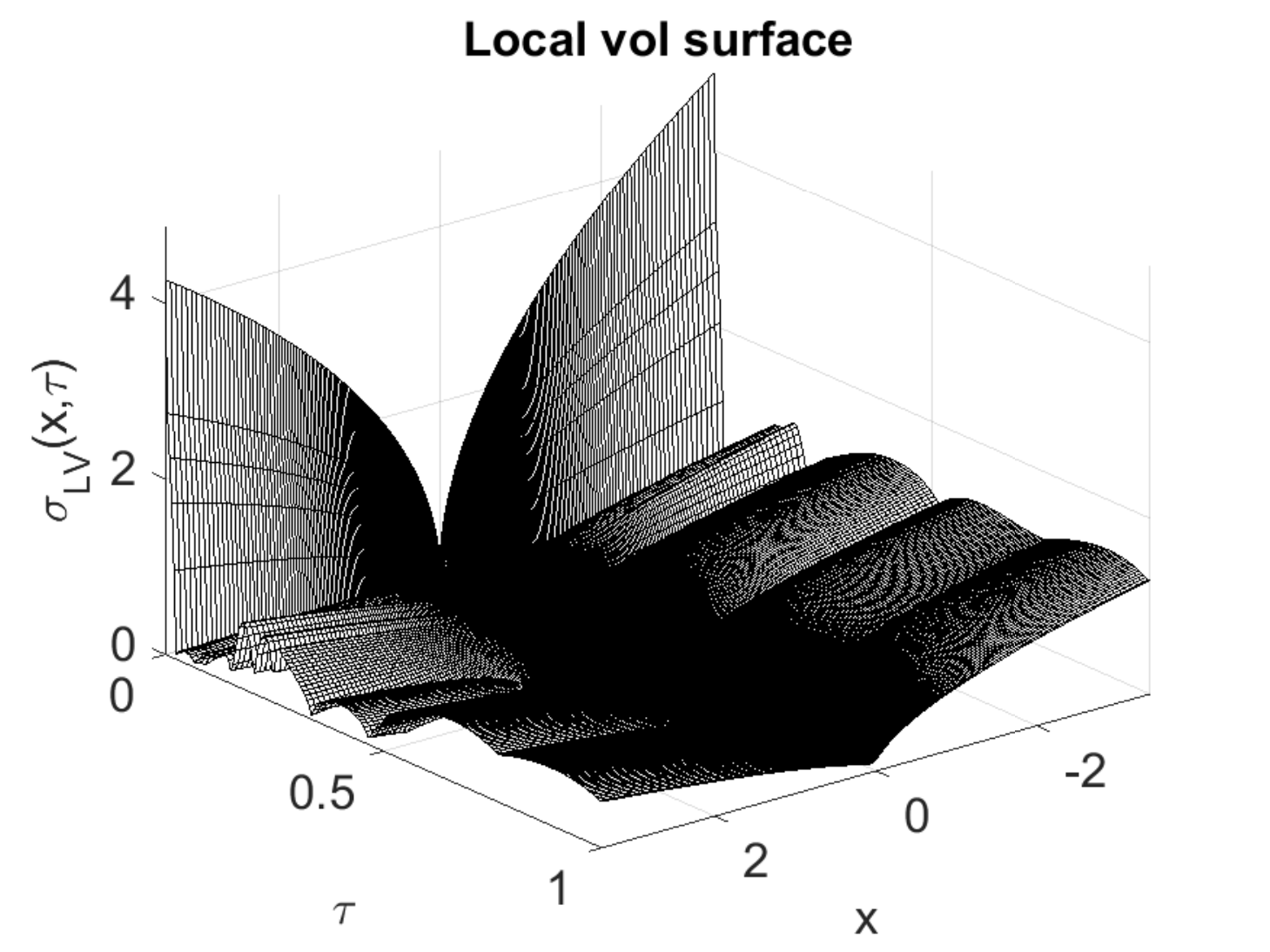} 
\includegraphics[scale=0.5]{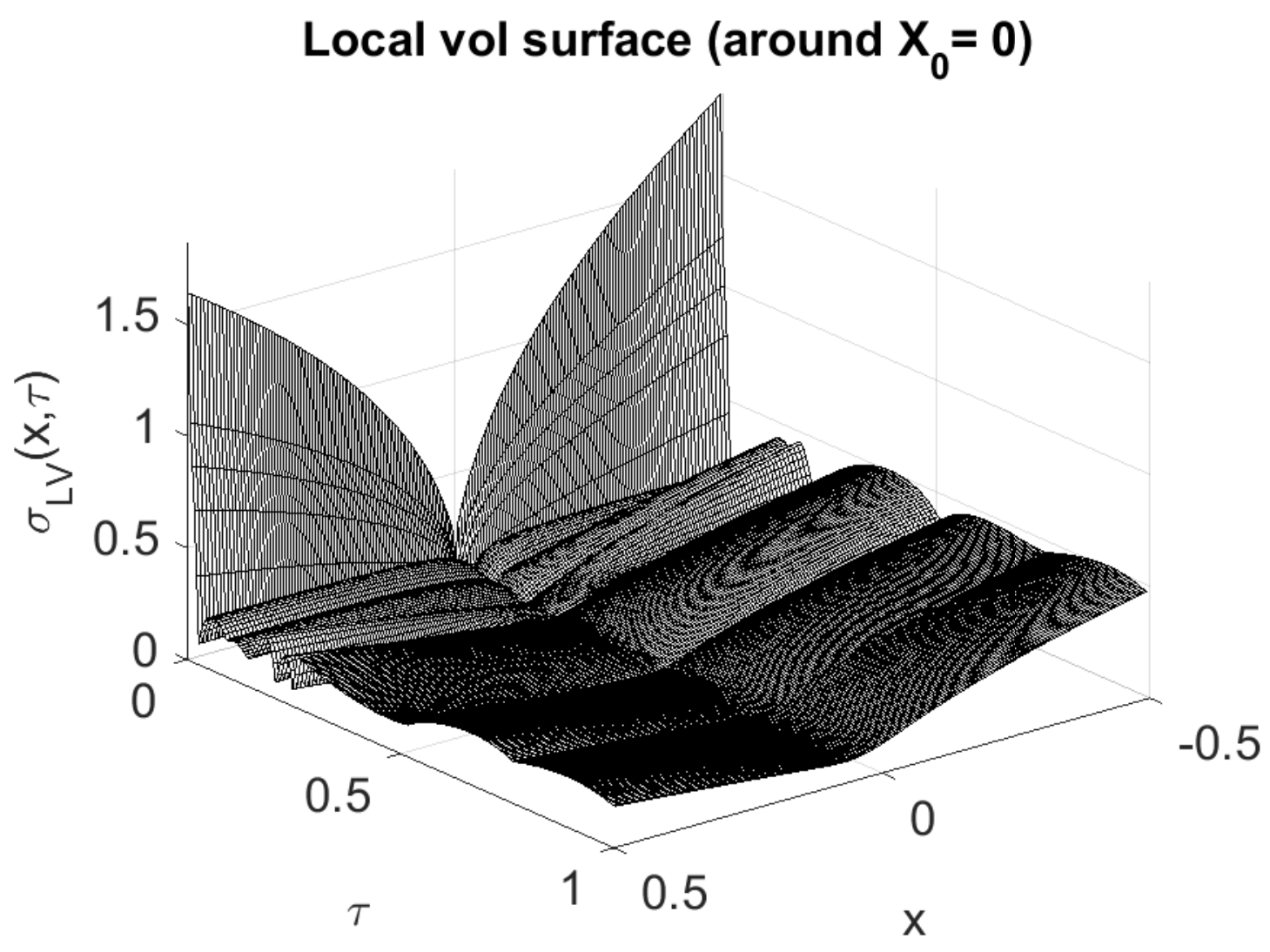} 
\caption{Local volatility function originating from actual EUR/USD vanilla option data (market data as of 2 March 2016).
The spot rate $S_{0} = 1.08815$.}
\label{fig:LVSurface}
\end{center}
\end{figure}

\begin{figure}
\begin{center}
\includegraphics[scale=0.5]{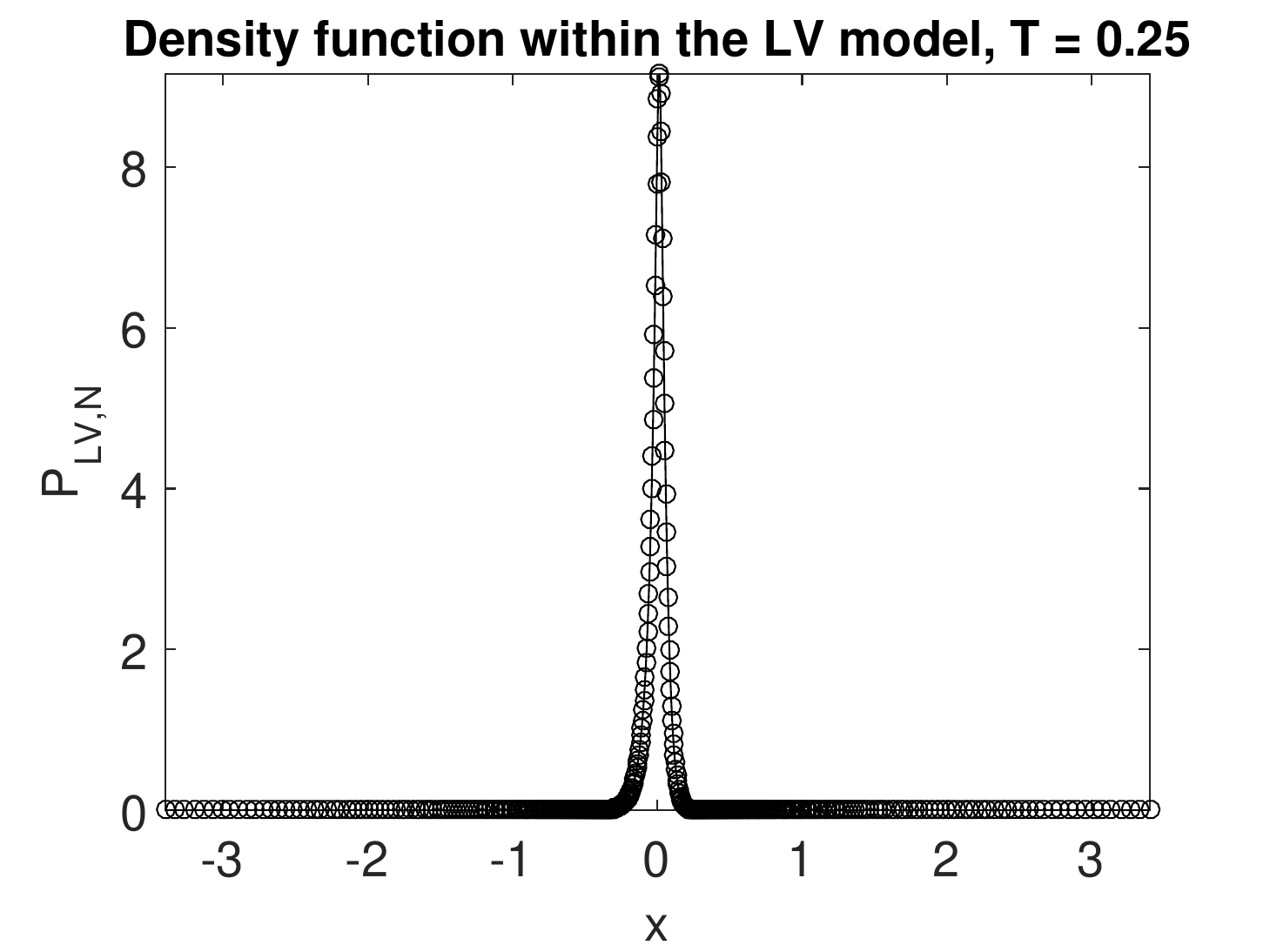} 
\includegraphics[scale=0.5]{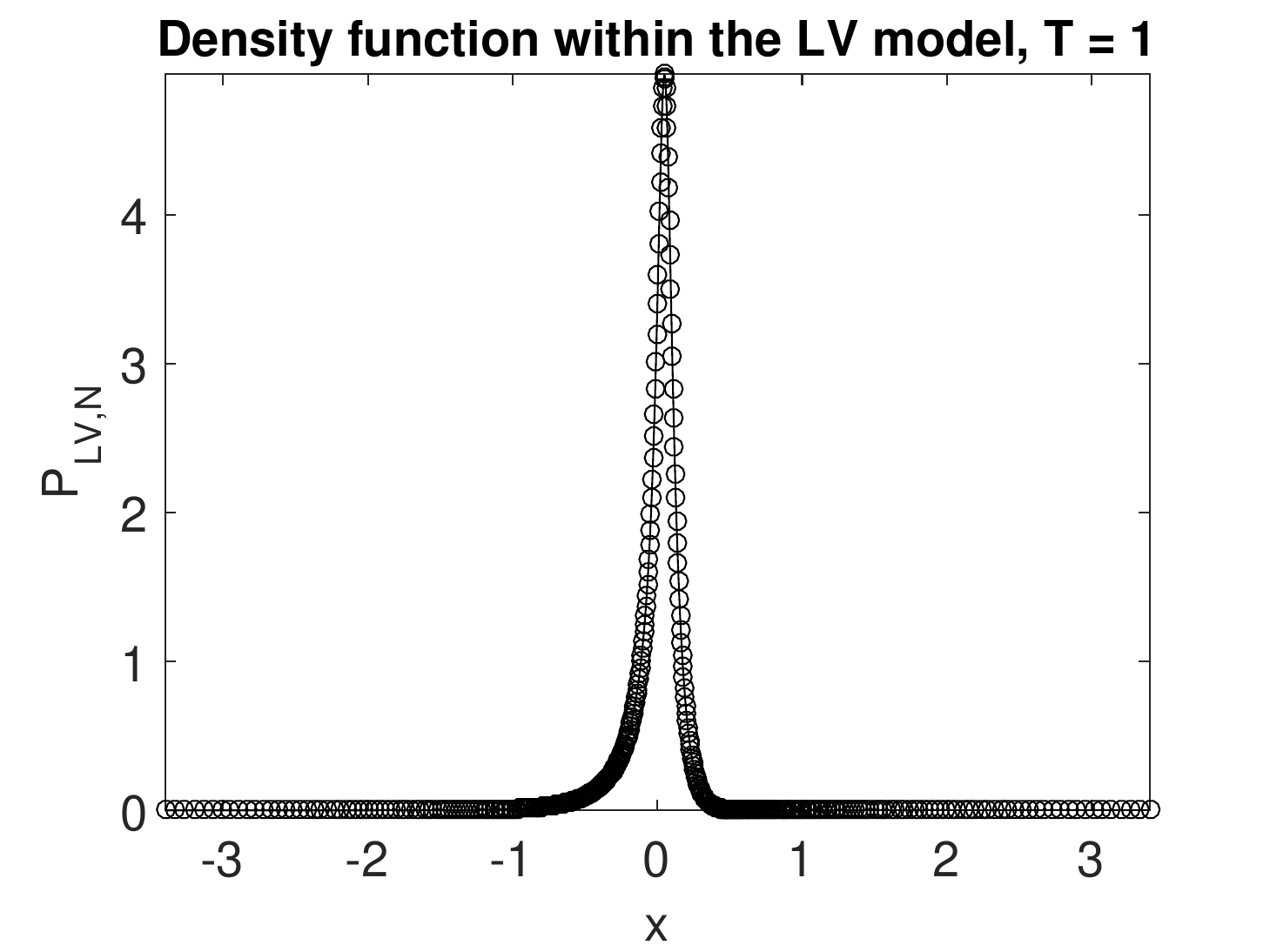}  
\caption{Approximation of the density function $p_{LV}(x,0.25)$ (left) and $p_{LV}(x,1)$ (right) by applying the FV discretization from Subsection \ref{1DKolmogorov} with $m=400$.}
\label{fig:LVDensity}
\end{center}
\end{figure}

For the actual experiments we consider the following sets of SV parameters:
\begin{center}
\begin{tabular}{|c| r r r r r r|}
\hline \newline
& $\kappa$ & $\eta$ & $\xi$ & $\rho$ & $T$ & $V_{0}$ \\
\hline \newline
Set E & $5$ & $ 0.16 $ & $ 0.9 $ & $ 0.1 $ &  $ 0.25 $ & $0.0625$ \\
\hline
\newline
Set F & $ 1.15 $ & $ 0.0348 $ & $ 0.39 $ & $ -0.64 $ &  $ 0.25 $ & $0.0348$ \\
\hline
\newline
Set G & $ 1.50 $ & $ 0.0154 $ & $ 0.24 $ & $ -0.11 $ &  $ 1$ & $0.0154$ \\
\hline
\end{tabular}
\end{center}
The Sets E and F correspond with the SV parameters from Sets C and D, and are taken from \cite{FO11}. 
Set G is taken from \cite{C11} and corresponds to the EUR/USD exchange rate for the pertinent maturity (market data as of 16 September 2008).
For Set E it holds that $q=\tfrac{2\kappa\eta}{\xi^2}-1=0.98$ and the process $V_{\tau}$ is strictly positive. For Set F, respectively Set G, it holds that $q=-0.47$, respectively $q=-0.20$, such that $V_{\tau}=0$ is attainable.
The LV model is completely determined by the LV function, the risk-free interest rates and the spot value $S_{0}$.
We assume that the risk-free interest rates are given by
$$ r_{d} = 0.02, \qquad r_{f} = 0.01, $$
and that the LV function is as displayed in Figure \ref{fig:LVSurface}. The pertinent LV function originates from actual EUR/USD vanilla option data (market data as of 2 March 2016) and is constructed by using an SSVI-type method for implied volatility interpolation, \cite{GJ14}.
The corresponding spot rate is given by
$$ S_{0} = 1.08815. $$
Note that, even if the LV function is given, there is often no analytical expression available for the density function $p_{LV}$ or the option values. It is well-known, however, that the density function satisfies the 1D forward Kolmogorov equation
\begin{equation*}
\tfrac{\partial}{\partial \tau} p_{LV} = \tfrac{\partial}{\partial x^{2}} \left( \tfrac{1}{2} \sigma^{2}_{LV} p_{LV} \right) - \tfrac{\partial}{\partial x} \left( (r_{d} - r_{f} - \tfrac{1}{2} \sigma^{2}_{LV})p_{LV} \right),
\end{equation*}
for $x \in \mathbb{R}, \tau >0$.
By applying the FV discretization described in Subsection \ref{1DKolmogorov} one defines approximations $P_{LV,i}(\tau)$ of the exact density values $p_{LV}(x_{i},\tau)$. 
Fully discrete approximations $P_{LV,N,i}$ of $p_{LV}(x_{i},T)$ are then obtained by applying a suitable time stepping method.
In Figure \ref{fig:LVDensity} the latter approximations are shown for $\tau=0.25$, respectively $\tau=1$, and the practical value $m=400$.

\begin{figure}
\begin{center}
\includegraphics[scale=0.5]{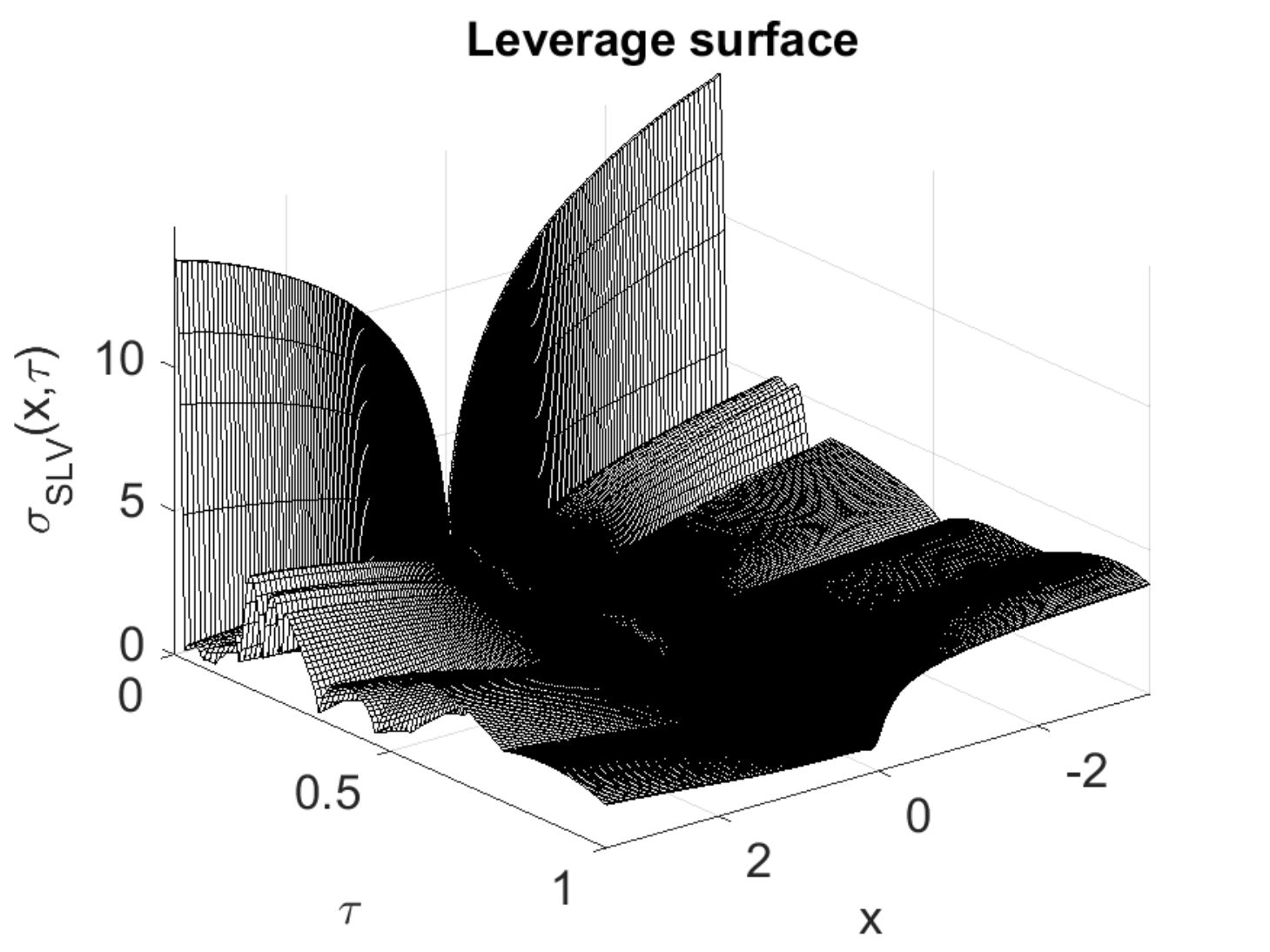} 
\includegraphics[scale=0.5]{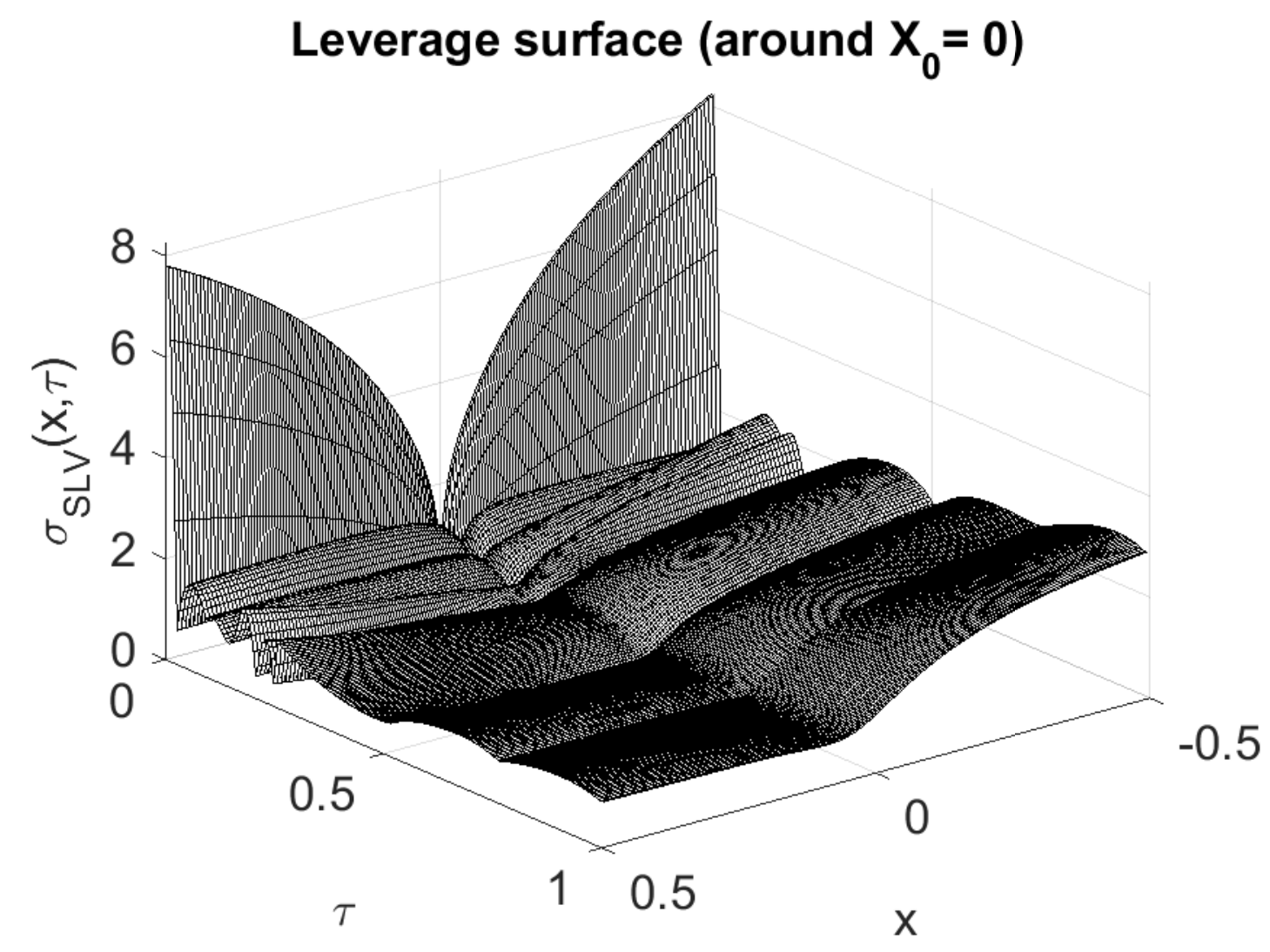} 
\caption{Leverage function stemming from the calibration procedure with local volatility function from Figure \ref{fig:LVSurface}, SV parameters from Set G and values $m_{1} = 400$, $m_{2} = 200$, $\Delta \tau = 1/200$, $\theta = \tfrac{1}{2}+\tfrac{1}{6}\sqrt{3}$, $Q =2.$}
\label{fig:SLVSurface_SetG}
\end{center}
\end{figure}

\begin{figure}
\begin{center}
\includegraphics[scale=0.5]{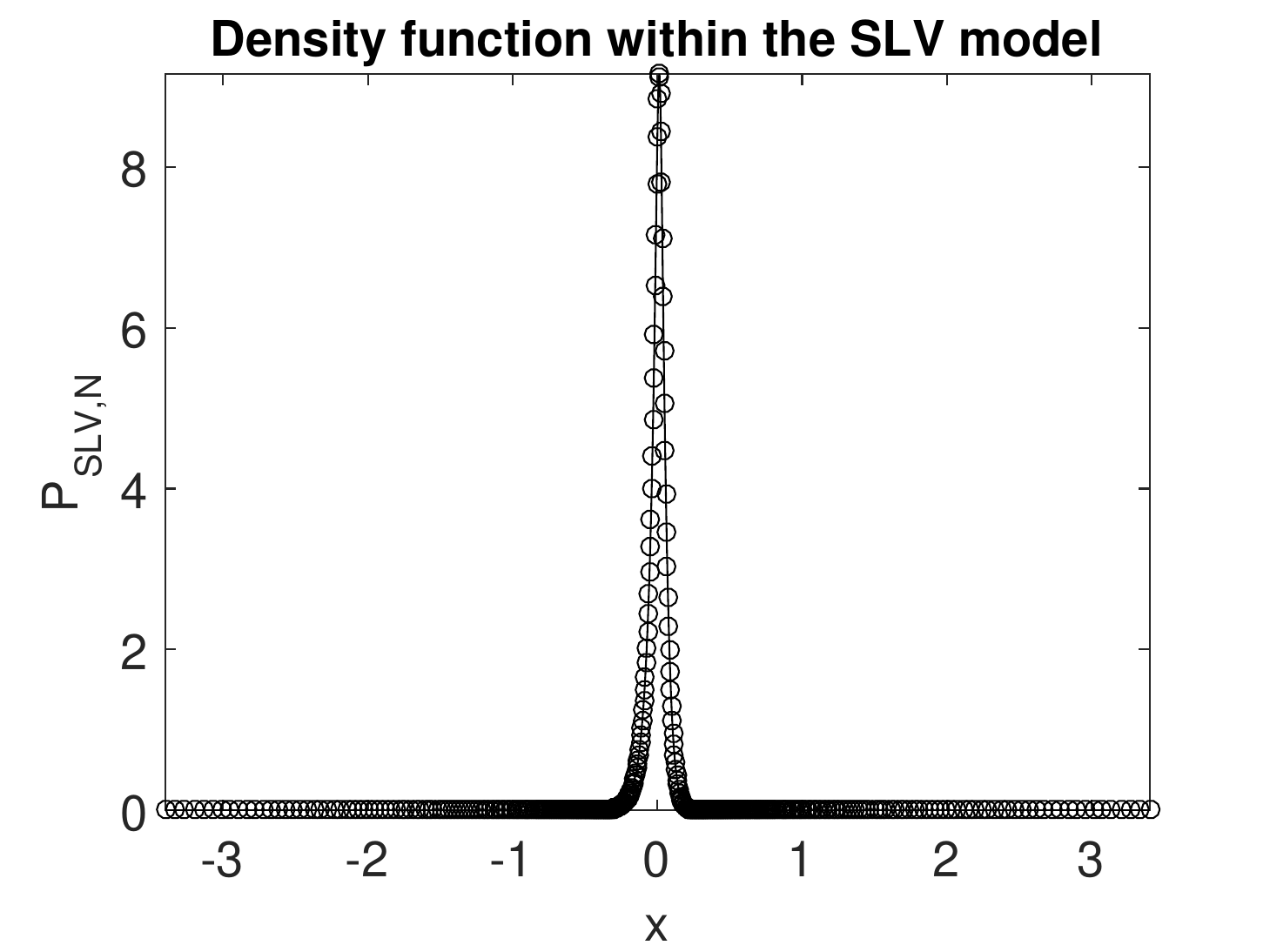}  
\includegraphics[scale=0.5]{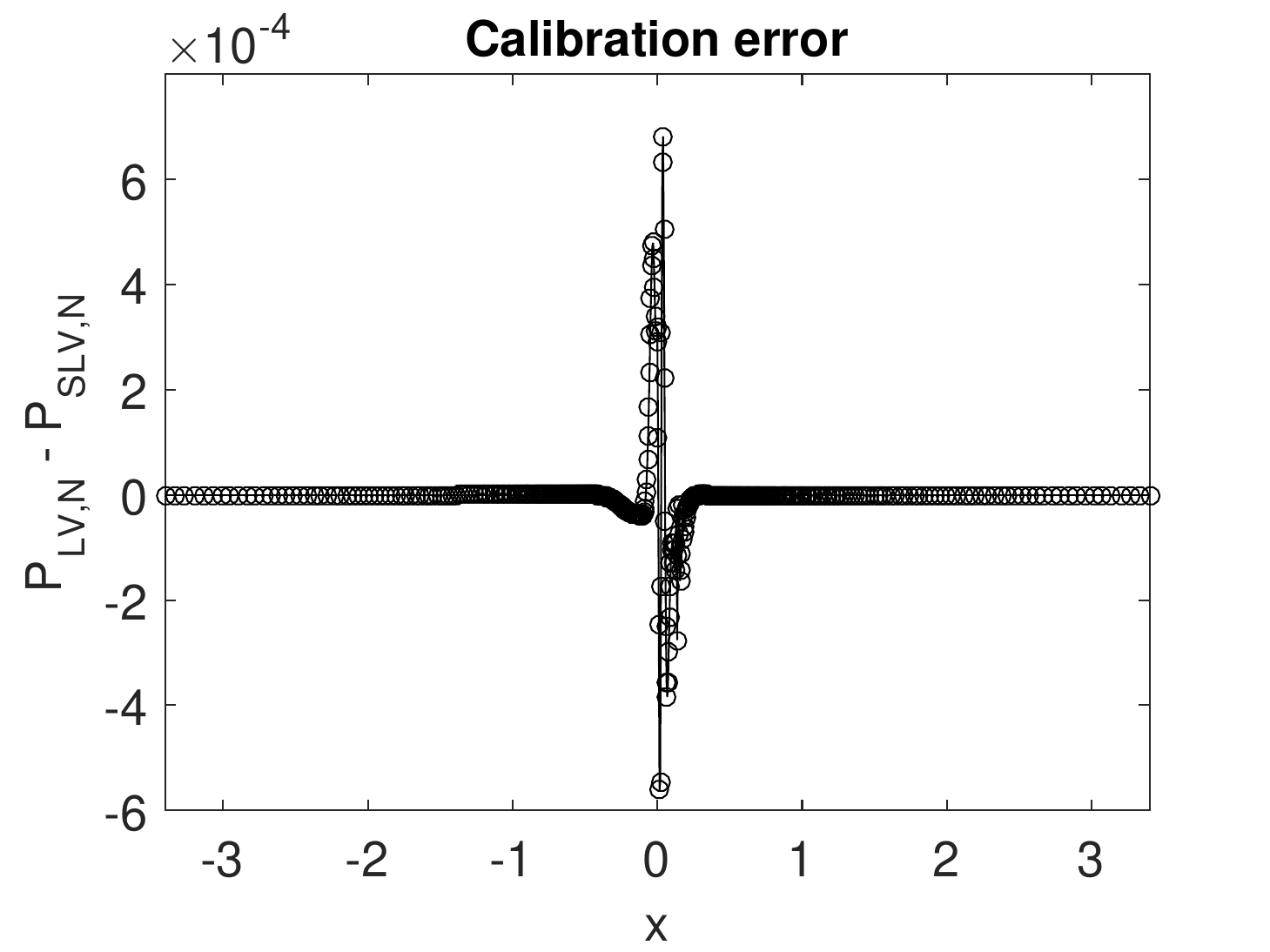} \\
Illustration of the density in case of Set E. \newline \newline \newline
\includegraphics[scale=0.5]{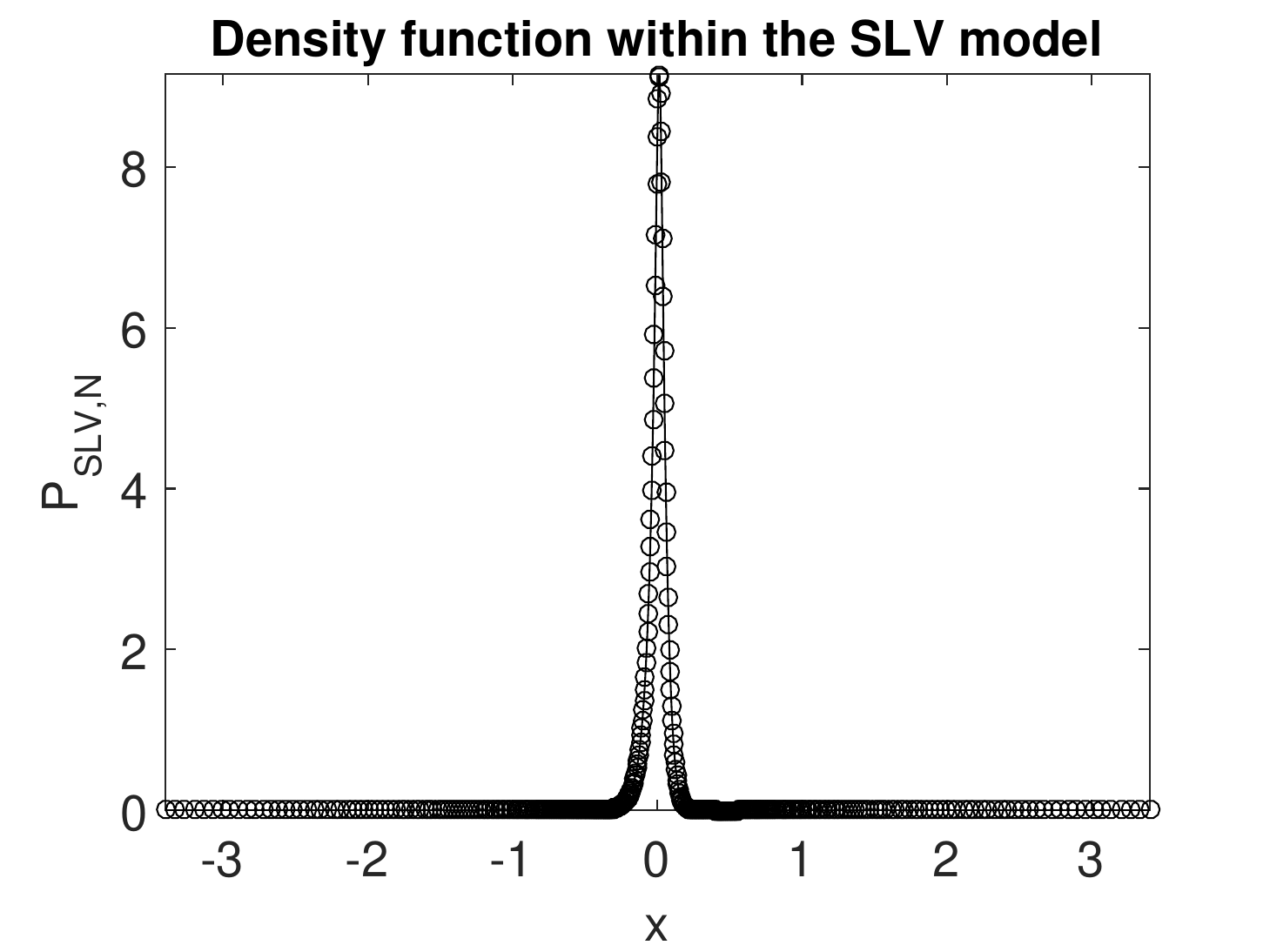}  
\includegraphics[scale=0.5]{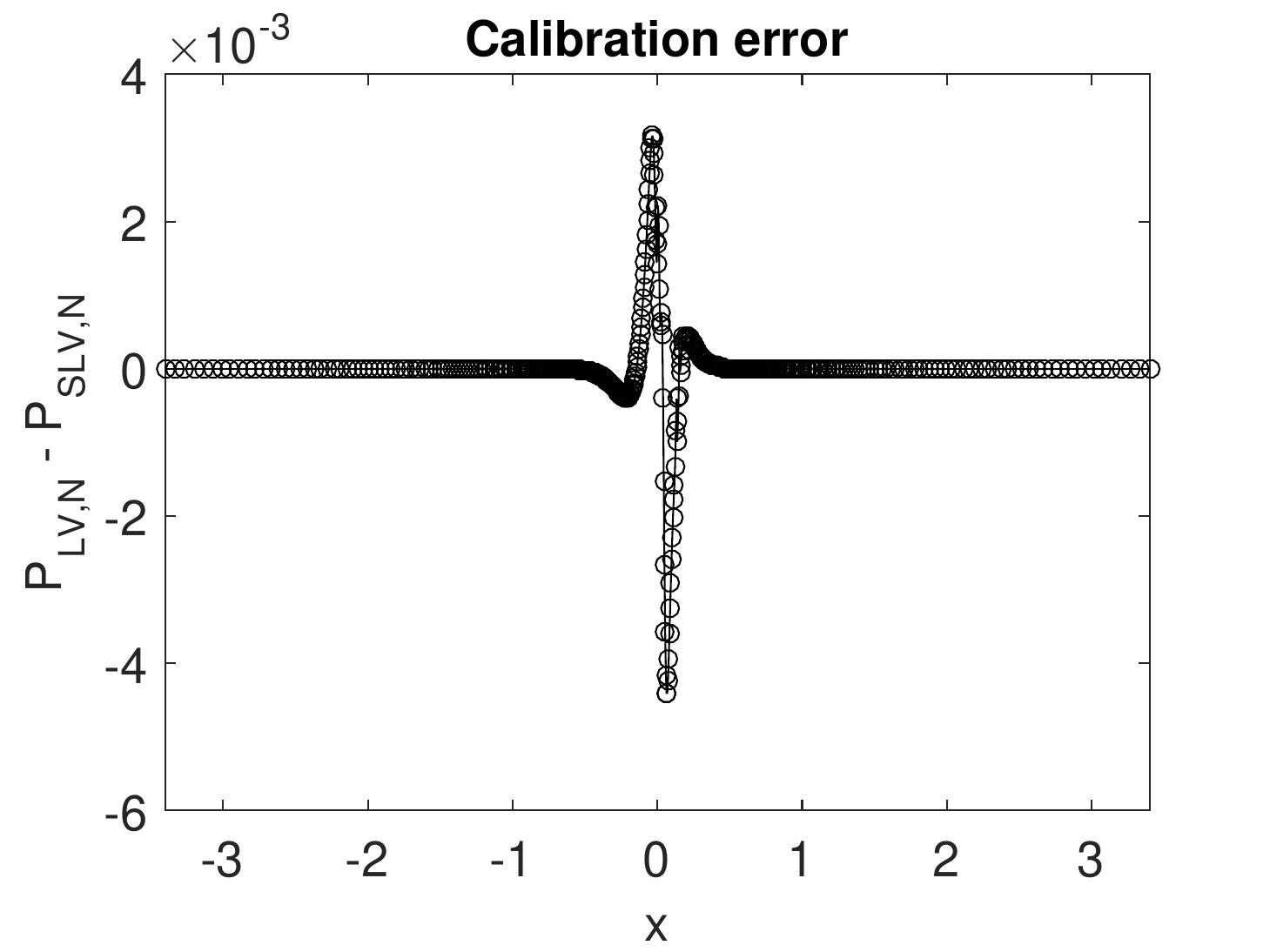} \\
Illustration of the density in case of Set F. \newline \newline \newline
\includegraphics[scale=0.5]{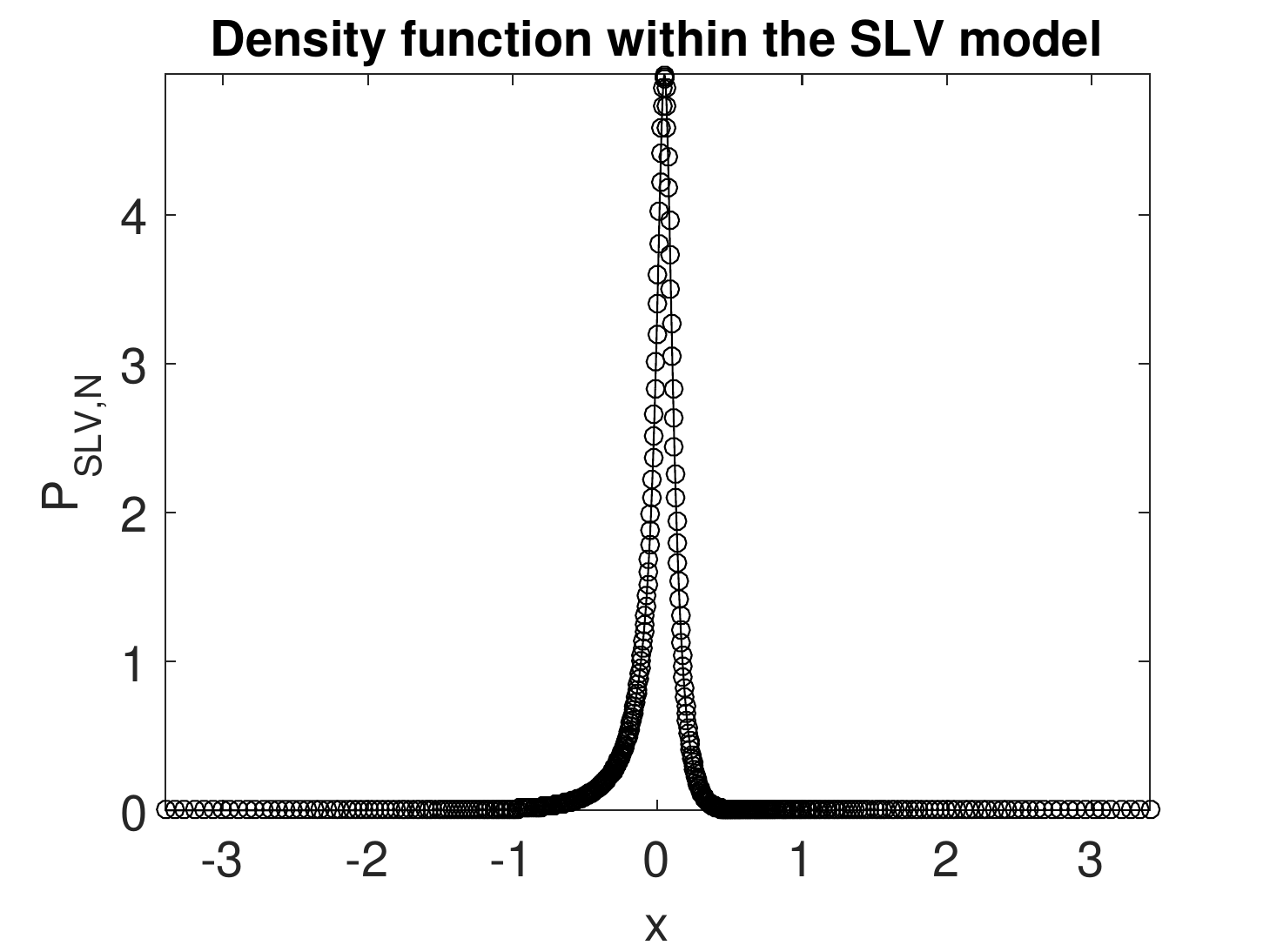}  
\includegraphics[scale=0.5]{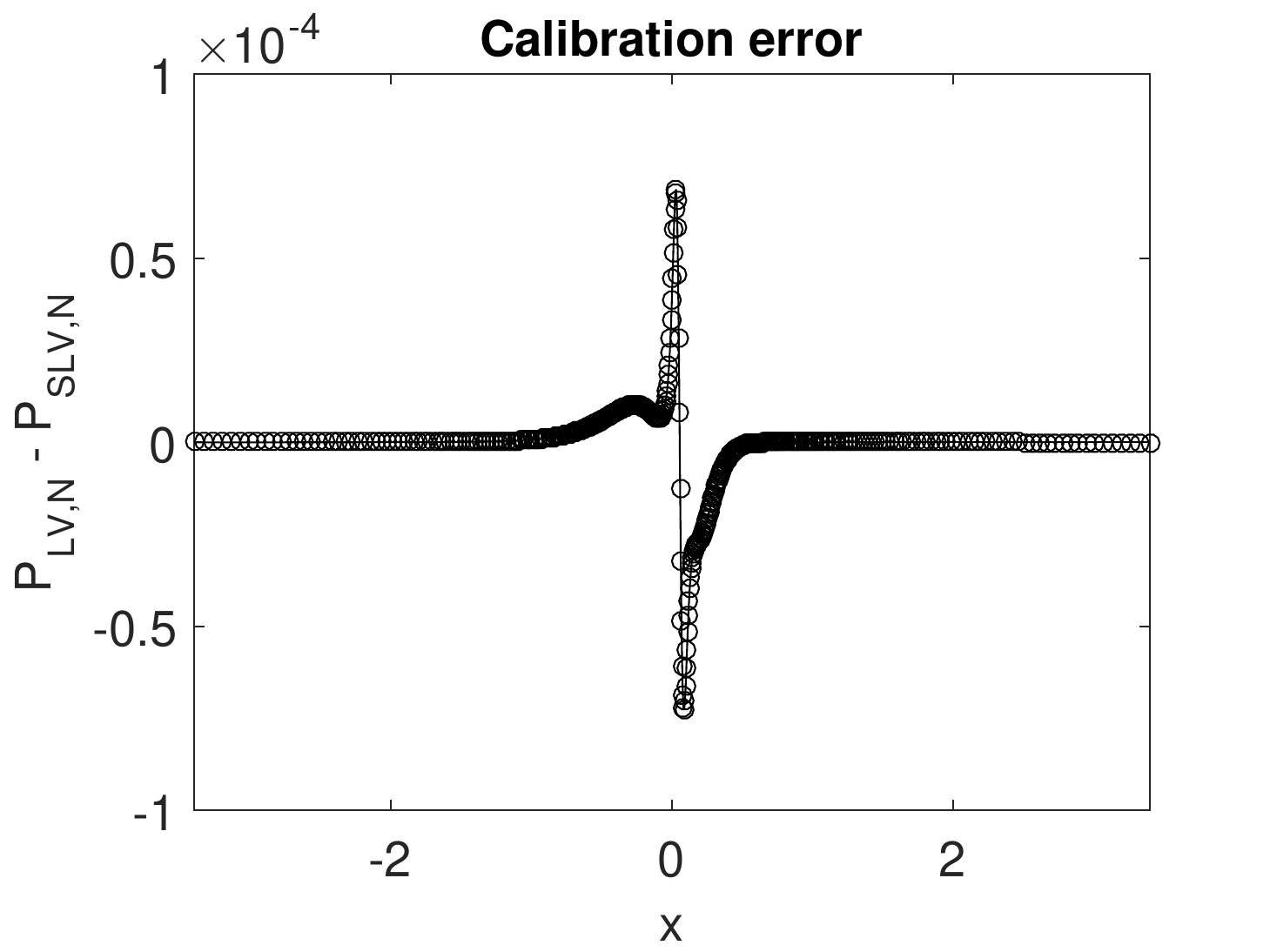} \\
Illustration of the density in case of Set G. \newline
\caption{Comparison of the fully discrete density functions $P_{LV,N}$ and $P_{SLV,N}$ for each of the parameters sets and for values $m_{1} = 400$, $m_{2} = 200$, $\Delta \tau = 1/200$, $\theta = \tfrac{1}{2}+\tfrac{1}{6}\sqrt{3}$, $Q=2$. }
\label{fig:SLVDensities}
\end{center}
\end{figure}

Once the underlying SV model and LV model are specified, the calibration procedure from Section \ref{Calibration} can be applied. For the actual experiments we consider the discretization parameters
$$ m_{1} = 400, \quad m_{2} = 200, \quad \Delta \tau = 1/200, \quad \theta = \tfrac{1}{2}+\tfrac{1}{6}\sqrt{3}, \quad Q=2. $$
In Figure \ref{fig:SLVSurface_SetG} the resulting discrete leverage function is shown for set G.
In order to illustrate the performance of the calibration, we first consider a discrete version of \eqref{eq:SLVmatchesLV}.
More precisely, the discrete numerical densities $P_{LV,N,i}$ from the LV model are compared with 
\begin{equation*}
P_{SLV,N,i} := \sum_{j=1}^{m_{2}} \boldsymbol{P}_{N,i,j} \tfrac{\Delta v_{j} + \Delta v_{j+1}}{2},
\end{equation*}
which can be seen as the fully discrete approximations of 
$$ \int_{0}^{\infty} p(x_{i},v,T)dv, $$
within the SLV model after applying the trapezoidal rule for numerical integration.
In the left plots of Figure \ref{fig:SLVDensities} the approximations $P_{SLV,N}$ are shown for each of the three sets of parameters. In the right plots, the corresponding differences $P_{LV,N}-P_{SLV,N}$ are plotted. 
Note that the final time $T=0.25$ for Set E and Set F, and $T=1$ for Set G. 
From Figure \ref{fig:SLVDensities} it is readily seen that the difference between the fully discrete numerical densities is very small and hence that the calibration procedure performs well.

The main goal of the calibration procedure is to define the leverage function in such a way that the LV model and the SLV model define the same fair values for non-path-dependent European options.
If the leverage function is defined by \eqref{eq:SigmatoMatch}, then it follows for the fair value (FaV) of such an option with payoff $u_{0}(x), x \in \R$, at maturity $T$ that
\begin{equation*}
\mathrm{FaV} = e^{-r_{d}T} \int_{-\infty}^{\infty} p_{LV}(x,T) u_{0}(x) dx = e^{-r_{d}T} \int_{0}^{\infty} \int_{-\infty}^{\infty} p(x,v,T) u_{0}(x) dx dv. 
\end{equation*}
Given the approximations $P_{LV,N}$ and $\boldsymbol{P}_{N}$, the pertinent FaV can easily be approximated by applying numerical integration with the trapezoidal rule. In case of the SLV model, it is readily seen that defining the approximated FaV via $\boldsymbol{P}_{N}$ and the trapezoidal rule is equivalent with defining the FaV via $P_{SLV,N}$ and the trapezoidal rule.
Denote by $\mathrm{FaV}_{LV}$, respectively $\mathrm{FaV}_{SLV}$, the approximated fair values obtained via $P_{LV,N}$, respectively $P_{SLV,N}$. We now compare these approximations for a set of European call options with a range of strikes given by
$$ K = 0.75 S_{0}, \ 0.8 S_{0}, \ 0.9 S_{0}, \ S_{0}, \ 1.1 S_{0}, \ 1.2 S_{0}, \ 1.25 S_{0}. $$
When the strike increases relatively to $S_{0}$, the fair value of European call options tends to zero and it is difficult to adequately compare approximations. In financial practice, European call and put options are often quoted in terms of \textit{implied volatility}. Let $\sigma_{imp,LV}$, respectively $\sigma_{imp,SLV}$, denote the implied volatility (in \%) corresponding to $\mathrm{FaV}_{LV}$, respectively $\mathrm{FaV}_{SLV}$. In the following we test the performance of the calibration procedure by calculating the absolute implied volatility errors
\begin{equation*}
\epsilon_{imp} = \vert \sigma_{imp,LV} - \sigma_{imp,SLV} \vert.
\end{equation*}
In Table \ref{Table:ImpliedVolErrors} these errors are presented for the different SV parameter sets, taking the same values of $m,m_{1}, m_{2}, \Delta \tau, \theta, Q$ as above. 
\begin{table}
\begin{center}
\begin{tabular}{|c||c|c|c||c|c|}
\hline
& $T = 0.25$ & Set E & Set F & $T=1$ & Set G \\
\hline
$K/S_{0}$ & $\sigma_{imp,LV}$ & $\epsilon_{imp}$ & $\epsilon_{imp}$ & $\sigma_{imp,LV}$ & $\epsilon_{imp}$ \\
\hline
$0.75$ & 19.18 & 0.1005 & 0.1208 & 21.94 & 0.0021 \\ 
$0.80$ & 18.40 & 0.0212 & 0.0454 & 20.20 & 0.0015 \\
$0.90$ & 15.01 & 0.0033 & 0.0154 & 16.65 & 0.0008 \\ 
$1.0$  & 11.26 & 0.0011 & 0.0030 & 13.14 & 0.0004 \\
$1.10$ & 11.59 & 0.0011 & 0.0153 & 11.38 & 0.0003 \\
$1.20$ & 13.20 & 0.0009 & 0.0937 & 11.77 & 0.0003 \\
$1.25$ & 14.03 & 0.0006 & 0.1888 & 12.12 & 0.0003 \\
\hline
\end{tabular}
\end{center}
\caption{Comparison of the approximated implied volatilities $\sigma_{imp,LV}$ and $\sigma_{imp,SLV}$ for values $m_{1} = 400$, $m_{2} = 200$, $\Delta \tau = 1/200$, $\theta = \tfrac{1}{2}+\tfrac{1}{6}\sqrt{3}$, $Q =2.$}
\label{Table:ImpliedVolErrors}
\end{table}
The somewhat larger values $\epsilon_{imp}$ for $T=0.25$ compared to $T=1$ can be explained from the fact that the implied volatility is more sensitive to changes in the fair value when the maturity is low.
The results in Table \ref{Table:ImpliedVolErrors} confirm that the calibration procedure performs well. They indicate that the fully discrete leverage surface is, indeed, defined such that the SLV model reproduces accurately the known market prices for European call options.

\setcounter{equation}{0}
\section{Conclusion}\label{Conclusion}

Stochastic local volatility (SLV) models constitute state-of-the-art models to describe asset price processes. Their calibration to the underlying local volatility model is, however, highly non-trivial. It incorporates the solution of non-linear forward Kolmogorov equations. 
In general, no analytical solution is available and one relies on numerical methods in order to approximate the exact solution. Here, we introduce a finite volume - alternating direction implicit method for the numerical solution of general one-dimensional and two-dimensional forward Kolmogorov equations.
The finite volume spatial discretization does not require a transformation of the PDE, which constitutes a main advantage in the calibration of SLV models, and handles the boundary conditions in a natural way. Moreover, the finite volume scheme preserves the crucial property that total mass of a density function is always equal to one. Our numerical experiments for relevant practical applications confirm that the pertinent spatial discretization is convergent.
Temporal discretization is performed by using the Hundsdorfer--Verwer ADI scheme. By splitting the semidiscrete system into different parts that represent spatial derivatives in the different spatial dimensions, and by handling spatial derivatives in only one spatial dimension in the implicit substeps, a major computational advantage can be achieved.
The non-linearity in the calibration procedure of stochastic local volatility models is handled by introducing an inner iteration.
Our numerical experiments reveal that the proposed calibration procedure performs well. The calibrated stochastic local volatility model matches the underlying local volatility model almost exactly, both in terms of the density function and of the implied volatilities of European call options.

\setcounter{equation}{0}
\section*{Funding} 
The work by the first author has been supported financially by a PhD Fellowship of the Research Foundation--Flanders.


\begin{thebibliography}{99}

\bibitem{AP06} L.~B.~G. Andersen and V.~V. Piterbarg, 
\textit{Moment explosions in stochastic volatility models},
Finance Stoch. \textbf{11} (2007), pp.\ 29--50.

\bibitem{C11} I.~J. Clarke, 
\textit{Foreign Exchange Option Pricing: A Practitioner's Guide},
John Wiley \& Sons, 2011.

\bibitem{CIR85} J.~C. Cox, J.~E. Ingersoll and S.~A. Ross,
\textit{A theory of the term structure of interest rates},
Econometrica \textbf{53} (1985), pp.\ 385--407.

\bibitem{D94} B. Dupire, 
\textit{Pricing with a smile}, 
RISK, January (1994), pp.\ 18--20.

\bibitem{EKO12} B. Engelmann, F. Koster and D. Oeltz,
\textit{Calibration of the Heston stochastic local volatility model: a finite volume scheme}, 
Available at SSRN 1823769 (2012).

\bibitem{FO08} F. Fang and C.~W. Oosterlee, 
\textit{A novel pricing method for European options based on Fourier-cosine series expansions},
SIAM J. Sci. Comput. \textbf{31} (2008), pp.\ 826--848.

\bibitem{FO11} F. Fang and C.~W. Oosterlee, 
\textit{A Fourier-based valuation method for Bermudan and barrier option under Heston's model},
SIAM J. Financial Math. \textbf{2} (2011), pp.\ 439--463.

\bibitem{GJ14} J. Gatheral and A. Jacquier, 
\textit{Arbitrage-free SVI volatility surfaces},
Quant. Financ. \textbf{14} (2014), pp.\ 59--71.

\bibitem{G86} I. Gy\"ongy, 
\textit{Mimicking the one-dimensional marginal distributions of processes having an Ito differential},
Probab. Th. Rel. Fields \textbf{71} (1986), pp.\ 501--516.

\bibitem{HH12} T. Haentjens and K.~J. in~'t~Hout, 
\textit{Alternating direction implicit finite difference schemes for the Heston--Hull--White partial differential equation},
J. Comp. Finan. \textbf{16} (2012), pp.\ 83--110.

\bibitem{H09} P. Henry-Labord\`ere, 
\textit{Calibration of local stochastic volatility models to market smiles},
Risk, September (2009), pp.\ 112--117.

\bibitem{HY02} V.~E. Henson and U.~M. Yang, 
\textit{BoomerAMG: A parallel algebraic multigrid solver and preconditioner},
Appl. Numer. Math. \textbf{41} (2002), pp.\ 155--177.

\bibitem{H93} S.~L. Heston, 
\textit{A closed-form solution for options with stochastic volatility with applications to bond and currency options},
Rev. Finan. Stud. \textbf{6} (1993), pp.\ 327--343.

\bibitem{IHF10} K.~J. in~'t~Hout and S. Foulon, 
\textit{ADI finite difference schemes for option pricing in the Heston model with correlation},
Int. J. Numer. Anal. Mod. \textbf{7} (2010), pp.\ 303--320.

\bibitem{IHM13} K.~J. in~'t~Hout and C. Mishra, 
\textit{Stability of ADI schemes for multidimensional diffusion equations with mixed derivative terms},
Appl. Numer. Math. \textbf{74} (2013), pp.\ 83--94.

\bibitem{IHW07} K.~J. in~'t~Hout and B.~D. Welfert, 
\textit{Stability of ADI schemes applied to convection-diffusion equations with mixed derivative terms},
Appl. Numer. Math. \textbf{57} (2007), pp.\ 19--35.

\bibitem{IHW09} K.~J. in~'t~Hout and B.~D. Welfert, 
\textit{Unconditional stability of second-order ADI schemes applied to multi-dimensional diffusion equations with mixed derivative terms},
Appl. Numer. Math. \textbf{59} (2009), pp.\ 677--692.

\bibitem{IHW15} K.~J. in~'t~Hout and M. Wyns, 
\textit{Convergence of the Hundsdorfer--Verwer scheme for two-dimensional convection-diffusion equations with mixed derivative term},
AIP Conf. Proc. \textbf{1648} (2015), 850054.

\bibitem{H02} W. Hundsdorfer, 
\textit{Accuracy and stability of splitting with Stabilizing Corrections},
Appl. Numer. Math. \textbf{42} (2003), pp.\ 213--233.

\bibitem{HV03} W. Hundsdorfer and J.~G. Verwer, 
\textit{Numerical Solution of Time-Dependent Advection-Diffusion-Reaction Equations},
Springer, Berlin, 2003.

\bibitem{L02} A. Lipton, 
\textit{The vol smile problem},
RISK, February (2002), pp.\ 61--65.

\bibitem{PVF03} D.~M. Pooley, K.~R. Vetzal and P.~A. Forsyth, 
\textit{Convergence remedies for non-smooth payoffs in option pricing},
J. Comp. Finan. \textbf{6} (2003), pp.\ 25--40.

\bibitem{R84} R. Rannacher,
\textit{Finite element solution of diffusion problems with irregular data},
Numer. Math. \textbf{43} (1984), pp.\ 309--327.

\bibitem{RMQ07} Y. Ren, D. Madan and M.~Q. Qian,
\textit{Calibrating and pricing with embedded local volatility models},
Risk, September (2007), pp.\ 138--143.

\bibitem{R89} H. Risken, 
\textit{The Fokker-Planck Equation: Methods of Solution and Applications,} 2nd edition,
Springer, Berlin, 1989.

\bibitem{RO12} M.~J. Ruijter and C.~W. Oosterlee,
\textit{Two-dimensional Fourier cosine series expansion method for pricing financial options},
SIAM J. Sci. Comput. \textbf{34} (2012), pp.\ B642--B671.

\bibitem{T11} R. Tachet, 
\textit{Non-Parametric Model Calibration in Finance},
PhD thesis, Ecole Centrale Paris, 2011.

\bibitem{TF10} G. Tataru and T. Fisher, 
\textit{Stochastic local volatility},
Technical report, Bloomberg (2010).

\bibitem{TR00} D. Tavella and C. Randall,
\textit{Pricing Financial Instruments: The Finite Difference Method},
John Wiley \& Sons, New York, 2000.

\bibitem{VGO14} A.~W. van der Stoep, L.~A. Grzelak and C.~W. Oosterlee,
\textit{The Heston stochastic-local volatility model: efficient Monte Carlo simulation},
Int. J. Theor. Appl. Finan. \textbf{17} (2014), pp.\ 1450045-1--1450045-30.

\bibitem{V99} J.~G. Verwer, E.~J. Spee, J.~G. Blom and W. Hundsdorfer,
\textit{A second-order Rosenbrock method applied to photochemical dispersion problems},
SIAM J. Sci. Comput. \textbf{20} (1999), pp.\ 1456--1480.

\bibitem{W16} M. Wyns,
\textit{Convergence analysis of the Modified Craig--Sneyd scheme for two-dimensional convection-diffusion equations with nonsmooth initial data},
To appear in IMA J. Numer. Anal. (2016), \url{http://dx.doi.org/10.1093/imanum/drw028}.

\end{thebibliography}
\end{document}